\documentclass[twoside]{amsart}
\usepackage{mathrsfs}
\usepackage{amsfonts}
\usepackage{amssymb}
\usepackage{amsmath}
\usepackage{amsthm}
\usepackage{mathdots}
\usepackage{todonotes}
\usepackage[all]{xy}

\usepackage[verbose,colorlinks,pagebackref]{hyperref}
\hypersetup{
  colorlinks   = true, 
  urlcolor     = blue, 
  linkcolor    = blue, 
  citecolor   = brown 
}

\usepackage[alphabetic,nobysame]{amsrefs}
\usepackage{lastpage}

\RequirePackage{amsmath} \RequirePackage{amssymb}
\usepackage{amscd,latexsym,amsthm,amsfonts,amssymb,amsmath,amsxtra}

\newcommand{\bpm}{\begin{pmatrix}}
\newcommand{\epm}{\end{pmatrix}}
\newcommand{\bsm}{\begin{smallmatrix}}
\newcommand{\esm}{\end{smallmatrix}}
\newcommand{\bspm}{\left(\begin{smallmatrix}}
\newcommand{\espm}{\end{smallmatrix}\right)}
\newcommand{\bbm}{\begin{bmatrix}}
\newcommand{\ebm}{\end{bmatrix}}
\usepackage{float}

    \newcommand{\BA}{{\mathbb {A}}} 
    \newcommand{\C}{{\mathbb {C}}}

     \newcommand{\Z}{{\mathbb {Z}}}

    \newcommand{\CA}{{\mathcal {A}}} \newcommand{\CB}{{\mathcal {B}}}
    \newcommand{\CC}{{\mathcal {C}}}

    \newcommand{\CM}{{\mathcal {M}}} 
    \newcommand{\CO}{{\mathcal {O}}} 
     
    \newcommand{\CS}{{\mathcal {S}}} 
     
    \newcommand{\CW}{{\mathcal {W}}}

       \newcommand{\fp}{{\mathfrak{p}}}

      \newcommand{\RB}{{\mathrm {B}}}

    \newcommand{\RI}{{\mathrm {I}}}

    \newcommand{\RU}{{\mathrm {U}}}

    \newcommand{\lenth}{{\mathrm {\lenth}}}

     \newcommand{\GL}{{\mathrm{GL}}}
    \newcommand{\Hom}{{\mathrm{Hom}}} \renewcommand{\Im}{{\mathrm{Im}}}
    \newcommand{\Ind}{{\mathrm{Ind}}}

      \newcommand{\Tr}{{\mathrm{Tr}}}

    \renewcommand{\Re}{{\mathrm{Re}}} 
    \newcommand{\Res}{{\mathrm{Res}}}

\newcommand{\Supp}{\mathrm{Supp}}

 \newcommand{\SO}{{\mathrm{SO}}}
 \newcommand{\Sym}{{\mathrm{Sym}}}

\newcommand{\vol}{{\mathrm{vol}}}

    \newcommand{\bx}{{\bf {x}}} \newcommand{\by}{{\bf {y}}} \newcommand{\bz}{{\bf {z}}} \newcommand{\bu}{{\bf {u}}} 
\newcommand{\bv}{{\bf {v}}}       \newcommand{\bt}{{\bf {t}}}  \newcommand{\Bs}{{\bf {s}}}\newcommand{\bW}{{\bf {W}}} 
\newcommand{\btau}{{\boldsymbol \tau}} 

\newcommand{\diag}{{\mathrm {diag}}}  \newcommand{\Mat}{{\mathrm {Mat}}}  

\newcommand{\V}{{\mathrm{V}}}

    \newcommand{\wt}{\widetilde} \newcommand{\wh}{\widehat} 
    
    \newcommand{\pair}[1]{\langle {#1} \rangle}
    \newcommand{\wpair}[1]{\left\{{#1}\right\}}

    \newcommand{\ov}{\overline}

    \newcommand{\bs}{\backslash}
    
    \theoremstyle{plain}

       \newtheorem*{theorem*}{Theorem}

    \newtheorem{theorem}{Theorem}[section] \newtheorem{corollary}[theorem]{Corollary}
    \newtheorem{lemma}[theorem]{Lemma}  \newtheorem{proposition}[theorem]{Proposition}
    \newtheorem {conj}[theorem]{Conjecture}  \newtheorem {indhypo}[theorem]{Inductive Hypothesis}

  \theoremstyle{definition}

      \theoremstyle{remark}
      \newtheorem{remark}[theorem]{Remark}

    \numberwithin{equation}{section}

\usepackage[top=1in, bottom=1in, left=1.25in, right=1.25in]{geometry}

\title[Product of Rankin-Selberg convolutions]{Product of Rankin-Selberg convolutions and a new proof of Jacquet's local converse conjecture}
\author{Pan Yan}
\address{Department of Mathematics, The University of Arizona, Tucson, AZ 85721, USA}
\email{panyan@arizona.edu}

\author{Qing Zhang}
\address{School of Mathematics and Statistics, Huazhong University of Science and Technology, Wuhan, 430074, China}
\email{qingzh@hust.edu.cn}
\thanks{The corresponding author is Qing Zhang.}
\thanks{The first named author is partially supported by an AMS-Simons Travel Grant. The second named author is partially supported by NSFC grant 12371010.}
\date{\today}
\subjclass[2020]{11F70, 22E50}
\keywords{Rankin-Selberg convolution, $L$-functions, gamma factors, local converse theorem}

\begin{document}

\dedicatory{Dedicated to our advisor Professor Jim Cogdell on the occasion of his 70th birthday.}
\begin{abstract}
This article constructs a family of integrals representing the product of Rankin–Selberg $L$-functions for $\GL_l \times \GL_m$ and $\GL_l \times \GL_n$ in the regime $m+n<l$. When $n=0$, these integrals reduce (up to a shift) to the classical Rankin-Selberg convolution integrals of Jacquet--Piatetski-Shapiro--Shalika (JPSS), thereby generalizing the JPSS integrals. We establish their basic properties and use them to define local gamma factors. As an application, we provide a new proof of Jacquet’s local converse conjecture.
\end{abstract}
\maketitle
\setcounter{tocdepth}{1}
\tableofcontents

\section{Introduction}
Let $F$ be a global field and $\BA$ its ring of adeles. The theory of Rankin-Selberg integrals for $L$-functions of $\GL_n(\BA)\times \GL_m(\BA)$, developed by Jacquet, Piatetski-Shapiro, and Shalika in a series of papers (see e.g., \cites{JPSS, Jacquet-Shalika-EulerI, Jacquet-Shalika-EulerII}, and the surveys by Cogdell \cites{Cogdell:Fields, Cogdell:IAS}), plays an important role in the study of $L$-functions and automorphic representations. A key application is the proof of converse theorems, such as those for $\GL_3$ in \cites{JPSS-GL3I, JPSS-GL3II} and for general $\GL_n$ in \cites{CPS:converse I, CPS:converse II}. These converse theorems have become an important tool in establishing the Langlands functoriality, as illustrated, for instance, in \cites{Ramakrishnan2000, CKPSS-functorial classical, Kim-Shahidi, Kim-functorial, CKPSS-functorial, CPSS-functoriality-quasisplit, CFK2018}.

The first goal of this article is to construct a family of Rankin-Selberg integrals that generalizes those of Jacquet, Piatetski-Shapiro, and Shalika. Concretely, we present a family of Rankin-Selberg type integrals that represent, up to a shift, the product Rankin-Selberg $L$-function
$$L(s_1,\pi\times \tau_1)L(s_2,\wt \pi\times \wt \tau_2),$$
where $\pi$ is an irreducible cuspidal automorphic representation of $\GL_l(\BA)$, $\tau_1$ (resp. $\tau_2$) is an irreducible cuspidal representation of $\GL_m(\BA)$ (resp. $\GL_n(\BA)$).  Here $l$ is a positive integer and $m,n$ are non-negative integers with $m+n<l$. If $n=0$, our integrals reduce to those defined by Jacquet--Piatetski-Shapiro--Shalika (JPSS for abbreviation). In this sense, our integrals indeed generalize the JPSS Rankin-Selberg convolution integrals.

 To give more details, we introduce some notation. For an integer $j$ with $0\le j\le l-m-n-1$, we set $k=l-m-n-1-j$ and consider the embedding $\iota_j:\GL_{m+n}\to \GL_{l}$ given by 
 $$\bpm a&b\\ c&d \epm \mapsto \bpm I_j &&& &\\ &a&&b&\\ &&1&&\\ &c&&d&\\ &&&&I_{k} \epm$$
for $a\in \Mat_{m\times m}, b\in \Mat_{m\times n}, c\in \Mat_{n\times m}, d\in \Mat_{n\times n}.$ Given an irreducible cuspidal automorphic representation $\pi$ (resp. $\tau_1,\tau_2$) of $\GL_l(\BA)$ (resp. $\GL_m(\BA), \GL_n(\BA)$), we consider the integral
$$I_j(\phi,f_{\Bs})=\int_{\GL_{m+n}(F)\bs \GL_{m+n}(\BA)}\phi^{\psi}_{Y_j}(\iota_j(h))E(h,f_{\Bs})dh.$$
Here $\phi\in \pi$ is a cusp form, and $\phi^{\psi}_{Y_j}$ is a certain Fourier coefficient of $\phi$ along a certain subgroup $Y_j\subset \GL_l$. Moreover, $\Bs=(s_1,s_2)$ is a pair of complex numbers, and $E(h,f_{\Bs})$ is the standard Eisenstein series on $\GL_{m+n}(\BA)$ associated with a section $f_{\Bs}$ in the representation induced from $\tau_1||^{s_1-1/2}\otimes \tau_2||^{-s_2+1/2}$ on the standard Levi subgroup of $\GL_{m+n}$ with partition $(m,n)$. See \S\ref{section: the global integrals} for the unexplained notations.
\begin{theorem}\label{theorem1-intro}
The integral $I_j(\phi,f_{\Bs})$ satisfies the following properties.
\begin{enumerate}
\item  The integral $I_j(\phi,f_{\Bs})$ converges absolutely and uniformly in vertical strips for each variable $s_1,s_2$ in $\Bs=(s_1,s_2)$, away from the poles of the Eisenstein series. 

\item  For decomposable data, when $\Re(s_1+s_2)>\frac{n+m}{2}+1$, we have an Euler product factorization
\begin{align*} I_j(\phi,f_{\Bs})=  \prod_v \Psi(\rho(s_{j,m,n})W_v,\xi_{f_{\Bs,v}};j)
\end{align*}
where $\rho(s_{j,m,n})$ is the right translation by the permutation matrix $s_{j,m,n}$ defined in \eqref{eq-permutation-matrix-sj}, $W_v\in \mathcal{W}(\pi_v,\psi_v)$ is a Whittaker function, 
$\xi_{f_{\Bs,v}}\in \Ind_{P_{m,n}(F_v)}^{\GL_{m+n}(F_v)}(\mathcal{W}(\tau_{1,v},\psi_v^{-1})|\cdot|^{s_1-\frac{1}{2}}\boxtimes \mathcal{W}(\tau_{2,v},\psi_v^{-1})|\cdot|^{-s_2+1/2})$, and the local zeta integral is given by
\begin{equation*}
\Psi(W_v, \xi_{f_{\Bs,v}}; j)=\int_{N_{m+n}(F_v)\bs \GL_{m+n}(F_v)}\int_{\ov U^{j,m,n}(F_v)}W_v\left(\ov u \gamma_{m,n}\bpm h&\\ &I_{l-m-n} \epm \right)\xi_{f_{\Bs,v}}(h)d\ov udh,\end{equation*}
with $N_{m+n}$ the upper triangular unipotent subgroup of $\GL_{m+n}$, and 
\begin{align*}
\ov{U}^{j,m,n}&=\wpair{ \bpm I_m&&&&\\ x&I_{j}&&&\\ &&1&&\\ &&&I_{k}&\\ &&&y& I_n  \epm \in \GL_l}, 
\quad 
\gamma_{m,n} =\bpm I_m&&\\ &&I_{l-m-n}\\ &I_n& \epm.
\end{align*}

\item At a finite place $v$ such that all data are unramified, the local zeta integral has the following expression
\begin{equation*}
   \Psi(W_v^0,\xi_{f_{\Bs,v}}^0;j)  =\frac{L(s_1+\frac{k-j}{2},\pi_v\times \tau_{1,v})L(s_2-\frac{k-j}{2},\wt \pi_v\times \wt\tau_{2,v})}{L(s_1+s_2,\tau_{1,v}\times \wt \tau_{2,v})}
\end{equation*}
where $W_v^0,\xi_{f_{\Bs,v}}^0$ are spherical vectors suitably normalized, and $\wt \pi_v$ (resp. $\wt \tau_{2,v}$) is the contragredient representation of $\pi_v$ (resp. $\tau_{2,v}$).

\item Over nonarchimedean local fields, there exist $W_v$ and $\xi_{f_{\Bs,v}}$, such that the local zeta integral $\Psi(W_v, \xi_{f_{\Bs,v}}; j)$ is absolutely convergent and equals $1$ for all $\Bs$. Over archimedean fields, for any $\Bs$, there are choices of data $(W_v^i, \xi_{f_{\Bs,v}}^i)$ such that a finite sum of archimedean local zeta integrals $\sum_i \Psi(W^i_v,\xi_{f_{\Bs,v}}^i; j)$ is holomorphic and nonzero in a neighborhood of $\Bs$.
\end{enumerate}
\end{theorem}

The four parts of Theorem \ref{theorem1-intro} are proved in Proposition~\ref{prop-convergence-global}, Theorem~\ref{theorem: unfolding of the left side}, Proposition~\ref{proposition-unramified-computation}, and Proposition~\ref{prop-convergence}, respectively. 

We highlight that our unramified calculation in part (3) of Theorem~\ref{theorem1-intro} uses the hive model of Littlewood-Richardson coefficients developed by Knutson and Tao \cite{Knutson-Tao}. It generalizes the unramified calculation of JPSS integrals and is applicable to the local integrals of Tamir and Ben-Artzi--Soudry for $\RU_{2r+1,E/F}\times \Res_{E/F}(\GL_n)$ in \cites{Tamir, BAS:Uodd} at split places.

We also prove that in the local non-archimedean case, there exists a local gamma function $\Gamma(\Bs,\pi_v,(\tau_{1,v},\tau_{2,v}),\psi_v;j)$ such that 
$$\Psi(W_v, M(\xi_{f_{\Bs,v}}); j)=\Gamma(\Bs,\pi_v,(\tau_{1,v},\tau_{2,v}),\psi_v; j)\Psi(W_v, \xi_{f_{\Bs,v}}; j).$$
Here, $M$ denotes an intertwining operator. See \S\ref{subsection-local-functional-equation} for more details. Moreover,  we show that, up to a normalizing factor, we have
\begin{equation}\label{eq: multiplicative property of gamma}\Gamma(\Bs,\pi_v,(\tau_{1,v},\tau_{2,v}),\psi_v;j)=\frac{\gamma(s_1+(k-j)/2,\pi_v\times \tau_{1,v},\psi_v)\gamma(s_2+(j-k)/2,\wt\pi_v\times \wt\tau_{2,v},\psi_v) }{\gamma(s_1+s_2,\tau_{1,v}\times \wt \tau_{2,v},\psi_v)}.\end{equation}
Here, the gamma factors on the right side are those defined by JPSS or by Shahidi \cites{Shahidi-on-certain-L-functions, Shahidi-gamma}. See Proposition \ref{proposition-gamma-comparison} for a more precise statement.

\begin{remark}
Here are some relationships between our integrals and the previous integrals. If $n=0$, our integral and gamma factors reduce to the case of the JPSS integral and gamma factors (up to certain twists). If $l=2r+1, m=n$ and $j=r-m$, then our local zeta integral and local gamma factor are those in \cite{BAS:Uodd} and \cite{ChengWang} at split places. If $l=3$ and $m=n=1$, the above integral $I_0(\phi, f_\Bs)$ is a special case of the family of integrals considered in \cite{BumpFurusawaGinzburg}.
\end{remark}

As we mentioned above, one important application of the JPSS Rankin-Selberg integrals is the proof of the converse theorems given by Cogdell and Piatetski-Shapiro in \cites{CPS:converseI, CPS:converseII}, which roughly says that for an admissible irreducible representation $\pi$ of $\GL_l(\BA)$ for a global field $F$, if $L(s,\pi\times \tau)$ is ``nice" (see \cite{CPS:converseI}*{page 165} for the definition) for all irreducible cuspidal automorphic representations $\tau$ of $\GL_m(\BA)$ with $1\le m\le l-2$, then  $\pi$ is cuspidal automorphic. In applications to the functoriality problems, it is desirable to reduce the number of twists used in the converse theorem. In this direction, one important open question is the following:
\begin{conj}[Jacquet's global converse conjecture, see \cite{CPS:converseII}*{\S8, Conjecture 1}]
Let $F$ be a global field and let $\pi=\otimes'_v\pi_v$ be an irreducible admissible generic representation of $\GL_l(\BA)$ such that its central character is trivial on $F^\times$ and its $L$-function $L(s,\pi)$ is convergent in some half plane. If $L(s,\pi\times \tau)$ is nice for all irreducible cuspidal automorphic representation $\tau$ of $\GL_m(\BA)$ with $1\le m\le [l/2]$, then $\pi$ is cuspidal automorphic.
\end{conj}
 Many years after the original proof given in \cites{CPS:converseI, CPS:converseII}, it seems very hard to use the original JPSS integral to attack the above conjecture. We expect that our new family of integrals might be useful in the above problem. In fact, using the property \eqref{eq: multiplicative property of gamma} of the gamma factors, the condition that $L(s,\pi\times \tau)$ is nice for all irreducible cuspidal automorphic representations $\tau$ of $\GL_m$ with $1\le m\le [l/2]$ will provide us with an additional new family of equalities of integrals besides those JPSS integrals. 

In this paper, we illustrate the above idea by providing a new proof of the following local counterpart of the above conjecture; hopefully, this new proof could shed some light on the global converse conjecture mentioned above.
\begin{conj}[Jacquet's local converse conjecture]
Let $F$ be a non-archimedean local field and let $\pi_1,\pi_2$ be two generic representations of $\GL_l(F)$ with the same central character. If $\gamma(s,\pi_1\times \tau,\psi)=\gamma(s,\pi_2\times \tau,\psi)$ for all irreducible generic representation $\tau$ of $\GL_m(F)$ with $1\le m\le [l/2]$, then $\pi_1\cong \pi_2$.
\end{conj}
As proved in \cite{JNS}, one can replace the ``generic" condition in the above conjecture by ``supercuspidal". In fact, what we prove is the following
\begin{theorem}[Theorem \ref{theorem: main}]\label{theorem2-intro}
Let $F$ be a non-archimedean local field of characteristic different from 2 and let $\pi_1,\pi_2$ be two irreducible supercuspidal representations of $\GL_l(F)$ with the same central character. If $\Gamma(\Bs,\pi_1\times (\tau_1,\tau_2),\psi; 0)=\Gamma(\Bs,\pi_2\times (\tau_1,\tau_2),\psi; 0)$ for all irreducible generic representations $\tau_1$ (resp. $\tau_2$) of $\GL_m(F)$ (resp. $\GL_n(F)$) with $0\le n\le [l/2], 0\le m\le [l/2]$, then $\pi_1\cong \pi_2$.
\end{theorem}

Local converse theorems for $\GL_l$ using twists up to $l-1$ and $l-2$ have been proved in \cites{Henniart, CPS: converseII, chen}. Jacquet's local converse conjecture has been proved in \cite{chai} and \cite{Jacquet-Liu} independently. Our new contribution here is to use the new family of integrals. A proof of Jacquet's local converse conjecture along this method was promised in \cite{Sp(2r)}*{\S8.2} and in \cite{U(2r+1)}*{Introduction}, where it was believed that the integrals of $\RU_{l,E/F}\times \Res_{E/F}(\GL_m)$ at split places for a quadratic extension $E/F$ as developed in \cite{BAS:Uodd} were enough. As explained above, these integrals are just our new family of integrals when $m=n$. It turns out that we need to use the entire new family of integrals. The proof of Theorem \ref{theorem2-intro} uses partial Bessel functions developed in \cite{CST} and is indeed similar to that outlined in \cite{Sp(2r)}*{\S8.2} and in \cite{U(2r+1)}*{Introduction}.  Similar methods have been successfully used in proving local converse theorems of other classical groups over local fields and $G_2$ over finite fields, see \cites{Sp(2r), U(2r+1), Liu-Zhang: classical, Liu-Zhang: G2, Jo, HL: p-adic, SO(4)}. See \cite{Liu-Zhang: G2} for more references on local converse problems. 

Here are some differences between our proof and that given in \cite{chai} and \cite{Jacquet-Liu}. The proof of \cite{chai} only proved the equality of two partial Bessel functions on the open Bruhat cell, which is enough in the non-archimedean local field case due to the smoothness of partial Bessel functions. However, to our understanding, this is not enough to cover the finite field case, as proved in \cite{Nien} because the topology of finite fields is discrete. Our proof, which can also cover the finite field case, proves the equality of two partial Bessel functions on all Bruhat cells. On the other hand, the proof given in \cite{Jacquet-Liu} depends on Kirillov models for representations of $\GL_l$, while our proof treats $\GL_l$ as a classical group and is thus independent of the existence of Kirillov models. Moreover, our proof provides a stronger result, Theorem \ref{theorem: inductive}, which gives a condition to detect when $\gamma(s,\pi_1\times \tau,\psi)=\gamma(s,\pi_2\times \tau,\psi)$ for $\tau$ runs over generic representations of $\GL_k$ for $1\le k\le m$ for any $m\le [l/2]$. We expect it will be useful to provide finer results on gamma factors and representations of $\GL_l$. For instance, for a fixed integer $m\le [l/2]$, it might be helpful to classify representations of $\GL_l$ that could be uniquely determined by their gamma factors twisted up to $\GL_m$.

A special case of the integral $I_j(\phi, f_\Bs)$ constructed in this paper has been used to prove an algebraicity result for special values of a product of Rankin-Selberg $L$-functions in \cite{JinYan}. As the JPSS integrals, we expect our integrals to have more applications.

In this paper, we only considered the integrals that represent the product of Rankin-Selberg $L$-functions of $\GL_l\times \GL_m$ and $\GL_l\times \GL_n$ when $m+n<l$. It is natural to ask if a similar construction is generalizable to the case when $m+n\ge l$. We will address this question in future work.

The paper is organized as follows. In \S\ref{section: the global integrals}, we introduce the global integrals and discuss absolute convergence, the functional equation, and the unfolding computations of the global integrals. \S\ref{section: local theory} is devoted to the local theory of the integrals. We prove the existence of a local gamma factor $\Gamma(\Bs, \pi, (\tau_1, \tau_2), \psi;j)$, and carry out the local unramified computation for the local integrals. In \S\ref{section-preparation}, we restate Theorem \ref{theorem2-intro} and prepare some necessary tools for the proof. In particular, we recall the notions of partial Bessel functions and a result from \cite{CST}. Theorem \ref{theorem2-intro} is proved in \S\ref{section-proof}. Actually, we prove a slightly more general result (see Theorem \ref{theorem: inductive}).

\section*{Acknowledgement} We thank our advisor, Jim Cogdell, for his guidance and support over the years. It is a pleasure to dedicate this paper to him on the occasion of his 70th birthday. Some ideas in this paper grew out of the second named author thesis work under Professor Cogdell's direction; we are grateful for his many fruitful communications related to this project. Our special thanks go to Professor Terence Tao, who answered a MathOverflow question on Littlewood–Richardson coefficients posed by the second‑named author and kindly allowed us to reproduce his answer in \S \ref{subsection: Proof of Tao's formula}. We also thank Dihua Jiang, Baiying Liu, Kazuiki Morimoto and Hang Xue for discussions and helpful communications related to this paper. Finally, we thank the referee for a careful reading and valuable suggestions on the manuscript. The second‑named author acknowledges support from a start‑up grant at Huazhong University of Science and Technology.

\section{The global integrals}\label{section: the global integrals}

\subsection{Some general notation}
We fix the notation that will be used throughout the paper. For a positive integer $k$, let $I_k$ be the $k\times k$ identity matrix. Let $B_k=T_k N_k\subset \GL_k$ be the standard upper triangular Borel subgroup, with $T_k$ the group of diagonal matrices and $N_k$ the upper triangular unipotent subgroup. Let $\overline{N}_k$ be the opposite of $N_k$, i.e., $\overline{N}_k$ is the lower triangular unipotent subgroup of $\GL_k$. For positive integers $m,n$, let $\Mat_{m\times n}$ be the set of $m\times n$ matrices. We consider the following subgroups of $\GL_{m+n}$ given by 
$$ M_{m,n}=\wpair{\bpm  g_1 &\\ &g_2 \epm, g_1\in \GL_m, g_2\in \GL_n}, N_{m,n}=\wpair{\bpm I_m & X\\ &I_n \epm, X\in \Mat_{m\times n}},$$
and $P_{m,n}=M_{m,n}N_{m,n}$. Denote $w_{m,n}=\bpm &I_m\\ I_n&\epm$.

In this section, let $F$ be a global field and $\BA$ be its ring of adeles. Let $\psi$ be a fixed nontrivial additive character of $\BA$, which is trivial on $F$, and we define a character on $N_k(\BA)$, still denoted by $\psi$, by $\psi(u)=\psi(\sum_{i=1}^{k-1} u_{i,i+1})$ for $u\in N_k(\BA)$.
\subsection{Eisenstein series}
Notice that the modulus character of $P_{m,n}$ is given by
$$\delta_{P_{m,n}}(\diag(a_1,a_2))=|\det(a_1)|^n |\det(a_2)|^{-m},  \quad a_1\in \GL_m, a_2\in \GL_n.$$ 
Let $\tau_1$ (resp. $\tau_2$) be an irreducible automorphic cuspidal representation of $\GL_{m}(\BA)$ (resp. $\GL_n(\BA)$), we write $\btau=(\tau_1,\tau_2)$. Given a pair of complex numbers $\Bs:=(s_1,s_2)$, we consider the normalized induced representation
$$\RI(\Bs,\btau):=\Ind_{P_{m,n}(\BA)}^{\GL_{m+n}(\BA)}(\tau_1|\det|^{s_1-\frac{1}{2}}\otimes \tau_2|\det|^{-s_2+\frac{1}{2}}).$$
Concretely, we associate with each $u\in \RI(\Bs,\btau)$ the function $f_\Bs(h)=(u(h))(1), h\in \GL_{m+n}(\BA)$. Thus, the space $\RI(\Bs,\btau)$ consists of all functions $f_\Bs: \GL_{m+n}(\BA)\to \C$ satisfying 
$$f_{\Bs}(\diag(a,b)uh)=|\det(a)|^{s_1+\frac{n-1}{2}}|\det(b)|^{-s_2+\frac{1-m}{2}}\varphi_h(a,b),$$
where $a\in \GL_m(\BA),b\in \GL_n(\BA), u\in N_{m,n}(\BA),h\in \GL_{m+n}(\BA)$ and, for a fixed $h$, the function $(a,b)\mapsto \varphi_h(a,b)$ is a cusp form in the space of $\tau=\tau_1\boxtimes \tau_2$ of the group $M_{m,n}(\BA)=\GL_m(\BA)\times \GL_n(\BA)$.

Denote $\wh \Bs:=(s_2,s_1), 1-\wh \Bs:=(1-s_2,1-s_1)$ and $\wh \btau:=(\tau_2,\tau_1)$. There is a standard intertwining operator 
$$M_{w_{m,n}}:\RI(\Bs,\btau)\to \RI(1-\wh\Bs,\wh\btau)$$
defined by the integral
$$M_{w_{m,n}}f_{\Bs}(g)=\int_{N_{n,m}(\BA)}f_{\Bs}\left( w_{m,n} ug \right)du.$$
By \cite{MW}*{Theorem IV.1.8}, the above integral is convergent in a cone of $\Bs$ and can be extended to a meromorphic function.

Notice that  $\RI(1-\wh\Bs,\wh \btau)$ is the induced representation
$$\Ind_{P_{n,m}(\BA)}^{\GL_{m+n}(\BA)}(\tau_2|\det|^{(1-s_2)-\frac{1}{2}}\otimes \tau_1|\det|^{-(1-s_1)+\frac{1}{2}}),$$
which consists of all functions $f_{1-\wh \Bs}$ satisfying 
$$f_{1-\wh \Bs}(\diag(a,b)uh)=|\det(a)|^{1-s_2+\frac{m-1}{2}}|\det(b)|^{-(1-s_1)-\frac{n-1}{2}}\varphi_h(a,b).$$
In the above equation, $\diag(a,b)\in M_{n,m}(\BA),u\in N_{n,m}(\BA),h\in \GL_{m+n}(\BA),$ and for a fixed $h$, the function $(a,b)\mapsto \varphi_h(a,b)$ is a cusp form in the space of $\wh \tau:=\tau_2\otimes \tau_1$ of the group $M_{n,m}(\BA).$

Given $f_{\Bs}\in \RI(\Bs,\btau)$, we consider the Eisenstein series
\begin{equation}\label{eq-eisenstein} E(h,f_{\Bs})=\sum_{\gamma \in P_{m,n}(F)\bs \GL_{m+n}(F)}f_\Bs(\gamma h).\end{equation}
The above series converges to an analytic function of $\Bs$ in an open cone. When we assume that $\tau_1, \tau_2$ are unitary, the open cone can be taken as $\wpair{(s_1,s_2)\in \C^2\mid \Re(s_1+s_2)> \frac{m+n}{2}+1}$. Moreover, the above series $E(h,f_{\Bs})$ can be meromorphically continued to $\C^2$. For the above statements, see \cite{MW}*{\S II.1.5, IV.1.8}. Similarly, we can also consider the Eisenstein series 
$$E(h,f_{1-\wh \Bs})=\sum_{\gamma \in P_{n,m}(F)\bs \GL_{m+n}(F)}f_{1-\wh\Bs}(\gamma h), $$
for $f_{1-\wh\Bs}\in \RI(1-\wh \Bs,\wh \btau)$. The series $E(h,f_{1-\wh \Bs})$ is again convergent in a cone and can be continued meromorphically to $\C^2.$

\subsection{Global integrals}
Fix a positive integer $l$. Let $m,n $ be non-negative integers such that $l>m+n$. For a non-negative integer $j$ with $0\le j \le l-m-n-1$, we set $k=l-m-n-1-j\ge 0$ and consider the embedding 
$$\iota_{j,m,n}:\GL_{m+n}\to \GL_l$$
$$\bpm a&b\\ c&d \epm \mapsto \bpm I_j &&& &\\ &a&&b&\\ &&1&&\\ &c&&d&\\ &&&&I_{k} \epm$$
for $a\in \Mat_{m\times m}, b\in \Mat_{m\times n}, c\in \Mat_{n\times m}, d\in \Mat_{n\times n}.$ We also consider the permutation matrix $s_{j,m,n}\in \GL_{l}$ defined by
\begin{equation}
\label{eq-permutation-matrix-sj}
s_{j,m,n}=\bpm 0&I_m &0&0&0\\ 0&0&0&I_n&0\\ I_{j}&0&0&0&0\\ 0&0&1&0&0\\ 0&0&0&0&I_{k} \epm.    
\end{equation}
The embedding $\iota_{j,m,n}:\GL_{m+n}\to \GL_{l}$ can then be written as 
$$\iota_{j,m,n}(h)=(s_{j,m,n})^{-1}\bpm h&\\ &I_{j+1+k}\epm s_{j,m,n}, \quad h\in \GL_{m+n}.$$
 Next, we consider the subgroup $Y_{j,m,n}$ of $\GL_l$ defined by 
$$Y_{j,m,n}=\wpair{\bspm v_1&x_1&x_2&x_3 &z\\ &I_{m}& &&y_3 \\ &&1&&y_2\\ &&&I_n&y_1\\ &&&&v_2 \espm, v_1\in N_{j}, v_2\in N_{k}}.$$
To ease the notation, if $m,n$ are understood, we usually drop $m,n$ from the subscripts of the above notations. For example, we may write $Y_{j,m,n}$ as $Y_j$, $\iota_{j,m,n}$ as $\iota_j$, and $s_{j,m,n}$ as $s_j$. We define a character $\psi_{j}$ on $Y_j(F)\bs Y_j(\BA)$ by 
\begin{equation}\label{eq-psi}\psi_j \left( \bspm v_1&x_1&x_2&x_3 &z\\ &I_{m}& &&y_3 \\ &&1&&y_2\\ &&&I_n&y_1\\ &&&&v_2 \espm \right)  =\psi( v_1) +\psi(v_2)+\psi(x_{2,j})+\psi(y_{2,1}),
\end{equation}
where we write $x_2=(x_{2,1}, \cdots, x_{2,j})^t, y_2=(y_{2,1}, \cdots, y_{2,j})$.
\begin{lemma}{\label{lemma: quasi-invariance under Y}}
For $h\in \GL_{m+n}(\BA)$, $y\in Y_{j}(\BA)$, we have
\begin{enumerate}
\item $\iota_j(h)^{-1}y\iota_j(h)\in Y_j$, and
\item $\psi_j(\iota_j(h)^{-1}y\iota_j(h))=\psi_j(y).$
\end{enumerate}
\end{lemma}
\begin{proof} This follows from a simple matrix calculation.
\end{proof}

 Let $\pi$ be an irreducible cuspidal automorphic representation of $\GL_{l}(\BA)$, and for $\phi\in V_\pi$, we consider the following Fourier coefficient of $\phi$ along $Y_j$:
$$\phi_{Y_j,\psi_j}(h)=\int_{Y_j(F)\bs Y_j(\BA)}\phi(y\iota_j(h))\psi_j^{-1}(y)dy,  \quad h\in \GL_{m+n}(\BA).$$
This is called a Bessel coefficient of $\phi$. 
By Lemma \ref{lemma: quasi-invariance under Y}, $\phi_{Y,\psi}$ is left $\GL_{m+n}(F)$-invariant. Thus, for $f_\Bs\in \RI(\Bs,\btau)$, we can consider the integral
\begin{equation}\label{eq-global-integral}
I_j(\phi,f_\Bs):=\int_{\GL_{m+n}(F)\bs \GL_{m+n}(\BA)}\phi_{Y_j,\psi_j}(h)E(h,f_\Bs)dh.
\end{equation}
Similarly, we can also consider $I_j(\phi, M_{w_{m,n}}(f_{\Bs}))$. Here we remark that if $l$ is odd and $m=n$, $I_{j}(\phi, f_{\Bs})$ is a Bessel period for $(\phi, E(\cdot,f_\Bs))$ considered in \cite{GGP}*{\S 13}.

\begin{proposition}
\label{prop-convergence-global}
The integral $I_j(\phi,f_{\Bs})$ converges absolutely and uniformly in vertical strips in $\C$ for each variable $s_1,s_2$ in $\Bs=(s_1,s_2)$, away from the poles of the Eisenstein series. Moreover, away from the poles of $E(h,f_\Bs)$ and $E(h,M_{w_{m,n}}(f_\Bs))$, we have 
$$I_j(\phi,f_{\Bs})=I_j(\phi, M_{w_{m,n}}(f_{\Bs})).$$
\end{proposition}
\begin{proof}
The second statement follows from the functional equation of the Eisenstein series; see \cite{MW}*{Theorem IV.1.10}. For the first statement, it is sufficient to show that $\phi_{Y_j,\psi_j}$ is rapidly decreasing. The proof is similar to other situations that have appeared elsewhere; see \cite{BAS: Uodd}*{Lemma 2.1} for one example. We provide some details below, following the same argument as in \cite{BAS: Uodd}*{Lemma 2.1}. 

Let $\Omega$ be a compact subset of $B_{m+n}(\BA)$. Let $c$ be a real number with $0<c<1$, and we define a set $A_c$ as follows. We consider the embedding $\mathbb{R}_{>0}\hookrightarrow \BA^\times$ given by $t\mapsto (t_v)_v$, where $t_v=t$ if $v$ is an archimedean place, and $t_v=1$ if $v$ is a nonarchimedean place. 
Denote the image of this embedding by $\mathbb{R}_{+, \Delta}$. Then $A_c$ is the set of all $\diag(t_1, \dots, t_{m+n})$, such that $t_i\in \mathbb{R}_{+, \Delta}$ and $t_1\ge ct_2\ge c^2t_3\ge \dots \ge c^{m+n-1}t_{m+n}$. Then $\mathcal{S}=\Omega A_c K_{\GL_{m+n}(\BA)}$ is a Siegel domain for $\GL_{m+n}(\BA)$. Similarly, let $\mathcal{S}^\prime=\Omega^\prime A_c^\prime K_{\GL_{l}(\BA)}$ be a Siegel domain for $\GL_{l}(\BA)$, where $\iota_{j}(\Omega)\subset\Omega^\prime$ is a compact subset of $B_{l}(\BA)$ and $A^\prime_c$ is defined similarly. We take $c$ small enough and $\Omega, \Omega^\prime$ large enough so that $\GL_{l}(\BA)=\GL_{l}(F)\mathcal{S}^\prime$ and $\GL_{m+n}(\BA)=\GL_{m+n}(F)\mathcal{S}$. Now, let $h=\omega a k\in \mathcal{S}$, where $\omega\in \Omega$, $a=  \diag(t_1, \dots, t_{m+n}) \in A_c$, and $k\in K_{\GL_{m+n}(\BA)}$. Associated with $a$, we define
\begin{equation*}
b= \diag( c^{j}t_1, c^{j-1}t_1, \dots, ct_1, I_m,  t_{m},  I_n, c^{-1}t_{m+n}, c^{-2}t_{m+n}, \dots, c^{-k}t_{m+n}).
\end{equation*}
Then $b\iota_j(a)\in A_c^\prime$. Let $\Omega_b^\prime=\Omega^\prime\cup \Omega^\prime\cdot b^{-1}$. For fixed $a\in A_c$, $\Omega_b^\prime$ is a compact subset of $B_{l}(\BA)$ which contains $\Omega^\prime$. Let $\mathcal{S}^\prime_b=\Omega^\prime_b A_c^\prime K_{\GL_{l}(\BA)}$. This is a Siegel domain for $\GL_{l}(\BA)$, which contains $\mathcal{S}^\prime$. Thus, $\iota_j( h)=(\iota_j(\omega) b^{-1}) (b \iota_j(a)) \iota_j(k)\in \mathcal{S}_b^\prime$. We fix a compact subset $Y_{j,0}\subset Y_j(\BA)$ such that $Y_j(\BA)=Y_j(F)Y_{j,0}$. We may assume that $Y_{j,0}\subset \Omega^\prime$. Then we have
\begin{equation}
|\phi_{Y_j,\psi_j}(h)|\le \int_{Y_{j,0}} |\phi(y  \iota_j(\omega) b^{-1}  b \iota_j(ak)) |dy.
\label{eq-absolute-convergence-1}
\end{equation}
Let $N>0$ be given.
Since $\phi$ is rapidly decreasing in $\mathcal{S}^\prime$, there exists a constant $c_0$ such that for all $\omega^\prime\in \Omega^\prime$, $a^\prime\in A_c^\prime$, and $k^\prime\in K_{\GL_{l}(\BA)}$, we have
\begin{equation}
|\phi(\omega^\prime a^\prime k^\prime)|\le c_0 \| a^\prime\|^{-N}.
\label{eq-absolute-convergence-2}
\end{equation}
Here, $\|\cdot\|$ is the norm on $\GL_{l}(\BA)$ defined by
\begin{equation*}
\|g\|=\prod_v \| g_v\|_v
\end{equation*}
where $g\in \GL_l(\BA)$, $v$ runs over all places of $F$, and $\| g_v\|_v$ is the local norm on $\GL_l(F_v)$ defined by
\begin{equation*}
\|g_v\|_v=\max\{ |(g_v)_{i,j}|_v,  |(g_v^{-1})_{i,j}|_v: 1\le i, j \le l\}.
\end{equation*}
When passing from the Siegel domain $\mathcal{S}^\prime$ to the Siegel domain $\mathcal{S}^\prime_b$, the constant $c_0$ in \eqref{eq-absolute-convergence-2} can be replaced by $c_0 \|b^{-1}\|^{N_0}=c_0 \|b\|^{N_0}$, for some positive number $N_0$, which does not depend on $b$ (see \cite{MW}*{Sec. I.2.10, I.2.11}). Thus, in the integrand in \eqref{eq-absolute-convergence-1}, we have
\begin{equation*}
|\phi(y  \iota_j(\omega) b^{-1}b \iota_j(ak)) | \le c_0 \|b\|^{N_0} \|b\iota_j(a)\|^{-N}.
\end{equation*}
Notice that
\begin{equation*}
\begin{split}
\|b\| = & \max\{ c^{j}t_1, c^{j-1}t_1, \dots, ct_1, t_{m},  c^{-1}t_{m+n}, c^{-2}t_{m+n}, \dots, c^{-k}t_{m+n}, \ \\
 & \ c^{-j}t_1^{-1}, c^{-j+1}t_1^{-1}, \dots, c^{-1}t_1^{-1},  t_{m}^{-1},  c t_{m+n}^{-1}, c^{2}t_{m+n}^{-1}, \dots, c^{k}t_{m+n}^{-1}   \}\\
 =& \max\{ ct_1, c^{-j}t_1^{-1}, t_{m}, t_{m}^{-1} , c^{-k}t_{m+n}, ct_{m+n}^{-1} \}\\
 \le & c^{\max\{1, -j, -k\}}\|a\|
\end{split}
\end{equation*}
and
\begin{equation*}
\begin{split}
\|b \iota_j(a)\| = & \max\{ c^{j}t_1, c^{j-1}t_1, \dots, ct_1,  t_1, t_2, \dots,   t_{m+n}, c^{-1}t_{m+n}, c^{-2}t_{m+n}\dots, c^{-k}t_{m+n}, \ \\
 & \ c^{-j}t_1^{-1}, c^{-j+1}t_1^{-1}, \dots, c^{-1}t_1^{-1},  t_1^{-1}, \dots, t_{m+n}^{-1},    c t_{m+n}^{-1}, c^2 t_{m+n}^{-1}, \dots, c^{k}t_{m+n}^{-1}\}  \\
  \ge &\max\{ t_1, t_2, \dots,   t_{m+n}, t_1^{-1}, t_2^{-1}, \dots,  t_{m+n}^{-1} \}\\
 = & \|a\|.
\end{split}
\end{equation*}
We conclude that 
\begin{equation}
|\phi(y \iota_j(\omega a k))| \le c_1 \|a\|^{N_0-N}
\label{eq-absolute-convergence-3}
\end{equation}
where $c_1$ is a positive constant depending on $c$ and $c_0$. Since $Y_{j,0}$ is compact, we combine \eqref{eq-absolute-convergence-2} and \eqref{eq-absolute-convergence-3} to conclude that $\phi_{Y_j, \psi_j}$ is rapidly decreasing in $\mathcal{S}$. This completes the proof. 
\end{proof}

\subsection{Unfolding of the global integral \texorpdfstring{$I_j(\phi,f_{\Bs})$}{}}

For integers $m,n\ge 0$, denote 
$$Z_{m,n}=\left\{\bpm I_m &0&z\\ &1&0\\ &&I_n \epm: z\in \Mat_{m\times n} \right\}\subset \GL_{m+n+1}.$$
For a cusp form  $\phi$ on $\GL_{m+n+1}(F)\bs \GL_{m+n+1}(\BA)$, we define its constant term along $Z_{m,n}$ by
$$\phi_{Z_{m,n}}(g)=\int_{Z_{m,n}(F)\backslash Z_{m,n}(\BA)}\phi\left(zg \right)dz.$$
We have the following expansion of $\phi_{Z_{m,n}}$.
\begin{lemma}\label{lemma: integral over Z}
For $\phi\in \CA_0(\GL_{m+n+1})$, the space of cusp forms on $\GL_{m+n+1}(F)\bs \GL_{m+n+1}(\BA)$, we have
\begin{align*}
\phi_{Z_{m,n}}(g)=\sum_{\substack{\gamma_1\in N_{m}(F)\backslash \GL_m(F), \\
\gamma_2\in N_{n}(F)\backslash \GL_n(F)}}W_\phi^{\psi}\left(\bpm \gamma_1 &&\\ &1&\\ &&\gamma_2  \epm g\right),
\end{align*}
with absolute and uniform convergence on compact subsets of $g.$ Here $W_\phi^{\psi}$ is the $\psi$-Whittaker function of $\phi$ defined by
\begin{equation}
\label{eq-Whittaker-function}
	W_\phi^{\psi}(g)=\int_{N_{m+n+1}(F)\backslash N_{m+n+1}(\BA)} \phi(ug)\psi^{-1}(u)du, \quad \psi(u)=\psi(\sum_{k=1}^{m+n} u_{k,k+1}).
\end{equation}
\end{lemma}
  Note that when $n=0$, the above expansion is just the usual Fourier expansion of cusp forms, due to Piatetski-Shapiro \cite{PS:Fourier} and Shalika \cite{Shalika:multiplicityone}. On the other hand, the above version expansion is an easy consequence of the results of Piatetski-Shapiro and Shalika. We give a sketch of the proof below.
  \begin{proof}
  Let $$Q_{m}=\wpair{\bpm g_1 &x\\ &1 \epm: g_1\in \GL_m, x\in \Mat_{m\times 1}}$$
  be the usual mirabolic subgroup of $\GL_{m+1}$. We consider the function $\phi_1=\phi_1^g$ on $Q_m(F)\backslash Q_m(\BA)$ defined by 
  $$\phi_1\left(\bpm g_1 &x\\ &1 \epm \right)=\phi_{Z_{m,n}}\left(\bpm g_1&x &\\ &1&\\ &&I_n \epm  g \right).$$
  Then $\phi_1$ is a cuspidal automorphic form on $ Q_m(F)\backslash Q_m(\BA)$ in the sense that for any parabolic subgroup $P=MU$ of $Q_m$ with unipotent radical $U$, we have 
  $$\int_{U(F)\backslash U(\BA)}\phi_1(uq)du=0, \quad  \forall q\in Q_m(\BA).$$
  This can be checked easily using cuspidality of $\phi$, see \cite{Cogdell:IAS}*{Lemma 2.2} for a similar situation. Thus by the Fourier expansion for $\phi_1$, see \cites{PS:Fourier, Shalika:multiplicityone}, we get that 
  $$\phi_1(I_{m+1})=\sum_{\gamma_1\in N_{m}(F)\backslash \GL_{m}(F)}W_{\phi_1}^{\psi}\left(\bpm \gamma_1 &\\ &1 \epm\right),$$
 with absolute and uniform convergence, where $W_{\phi_1}^\psi$ is the standard $\psi$-Whittaker function of $\phi_1$. Plugging in the definitions, we get that
  \begin{align*}
  \phi_{Z_{m,n}}(g)&=\sum_{ \gamma_1\in N_{m}(F)\backslash \GL_{m}(F)}\int \phi\left(\bpm u &x&z\\ &1&\\ &&I_n \epm \bpm \gamma_1&&\\&1&\\ &&I_n \epm g\right)\psi^{-1}(u)\psi^{-1}(x_m)dudxdz,
  \end{align*}
  where $u=(u_{ij})\in  N_{m}(\BA)$, $\psi^{-1}(u)=\psi^{-1}(\sum_i u_{i,i+1})$ and $x_m$ is the last component of $x$. Similarly, we consider the mirabolic subgroup $Q_n'$ of $\GL_{n+1}$ of the form
  $$Q_n'=\wpair{\bpm 1&y\\ 0&g_2 \epm, y\in \Mat_{1\times n}, g_2\in \GL_n}.$$
 For fixed $\gamma_1$ and $g$, we consider the function $\phi_2$ on $Q_n'(F)\backslash Q_n'(\BA)$ defined by 
 \begin{align*}
 \phi_2\left(\bpm 1&y\\ 0&g_2 \epm \right)&=\int \phi\left(\bpm u &x&z\\ &1&y\\ &&g_2 \epm \bpm \gamma_1&&\\&1&\\ &&I_n \epm g\right)\psi^{-1}(u)\psi^{-1}(x_m)dudxdz.
 \end{align*}
 Again, $\phi_2$ is a cusp form on $Q_n'(F)\backslash Q_n'(\BA)$. By a slightly variant form of the Fourier expansion, see for example \cite{CPS:converse II}*{\S1, Proposition}, we have 
 $$\phi_2(I_{n+1})=\sum_{\gamma_2\in N_{n}(F)\backslash \GL_{n}(F)}W_{\phi_2}^{\psi}\left(\bpm 1&\\&\gamma_2 \epm\right).$$
 The above series is also absolutely convergent.
 Note that 
 \begin{align*}W_{\phi_2}^\psi\left(\bpm 1&\\&\gamma_2 \epm \right)
 &=\int \phi\left( \bpm u&x&z\\ &1&y\\ &&v \epm  \bpm \gamma_1&&\\&1&\\ &&\gamma_2\epm g\right)\psi^{-1}(u)\psi^{-1}(v)\psi^{-1}(x_m+y_1)dxdydudvdz\\
 &=W_\phi^\psi\left( \bpm \gamma_1&&\\&1&\\ &&\gamma_2\epm g \right),
 \end{align*}
 where $y_1$ in the first integral is the first component of $y$. The result follows.
  \end{proof}

\begin{theorem}\label{theorem: unfolding of the left side}
The integral $I_j(\phi,f_{\Bs})$ is Eulerian. More precisely, when $\Bs$ is in the cone of absolute convergence of the Eisenstein series \eqref{eq-eisenstein}, we have
\begin{align}\label{eq-unfolding}I_j(\phi,f_{\Bs})=\int_{N_{m+n}(\BA)\bs \GL_{m+n}(\BA)}\int_{\ov{U}^{j,m,n}(\BA)}W_{\phi}^\psi\left( \ov u \eta_{j} \iota_{j}(h)\right)\xi_{f_\Bs}^{\psi^{-1}}(h)d\ov udh
\end{align}
where 
\begin{align*}
\ov{U}^{j,m,n}&=\wpair{\ov u(x,y)=\bpm I_m&&&&\\ x&I_{j}&&&\\ &&1&&\\ &&&I_{k}&\\ &&&y& I_n  \epm :  x\in \Mat_{j\times m}, y\in \Mat_{n\times k}},\\
\eta_j&=\eta_{j,m,n}=\bpm &I_m&&&\\ I_{j}&&&&\\ &&1&&\\ &&&&I_{k}\\ &&&I_n& \epm,\\
\xi_{f_\Bs}^{\psi^{-1}}(h)&=\int_{N_{m}(F)\bs N_{m}(\BA)\times N_{n}(F)\bs N_{n}(\BA)}
f_{\Bs}\left(\bpm u_1 & \\ & u_2 \epm h \right)\psi(u_1)\psi(u_2)du_1 du_2.
\end{align*}
\end{theorem}

\begin{proof}
For notational convenience, we set $[G]:=G(F)\bs G(\BA)$ for an algebraic group $G$ over $F$. In what follows, we assume that $\Bs$ is in the cone of absolute convergence, so that the Eisenstein series $I(h,f_{\Bs})$ is given by the convergent sum \eqref{eq-eisenstein}. By plugging \eqref{eq-eisenstein} into the definition of $I_j(\phi, f_{\Bs})$ (see \eqref{eq-global-integral}), we get
\begin{align}\label{eq: I(phi,f)}
I_j(\phi,f_s)&=\int_{P_{m,n}(F)\bs\GL_{m+n}(\BA)}\phi_{Y_j,\psi_j}(h)f_\Bs(h)dh \\
&=\int_{M_{m,n}(F)N_{m,n}(\BA)\bs\GL_{m+n}(\BA)}\int_{[N_{m,n}]}\phi_{Y_j,\psi_j}(uh)duf_\Bs(h)dh \nonumber\\
&=\int_{M_{m,n}(F)N_{m,n}(\BA)\bs\GL_{m+n}(\BA)}\phi_{Y_j,\psi_j,N_{m,n}}(h)f_\Bs(h)dh,\nonumber
\end{align}
where 
\begin{align*}
\phi_{Y_j,\psi_j,N_{m,n}}(h)&=\int_{[N_{m,n}]}\phi_{Y_j,\psi_j}(uh)du.
\end{align*}
By the definition of $\phi_{Y_j,\psi_j}$, together with the left $\GL_{l}(F)$ invariance of $\phi$, we get 
\begin{align}\label{eq-YN-coefficient}
\begin{split}
\phi_{Y_j,\psi_j,N_{m,n}}(h)&=\int_{[Y_j]\times [N_{m,n}]}\phi(y\iota_j(uh))\psi_j^{-1}(y)dudy\\
&=\int_{[Y_j]\times [\iota_j(N_{m,n})]}\phi(\eta_j y\iota_j(u)\eta_j^{-1} \eta_j \iota_j (h))\psi_j^{-1}(y)dudy.
\end{split}
\end{align}
Write
\begin{align*} y=\bspm v_1&x_1&x_2&x_3 &z\\ &I_{m}& &&y_3 \\ &&1&&y_2\\ &&&I_n&y_1\\ &&&&v_2 \espm \in Y_j(\BA), \quad  u=\bpm I_m & t\\ &I_n\epm \in N_{m,n}(\BA),\end{align*}
then we have 
\begin{align*}\eta_{j} y \iota_j(u) \eta_{j}^{-1}=\bspm I_m&0&0&y_3&t\\ x_1 &v_1&x_2&z&x_3\\ 0 &0&1&y_2&0\\ 0 &0&0&v_2&0\\ 0&0&0&y_1& I_n \espm,\end{align*}
where $v_1\in N_{j}(\BA),v_2\in N_{k}(\BA), (x_1,x_2,x_3)\in \Mat_{j\times (m+n+1)}(\BA), z\in \Mat_{j\times k}(\BA), $ $y_3\in \Mat_{m\times k}(\BA), $ $y_2\in \Mat_{1\times k}(\BA), y_1\in \Mat_{n\times k}(\BA)$  and  $t\in \Mat_{m\times n}(\BA)$. 
Then \eqref{eq-YN-coefficient} becomes
\begin{equation}\label{eq: phi(Y,psi,N)}\phi_{Y_j,\psi_j,N_{m,n}}(h)=\int_{[Y_j]\times [\iota_j(N_{m,n})]} \phi\left(\bspm I_m&0&0&y_3&t\\ x_1 &v_1&x_2&z&x_3\\ 0 &0&1&y_2&0\\ 0 &0&0&v_2&0\\ 0&0&0&y_1& I_n \espm \eta_{j} \iota_j(h) \right)\psi_j^{-1}(y)dydu.\end{equation}
 Write 
$$Z=\bpm y_3 &t\\ z& x_3 \epm \in \Mat_{(m+j)\times (n+k)}(\BA).$$
The right side integral of \eqref{eq: phi(Y,psi,N)} contains the following inner integral
\begin{align*}
\int_{[\Mat_{(m+j)\times (n+k)}]}\phi\left( \bpm I_{m+j} &&Z\\ &1&\\ &&I_{n+k} \epm  g\right)dZ,
\end{align*}
which is 
\begin{align}\label{eq: inner integral over Z}
\sum_{\substack{\gamma_1\in N_{m+j}(F)\bs \GL_{m+j}(F) \\ \gamma_2\in N_{n+k}(F)\bs \GL_{n+k}(F)}}W_\phi^\psi\left( \bpm \gamma_1 &&\\ &1&\\ &&\gamma_2 \epm g \right)
\end{align}
by Lemma \ref{lemma: integral over Z}. Plugging \eqref{eq: inner integral over Z} into \eqref{eq: phi(Y,psi,N)}, we get 
\begin{align}\label{eq: phi Y psi N 2}
\phi_{Y_j,\psi_j,N_{m,n}}(h)&=\sum_{\gamma_1,\gamma_2}\int W_\phi^{\psi}\left(\bpm \gamma_1 &&\\ &1&\\ &&\gamma_2 \epm \bspm I_m&0&0&0&0\\ x_1 &v_1&x_2&0&0\\ 0 &0&1&y_2&0\\ 0 &0&0&v_2&0\\ 0&0&0&y_1& I_n \espm \eta_{j} \iota_j(h) \right)\psi_j^{-1}(y)dy.
\end{align}
To simplify the above integral \eqref{eq: phi Y psi N 2}, we consider its inner integral with respect to $x_2=[x^1,\dots, x^j]^t\in \Mat_{j\times 1}(\BA)$ first, which is 
\begin{align*}
\int_{(F\bs \BA)^j}W_\phi^\psi\left(\bpm \gamma_1 &&\\ &1&\\ &&\gamma_2 \epm \bspm I_m&0&0&0&0\\ 0 &I_j&x_2&0&0\\ 0 &0&1&0&0\\ 0 &0&0&I_k&0\\ 0&0&0&0& I_n \espm \bspm I_m&0&0&0&0\\ x_1 &v_1&0&0&0\\ 0 &0&1&y_2&0\\ 0 &0&0&v_2&0\\ 0&0&0&y_1& I_n \espm \eta_{j} \iota_j(h) \right)\psi^{-1}(x^j)dx_2.
\end{align*}
Write $\gamma_1=(\gamma_{pq})_{1\le p,q\le m+j}$, then we have 
$$\gamma_1\bpm 0\\ x_2 \epm=\bpm *\\ *\\ \vdots\\ \gamma_{m+j,m+1}x^1+\gamma_{m+j,m+2}x^2+\dots+\gamma_{m+j,m+j}x^j \epm.$$
After conjugation, we get 
\begin{align*}
 &W_\phi^\psi\left(\bpm \gamma_1 &&\\ &1&\\ &&\gamma_2 \epm \bspm I_m&0&0&0&0\\ 0 &I_j&x_2&0&0\\ 0 &0&1&0&0\\ 0 &0&0&I_k&0\\ 0&0&0&0& I_n \espm g \right)\\
& \qquad =\psi(\gamma_{m+j,m+1}x^1+\dots+\gamma_{m+j,m+j}x^j )
 \cdot W_\phi^{\psi}\left(\bpm \gamma_1&&\\ &1&\\ &&\gamma_2  \epm g\right),
\end{align*}
with 
$$g= \bspm I_m&0&0&0&0\\ x_1 &v_1&0&0&0\\ 0 &0&1&y_2&0\\ 0 &0&0&v_2&0\\ 0&0&0&y_1& I_n \espm \eta_{j} \iota_j(h).$$
Thus the inner integral of \eqref{eq: phi Y psi N 2} with respect to $x_2$ is 
\begin{align*}
\int_{(F\bs \BA)^j}\psi(\gamma_{m+j,m+1}x^1+\dots+(\gamma_{m+j,m+j}-1)x^j)dx^1\dots dx^j W_\phi^{\psi}\left(\bspm \gamma_1&&\\ &1&\\ &&\gamma_2  \espm g\right).
\end{align*}
The above integral over $x^1,\dots,x^j$ is $1$ if $\gamma_{m+j,m+1}=\dots=\gamma_{m+j,m+j-1}=0$ and  $\gamma_{m+j,m+j}=1$, and is zero otherwise. Note that if $\gamma_{m+j,m+1}=\dots=\gamma_{m+j,m+j-1}=0$, as an element of the coset $N_{m+j}(F)\bs \GL_{m+j}(F)$, we can write 
$$\gamma_1=\bpm \gamma_1' &\\ &1 \epm \bpm I_{m} &&\\ &I_{j-1} &\\ \xi && 1 \epm,$$
with $\gamma_1'\in N_{m+j-1}(F)\bs \GL_{m+j-1}(F), \xi\in\Mat_{1\times m}(F).$ By changing the summation notation, integral \eqref{eq: phi Y psi N 2} becomes
\begin{align}\label{eq: phi Y psi N 3}
&\phi_{Y_j,\psi_j,N_{m,n}}(h)=\sum_{\substack{\gamma_1\in N_{m+j-1}(F)\bs \GL_{m+j-1}(F)\\ \gamma_2\in N_{n+k}(F)\bs \GL_{n+k}(F)}}\sum_{\xi\in F^m}\\
& \qquad \int W_\phi^{\psi}\left(\bspm \gamma_1 &&\\ &I_2&\\ &&\gamma_2 \espm \bspm I_m &&&&\\ &I_{j-1} &&&\\ \xi &&1&&\\ &&&1&\\ &&&&I_{n+k} \espm \bspm I_m&0&0&0&0\\ x_1 &v_1&0&0&0\\ 0 &0&1&y_2&0\\ 0 &0&0&v_2&0\\ 0&0&0&y_1& I_n \espm \eta_{j} \iota_j(h) \right) \nonumber\\
& \qquad \cdot \psi^{-1}(v_1)\psi^{-1}(v_2)\psi^{-1}(y^1)dx_1dy_1dy_2dv_1dv_2.\nonumber
\end{align}
Here $y^1$ is the first component of the element $y_2\in \Mat_{1\times k}(\BA).$ In other words, $y_2=[y^1,\dots,y^k]\in \Mat_{1\times k}(\BA).$ In \eqref{eq: phi Y psi N 3}, the summation over $\xi$ could be absorbed  into the integral over the last row of $x_1$. Thus, we get
\begin{align}\label{eq: phi Y psi N 4}
\phi_{Y_j,\psi_j,N_{m,n}}(h)=&\sum_{\substack{\gamma_1\in N_{m+j-1}(F)\bs \GL_{m+j-1}(F)\\ \gamma_2\in N_{n+k}(F)\bs \GL_{n+k}(F)}}\int_{(F\bs \BA)^*}\int_{\BA^m}\\
&\int W_\phi^{\psi}\left(\bpm \gamma_1 &&\\ &I_2&\\ &&\gamma_2 \epm \bspm I_m&0&0&0&0&0\\ x_1' &v_1'&p&0&0&0\\ (x_{j1},\dots, x_{jm}) &0&1&0&0&0 \\ 0&0 &0&1&y_2&0\\ 0&0 &0&0&v_2&0\\ 0&0&0&0&y_1& I_n \espm \eta_{j} \iota_j(h) \right) \nonumber\\
&\quad \cdot \psi^{-1}(v_1)\psi^{-1}(v_2)\psi^{-1}(y^1)(\prod_{t=1}^m dx_{jt})dx_1'dy_1dy_2dv_1dv_2,\nonumber
\end{align}
where we wrote $x_1=\bpm x_1'\\ (x_{j1},\dots, x_{jm}) \epm$, $v_1=\bpm v_1'& p \\ &1 \epm$ with $p\in \Mat_{(j-1)\times 1}(\BA)$, and $*$ in $(F\bs \BA)^*$ denotes the number of variables other than the part in $(x_{j1},\dots,x_{jm})$. We next compute the inner integral over the $p$-part, which is similar as above. Note that $\psi(v_1)=\psi(v_1')\psi(p^{j-1})$, where $p=(p^1,\dots,p^{j-1})^t$. For $\gamma_1=(\gamma_{ab})_{1\le a,b\le m+j-1}\in \GL_{m+j-1}(F)$, and $p=(p^1,\dots,p^{j-1})^t$ we have
$$\gamma_1\bpm 0_{m\times 1} \\ p \epm =\bpm *\\ \vdots\\ *\\ \gamma_{m+j-1,m+1}p^1+\dots \gamma_{m+j-1,m+j-1}p^{j-1}  \epm.$$
Thus the inner integral over $p$ in \eqref{eq: phi Y psi N 4} is 
\begin{align*}
\int_{(F\bs \BA)^{j-1}}\psi( \gamma_{m+j-1,m+1}p^1+\dots+(\gamma_{m+j-1,m+j-1}-1)p^{j-1})\prod_t dp^t W_\phi^{\psi}\left( \bpm \gamma_1&&\\ &1&\\ &&\gamma_2 \epm g \right),
\end{align*}
for certain appropriate $g$ which should be self-evident from the context. The above integral is 1 if $\gamma_{m+j-1,m+1}=\dots=\gamma_{m+j-1,m+j-2}=0$ and $\gamma_{m+j-1,m+j-1}=1$, and is zero otherwise. In this case, we can write that $$\gamma_1=\bpm \gamma_1' &\\ &1 \epm \bpm I_m &&\\ &I_{j-2} &\\ \xi&&1 \epm$$
as an element in the coset $N_{m+j-1}(F)\bs \GL_{m+j-1}(F)$, where $\gamma_1'\in N_{m+j-2}(F)\bs \GL_{m+j-2}(F), \xi\in F^m$. Similarly as above, by absorbing the summation over $\xi$, we get that 
\begin{align*}
&\phi_{Y_j,\psi_j,N_{m,n}}(h)=\sum_{\substack{\gamma_1\in N_{m+j-2}(F)\bs \GL_{m+j-2}(F)\\ \gamma_2\in N_{n+k}(F)\bs \GL_{n+k}(F)}}\int_{(F\bs \BA)^*}\int_{\BA^{2m}}\\
&\qquad \int W_\phi^{\psi}\left(\bpm \gamma_1 &&\\ &I_3&\\ &&\gamma_2 \epm \bpm I_m&0&0&0&0&0&0\\ x_1'' &v_1''&p'&0&0&0&0\\  (x_{j-1,1},\dots,x_{j-1,m}) &0 &1 &0&0&0&0   \\ (x_{j1},\dots, x_{jm}) &0&0&1&0&0&0 \\ 0&0 &0&0&1&y_2&0\\ 0&0 &0&0&0&v_2&0\\ 0&0&0&0&0&y_1& I_n \epm \eta_{j} \iota_j(h) \right) \nonumber\\
&\qquad \quad \cdot \psi^{-1}(v_1')\psi^{-1}(v_2)\psi^{-1}(y^1)(\prod_{i=j-1}^j\prod_{t=1}^m dx_{it})dx_1''dy_1dy_2dv_1'dv_2,\nonumber
\end{align*}
where $v_1'=\bpm v_1''&p'\\ &1 \epm$.
An inductive argument shows that
\begin{align*}
\phi_{Y_j,\psi_j,N_{m,n}}(h)=&\sum_{\substack{\gamma_1\in N_{m}(F)\bs \GL_{m}(F)\\ \gamma_2\in N_{n+k}(F)\bs \GL_{n+k}(F)}}\int_{(F\bs \BA)^*}\int_{\Mat_{j\times m}(\BA)}\\
& W_\phi^\psi\left( \bpm \gamma_1&&\\ &I_{j+1}&\\ &&\gamma_2 \epm \bpm I_m &&&&&\\ x& I_{j} &&&&\\ &&1&&&\\ &&&1&y_2&\\ &&&&v_2&\\ &&&&y_1&I_n  \epm \eta_{j} \iota_j(h) \right)\\
&\cdot \psi^{-1}(v_2)\psi^{-1}(y^1)dxdy_1dy_2 dv_2.
\end{align*}
The integral over $y_1,y_2,v_2$ can be done similarly and we have
\begin{align*}
\phi_{Y_j,\psi_j,N_{m,n}}(h)=&\sum_{\substack{\gamma_1\in N_{m}(F)\bs \GL_{m}(F)\\ \gamma_2\in N_n(F)\bs \GL_n(F)}} \int_{\Mat_{j\times m}(\BA)} \int_{\Mat_{n\times k}(\BA)}\\
& \quad W_\phi^\psi\left( \bpm \gamma_1&&\\ &I_{j+k+1}&\\ &&\gamma_2 \epm \bpm I_m &&&&\\ x& I_{j} &&&\\ &&1&&\\ &&&I_{k}&\\ &&&y &I_n \epm \eta_{j} \iota_j(h) \right)dydx\\
=&\sum_{\substack{\gamma_1\in N_{m}(F)\bs \GL_{m}(F)\\ \gamma_2\in N_n(F)\bs \GL_n(F)}} \int_{\Mat_{j \times m}(\BA)} \int_{\Mat_{n\times k}(\BA)}\\
&\quad W_\phi^\psi\left( \bpm I_m &&&&\\ x& I_{j} &&&\\ &&1&&\\ &&&I_{k}&\\ &&&y &I_n \epm \eta_{j} \iota_j\left(\bpm \gamma_1&\\ &\gamma_2 \epm h \right) \right)dydx
\end{align*}
We now plug the above formula into \eqref{eq: I(phi,f)} to get
\begin{equation}
\begin{split}
I_j(\phi,f_{\Bs})=&\int_{M_{m,n}(F)N_{m,n}(\BA)\bs \GL_{m+n}(\BA)}\phi_{Y_j,\psi_j,N_{m,n}}(h)f_{\Bs}(h)dh\\
=&\int_{(N_m(F)\times N_n(F))N_{m,n}(\BA)\bs \GL_{m+n}(\BA)}\int_{\ov U^{j,m,n}(\BA)}
W_\phi^\psi\left( \ov u \eta_{j} \iota_j(  h ) \right)f_{\Bs}(h)d\ov udh.
\end{split}
\label{eq-unfolding-eq1}
\end{equation}
In order to justify this step, we need to show that the double integral in the second line of \eqref{eq-unfolding-eq1} converges absolutely. This will be done in Subsection~\ref{subsection-convergence-and-justifications}. From \eqref{eq-unfolding-eq1}, we obtain
\begin{align*}
I_j(\phi,f_{\Bs})=&\int_{N_{m+n}(\BA)\bs \GL_{m+n}(\BA)}\int_{\ov U^{j,m,n}(\BA)}\int_{N_m(F)\bs N_m(\BA)}\int_{N_n(F)\bs N_n(\BA)}W_\phi^\psi\left(\ov u \eta_{j}\iota_j\left(\bpm u_1 &\\ &u_2 \epm h \right)\right)\\
&\cdot f_\Bs\left(\bpm u_1 &\\ &u_2 \epm h \right)du_2du_1dydxdh\\
=&\int_{N_{m+n}(\BA)\bs \GL_{m+n}(\BA)}\int_{\ov U^{j,m,n}(\BA)}\int_{N_m(F)\bs N_m(\BA)}\int_{N_n(F)\bs N_n(\BA)}W_\phi^\psi\left(\ov u \eta_{j}\iota_j\left( h \right)\right)\\
&\cdot \psi^{-1}(u_1)\psi^{-1}(u_2) f_\Bs\left(\bpm u_1 &\\ &u_2 \epm h \right)du_2du_1dydxdh\\
=&\int_{N_{m+n}(\BA)\bs \GL_{m+n}(\BA)}\int_{\ov{U}^{j,m,n}(\BA)}W_{\phi}^\psi\left( \ov u \eta_{j} \iota_{j}(h)\right)\xi_{f_\Bs}^{\psi^{-1}}(h)d\ov udh.
\end{align*}
The result follows.
\end{proof}
\begin{remark}\label{remark-unfolding-left}
Denote 
\begin{equation}
\label{eq-gamma-mn}
\gamma_{m,n}:=\eta_{j,m,n}s_{j,m,n}^{-1}=\bpm I_m&&\\ &&I_{l-m-n}\\ &I_n& \epm.
\end{equation}
Then \eqref{eq-unfolding} can be rewritten as 
\begin{align*}I_j(\phi,f_{\Bs})=\int_{N_{m+n}(\BA)\bs \GL_{m+n}(\BA)}\int_{\ov{U}^{j,m,n}(\BA)}W_{\phi}^\psi\left( \ov u \gamma_{m,n}\bpm h&\\ & I_{l-m-n} \epm s_{j,m,n}\right)\xi_{f_\Bs}^{\psi^{-1}}(h)d\ov udh.
\end{align*}
\end{remark}

\subsection{Unfolding of \texorpdfstring{$I_j(\phi,M_{w_{m,n}}(f_\Bs))$}{}}
\begin{theorem}\label{theorem: unfolding of the right side}
The integral $I_j(\phi,M_{w_{m,n}}(f_\Bs))$ is Eulerian. More precisely, in the region of absolute convergence of the Eisenstein series $ E(h, M_{w_{m,n}}(f_\Bs))$, we have
\begin{align*}
I_j(\phi,\wt f_{\Bs})=\int_{N_{n+m}(\BA)\bs \GL_{n+m}(\BA)} \int_{\ov V^{j,m,n}} W_\phi^\psi\left(\ov u \gamma_{n,m}\bpm h&\\ &I_{l-m-n}\epm s_{j,m,n}\right) \xi_{\wt f_\Bs}^{\psi^{-1}}(h)d\ov udh,
\end{align*}
where 
\begin{align*}
\wt f_\Bs&=M_{w_{m,n}}(f_\Bs),\\
\ov V^{j,m,n}&=\wpair{\bpm I_n&0&0&0&0\\ x&I_{j}&0&0&0\\ &&1&0&0\\ &&&I_{k}&0\\ &&&y&I_{m}  \epm : x\in \Mat_{j\times n}, y\in \Mat_{m\times k}}=\ov{U}^{j,n,m},\\
\gamma_{n,m}&=\bpm I_n&&\\ &&I_{l-m-n}\\ &I_m& \epm, \\
\xi_{\wt f_\Bs}^{\psi^{-1}}(h)&=\int_{N_n(F)\bs N_n(\BA)\times N_m(F)\bs N_m(\BA)}\wt f_{\Bs}\left(\bpm u_1&\\ &u_2 \epm h \right)\psi(u_1)\psi(u_2)du_1du_2.
\end{align*}
\end{theorem}
The proof is similar to the proof of Theorem \ref{theorem: unfolding of the left side}.  We give some details for completeness. 
\begin{proof}
In the following, we assume that $m\ge n$. If $n\le m$, the matrix calculation performed below is a little bit different, but other parts of the proof go through, and the result is the same. Note that $\wt f_\Bs\in \RI(1-\wh \Bs, \wh \btau)$ is left invariant under $N_{n,m}(\BA)$. When the Eisenstein series $E(h, \wt f_{\Bs}) $ is absolutely convergent, we can replace it by its infinite series definition to obtain  
\begin{equation}
\begin{split}
I_j(\phi,\wt f_\Bs)&=\int_{P_{n,m}(F)\bs \GL_{n+m}(\BA)}\phi_{Y_j,\psi_j}(h)\wt f_{\Bs}(h)dh \\
&=\int_{M_{n,m}(F)N_{n,m}(\BA)\bs \GL_{n+m}(\BA)}\int_{N_{n,m}(F)\bs N_{n,m}(\BA)}\phi_{Y_j,\psi_j}(uh)\wt {f}_\Bs(h)dudh \\
&=\int_{M_{n,m}(F)N_{n,m}(\BA)\bs \GL_{n+m}(\BA)}\phi_{Y_j,\psi_j,N_{n,m}}(h)\wt {f}_\Bs(h)dh,
\end{split}
\label{eq: I(phi,wt f)}
\end{equation}
where
\begin{align*}\phi_{Y_j,\psi_j,N_{n,m}}(h):&=\int_{N_{n,m}(F)\bs N_{n,m}(\BA)}\phi_{Y_j,\psi_j}(uh) du\\
&=\int_{N_{n,m}(F)\bs N_{n,m}(\BA)}\int_{Y_j(F)\bs Y_j(\BA)}\phi(y\iota_j(u)\iota_j(h))\psi_j^{-1}(y)dydu.
\end{align*}
Since $\phi$ is left $\GL_{l}(F)$-invariant, we have
$$\phi(y\iota_j(u)\iota_j(h))=\phi\left(\gamma_{n,m}s_{j,m,n}y\iota_j(u)s_{j,m,n}^{-1} \gamma_{n,m}^{-1} \gamma_{n,m}\bpm h&\\ &I_{l-m-n} \epm s_{j,m,n} \right) .$$
Write 
\begin{align*}
y=\bpm v_1&x_1&x_1'&x_2&x_3&z\\ &I_n&&&&y_3\\ &&I_{m-n}&&&y_3'\\ &&&1&&y_2\\ &&&&I_n&y'_1\\ &&&&&v_2 \epm \in Y(\BA), \quad u=\bpm I_n&t_1&t_2\\ &I_{m-n}&\\ &&I_n \epm \in N_{n,m}(\BA),
\end{align*}
with $v_1\in N_{j}(\BA), v_2\in N_{k}(\BA)$ and other variables in appropriate matrices spaces. A matrix calculation shows that 
$$\gamma_{n,m}s_{j,m,n}y\iota_j(u)s_{j,m,n}^{-1} \gamma_{n,m}^{-1}=\bpm I_n&0&0&y_3&t_1&t_2\\ x_1&v_1&x_2&z&x_1'+x_1t_1&x_3+x_1t_2\\ &&1&y_2&0&0\\ &&&v_2&0&0\\ &&&y_3'&I_{m-n}&0\\ &&&y'_1&0&I_n \epm.$$
Thus we get
\begin{align*}
&\phi_{Y_j,\psi_j,N_{n,m}}(h)\\
&\qquad =\int_{[Y_j]\times [N_{n,m}]}\phi\left(\bspm I_n&0&0&y_3&t_1&t_2\\ x_1&v_1&x_2&z&x_1'&x_3\\ &&1&y_2&0&0\\ &&&v_2&0&0\\ &&&y_3'&I_{m-n}&0\\ &&&y'_1&0&I_n \espm  \gamma_{n,m}\bpm h&\\ &I_{l-m-n} \epm s_{j,m,n} \right)\psi_j^{-1}(y)dydu.
\end{align*}
Denote $$Z=\bpm y_3&t_1&t_2\\ z&x_1'&x_3 \epm\in [\Mat_{(n+j)\times (m+k)}].$$
Then inside the integral $ \phi_{Y_j,\psi_j,N_{n,m}}(h)$, there is an inner integral
$$ \int_{[\Mat_{(n+j)\times (m+k)}]}\phi\left( \bpm I_{n+j} &&Z\\ &1&\\ &&I_{m+k} \epm  g\right)dZ,$$
which, by Lemma \ref{lemma: integral over Z}, is
\begin{align*}
\sum_{\substack{\gamma_1\in N_{n+j}(F)\bs \GL_{n+j}(F)\\ \gamma_2\in N_{m+k}(F)\bs \GL_{m+k}(F)}}W_\phi^\psi\left(\bpm \gamma_1 &&\\ &1&\\ &&\gamma_2 \epm g \right).
\end{align*}
Thus we get
\begin{align*}
\phi_{Y_j,\psi_j,N_{n,m}}(h)=\sum_{\substack{\gamma_1\in N_{n+j}(F)\bs \GL_{n+j}(F)\\ \gamma_2\in N_{m+k}(F)\bs \GL_{m+k}(F)}}\int W_\phi^\psi\left( \bpm \gamma_1 &&\\ &1&\\ &&\gamma_2 \epm \bspm I_n&0&0&0&0\\ x_1&v_1&x_2&0&0\\ &&1&y_2&0\\ &&&v_2&0\\ &&&y_1&I_{m}\\  \espm  \gamma_{n,m}h s_{j,m,n}  \right)
\end{align*}
where $y_1=\bpm y_3' \\ y_1' \epm \in [\Mat_{m\times k}],$ and $h=\bpm h&\\ &I_{l-m-n} \epm$. Note that the above formula is similar to \eqref{eq: phi Y psi N 2}. By the same method as in the proof of Theorem \ref{theorem: unfolding of the left side}, we get that 
\begin{align*}
\phi_{Y_j,\psi_j,N_{n,m}}(h)&=\sum_{\substack{\gamma_1\in N_n(F)\bs \GL_n(F)\\ \gamma_2\in N_m(F)\bs \GL_m(F)}}\int_{\ov V^{j,m,n}(\BA)} W_\phi^\psi \left( \bspm \gamma_1&&\\ &I_{l-m-n}&\\ &&\gamma_2 \espm \ov v  \gamma_{n,m}h s_{j,m,n}  \right)d\ov v\\
&=\sum_{\substack{\gamma_1\in N_n(F)\bs \GL_n(F)\\ \gamma_2\in N_m(F)\bs \GL_m(F)}}\int_{\ov V^{j,m,n}(\BA)} W_\phi^\psi \left(  \ov v  \gamma_{n,m}\bspm\gamma_1&&\\ &\gamma_2&\\ &&I_{l-m-n} \espm h s_{j,m,n}  \right)d\ov v.
\end{align*}
Plugging the above equation into \eqref{eq: I(phi,wt f)}, we get that
\begin{align*}
I(\phi,\wt f_\Bs)&=\int_{\{(N_n(F)\times N_m(F))N_{n,m}(\BA)\}\bs \GL_{n+m}(\BA)}\int_{\ov V^{j,m,n}(\BA)}W_\phi^\psi(\bar v \gamma_{n,m} hs_{j,m,n})\wt f_\Bs(h)d\ov v dh\\
&=\int_{N_{n+m}(\BA)\bs \GL_{n+m}(\BA)}\int_{\ov V^{j,m,n}(\BA)}W_\phi^\psi(\ov v \gamma_{n,m} h s_{j,m,n})\\
&\quad\quad\cdot \int_{[N_n]\times [N_m]}\wt f_s\left(\bpm u_1 &\\ &u_2 \epm \right)\psi(u_1)\psi(u_2)du_1du_2 d\ov v dh\\
&=\int_{N_{n+m}(\BA)\bs \GL_{n+m}(\BA)}\int_{\ov V^{j,m,n}(\BA)}W_\phi^\psi(\ov v \gamma_{n,m} h s_{j,m,n})\xi^{\psi^{-1}}_{\wt f_\Bs}(h)d\ov v dh.
\end{align*}
The result follows.
\end{proof}

\subsection{Convergence and justifications}
\label{subsection-convergence-and-justifications}
In this subsection, we prove the convergence of the double integral in \eqref{eq-unfolding-eq1}, for $\Re(s_1)\gg 0, \Re(s_2)\gg 0$. This is standard and similar to many other situations like \cite{BAS: Uodd}.  In the following, we will only consider the number field case since the proof in the function field case is easier. Using the Iwasawa decomposition and the fact that $(N_m(F)\times N_n(F))\backslash (N_m(\BA)\times N_n(\BA))$ is compact, the convergence of the double integral in \eqref{eq-unfolding-eq1}, for $\Re(s_1)\gg 0, \Re(s_2)\gg 0$, quickly reduces to the convergence of
\begin{equation*}
\int_{T_{m+n}(\BA)}\|t\|^{N_0}  |\det(a)|^{\Re(s_1)+c_1}|\det(b)|^{-\Re(s_2)+c_2}
\int_{\ov U^{j,m,n}(\BA)} |W_\phi^\psi\left(\bt_{m,n}(a, b) \ov u  \right)| d\ov udt
\end{equation*}
where $t=\diag(a, b)$ with $a\in T_m(\BA), b\in T_n(\BA)$ and $\bt_{m,n}(a,b)=\eta_j \iota_j( t) \eta_j^{-1}=\diag(a,I_{l-m-n},b)$. Here, $N_0$, $c_1$ and $c_2$ are fixed given positive numbers. Note that the integration over $K_{\GL_l(\BA)}$ is dropped, using reasoning similar to that in \cite{BAS:Uodd}*{Remark 4.7} (in conjunction with Lemma~\ref{lemma-archimedean-gauge-estimate} and \eqref{eq-convergence-gauge-arch}).
We may assume that the Whittaker function $W_\phi^\psi$ decomposes as $\prod_{v}W_v$, where $W_v$ is a local $\psi_v$-Whittaker function of $\pi_v$, such that outside of a finite set $S$ of places (including the archimedean ones), $\pi_v$ is unramified, and $W_v=W_v^0$ is the normalized unramified $\psi_v$-Whittaker function of $\pi_v$ whose value at the identity is equal to 1. At a place $v\in S$, we assume that $W_v(I_{l})=1$. We assume that for $v$ outside of $S$, $\psi_v$ is unramified. It suffices to prove, for $\Re(s_1)\gg 0, \Re(s_2)\gg 0$, that we have
\begin{equation}
\prod_{v} \int_{T_{m+n}(F_v)}\|t\|^{N_0}_v  |\det(a)|_v^{\Re(s_1)+c_1}|\det(b)|_v^{-\Re(s_2)+c_2}
\int_{\ov U^{j,m,n}(F_v)} |W_v\left(\bt_{m,n}(a,b) \ov u  \right)| d\ov udt<\infty.
\label{eq-convergence-eq2}
\end{equation}

\begin{lemma}\label{lemma-convergence-support-Ujmn}
Let $v$ be a finite place of $F$.
For fixed $W_v\in \CW(\pi_v,\psi_v)$, and $t=\diag(a,b)$ with $a\in T_m(F_v), b\in T_n(F_v)$, the function
$$\ov u\mapsto W_v(\bt_{m,n}(a,b)\ov u ), \quad \ov u\in \ov U^{j,m,n}(F_v)$$
has compact support in $\ov U^{j,m,n}(F_v)$, which is independent of $t$. If $W_v=W_v^0$, then this support is in $\ov U^{j,m,n}(\CO_v)$.
\end{lemma}
\begin{proof}
The proof is a standard ``root killing" argument, and it is similar to the proof of \cite{BAS: Uodd}*{Lemma 4.1}. We omit the details.
\end{proof}

Suppose $v$ is finite. By Lemma~\ref{lemma-convergence-support-Ujmn}, to prove that the local integral of \eqref{eq-convergence-eq2} at the place $v$ converges, it suffices to show
\begin{equation}
\int_{T_{m+n}(F_v)}\|t\|^{N_0}_v  |\det(a)|_v^{\Re(s_1)+c_1}|\det(b)|_v^{-\Re(s_2)+c_2} |W_v\left(\bt_{m,n}(a,b)   \right)|  dt<\infty
\label{eq-convergence-eq3}
\end{equation}
for $\Re(s_i)\gg 0$. Now we recall gauge estimates on Whittaker functions in \cite{JPSS-GL3I}*{Section 2}. A gauge on $\GL_l(F_v)$ is a function $\xi$ on $\GL_l(F_v)$ which is invariant on the left under $N_l(F_v)$, on the right under $\GL_l(\CO_v)$, and which on $T_l(F_v)$ has the form
\begin{equation}
\label{eq-convergence-gauge1}
\xi(t)=|t_1t_2\cdots t_{l-1}|_v^{-c}\Phi(t_1, t_2, \cdots, t_{l-1})
\end{equation}
for 
\begin{equation*}
t=\diag(t_1t_2\cdots t_l, t_2\cdots t_l, \cdots, t_{l-1}t_l, t_l)\in T_l(F_v),
\end{equation*}
where $c\ge 0$ is a real number and $\Phi\ge 0$ is a Schwartz-Bruhat function on $F_v^{l-1}$. In particular, $\xi$ is invariant under the center of $\GL_l(F_v)$.
Write $a\in T_m(F_v)$ and $b\in T_n(F_v)$ as
\begin{equation}
\begin{split}
\label{eq-convergence-abcoordinates}
a&=\diag(a_1\cdots a_m, a_2\cdots a_m, \dots, a_{m-1}a_m, a_m),\\
 b&=\diag(b_1^{-1}, b_1^{-1}b_2^{-1}, \dots , b_1^{-1}b_2^{-1}\cdots b_{n}^{-1}),
\end{split}
\end{equation}
with $a_i\in F_v^\times, b_j\in F_v^\times$. Then
\begin{equation*}
\begin{split}
|\det(a)|_v&=|a_1a_2^2\cdots a_{m-1}^{m-1}a_m^m|_v, \\
|\det(b)|_v&=|b_1^nb_2^{n-1}\cdots b_{n-1}^2 b_n|_v^{-1},
\end{split}
\end{equation*}
and
\begin{align*}
\bt_{m,n}&(\diag(a,b))=\\
&\diag(a_1\cdots a_m, a_2\cdots a_m, \dots, a_{m-1}a_m, a_m, 1, 1, \dots, 1, b_1^{-1}, b_1^{-1}b_2^{-1}, \dots , b_1^{-1}b_2^{-1}\cdots b_{n}^{-1}).
\end{align*}
Then for a gauge $\xi$ on $\GL_l(F_v)$, it follows from \eqref{eq-convergence-gauge1} that there is some real number $c\ge 0$ and a Schwartz-Bruhat function $\Phi$ on $F_v^{l-1}$ such that
\begin{equation}
\label{eq-convergence-gauge2}
\xi(\bt_{m,n}(\diag(a,b)))= |a_1\cdots a_m b_1\cdots b_n|_v^{-c} \Phi(a_1, a_2, \dots, a_m, 1, \dots, 1, b_1, b_2, \dots, b_n).
\end{equation}
Write $|\omega_\pi |_v=\alpha^{c_0}$, where $\alpha$ is a non-negative real-valued function on $F_v^\times$ and $c_0$ is a real number. By  \cite{JPSS-GL3I}*{Proposition 2.3.6}, for any Whittaker function $W_v\in \mathcal{W}(\pi_v,\psi_v)$, there is a gauge $\xi$ such that 
\begin{equation}
\label{eq-convergence-eq4}
|W_v\otimes \alpha^{-c_0/l}|\le \xi.
\end{equation}
Then \eqref{eq-convergence-eq3} follows from \eqref{eq-convergence-gauge2} and the estimate \eqref{eq-convergence-eq4}. This proves that the local integral from \eqref{eq-convergence-eq2} over a finite place $v$ is convergent. 

Now we turn to the archimedean places. Let $v$ be an archimedean place, so $F_v$ is either $\mathbb{R}$ or $\mathbb{C}$. 
We recall the notion of gauge \cite{JacquetShalika-archimedean} in this setting, which is slightly different from the non-archimedean case. Let $\chi$ be a sum of positive characters of $T_l(F_v)$ trivial on the center of $\GL_l(F_v)$. An homogeneous gauge on $\GL_l(F_v)$ is a function $\xi$ on $\GL_l(F_v)$ of the form
\begin{equation*}
\xi(ntk)=\chi(t)\Phi(t_1, t_2,\cdots, t_{l-1}),
\end{equation*}
where $n\in N_l(F_v)$, $t=\diag(t_1, \cdots, t_l)\in T_l(F_v)$, $k$ is in the maximal compact subgroup $K_l$ of $\GL_l(F_v)$, and $\Phi>0$ is a rapidly decreasing function in $l-1$ variables. Here, $\Phi$ being rapidly decreasing means that, for every set of integers $N_i$, $1\le i \le l-1$, there is a constant $C>0$ such that
\begin{equation*}
\Phi(t_1, t_2, \cdots, t_{l-1})\le C\prod_i (1+ |t_i|_v^2)^{-N_i}.
\end{equation*}
We have the following estimate.

\begin{lemma}
\label{lemma-archimedean-gauge-estimate}
Let $v$ be an archimedean place. Let $\xi$ be a homogeneous gauge on $\GL_l(F_v)$. Let $a\in T_m(F_v)$ and $b\in T_n(F_v)$, with Iwasawa decompositions
\begin{equation*}
\begin{split}
a=n_1 \diag(t_1, \cdots, t_{m})  k_1,\quad b =n_2 \diag(t_{l-n+1}, \cdots, t_{l}) k_2,
\end{split}
\end{equation*}
where $t_1, \cdots, t_{m}, t_{l-n+1}, \cdots, t_l$ are positive real numbers. Set $t_{m+1}=t_{m+2}=\cdots =t_{l-n}=1$. Given positive integers $M_1, \cdots, M_j, N_1, \cdots, N_n$, $L_1, \cdots, L_{l-1}$, there exists a positive constant $C>0$ such that 
\begin{equation*}
\begin{split}
&   \xi \left(  \begin{pmatrix} a&&&&\\ &I_j&&&\\&&1&&\\&&&I_k&\\&&&&b \end{pmatrix}  \begin{pmatrix} I_m&&&&\\ x&I_j&&&\\&&1&&\\&&&I_k&\\&&&y&I_n \end{pmatrix}\right) \\
 \le& 
C  \prod_{i=1}^j  (1+||x_{i}||^2)^{-M_{i}} \prod_{i=1}^n  (1+||y_{i}||^2)^{-N_{i}}  \prod_{i=1}^{l-1} \left(1+\left| \frac{t_i}{t_{i+1}}\right|^2\right)^{-L_i}\chi(t_1, t_2, \cdots, t_l).
\end{split}
\end{equation*}
Here, $\chi$ is a fixed sum of positive characters of $T_l(F_v)$.
\end{lemma}
\begin{proof}
The proof is similar to that of \cite{JacquetShalika-archimedean}*{Lemma 5.2}. See also \cite{BAS:Uodd}*{Lemma 4.6}. We omit the details.
\end{proof}

By \cite{JacquetShalika-archimedean}*{Proposition 2.1}, for any Whittaker function $W_v\in \mathcal{W}(\pi_v, \psi_v)$, there is a gauge $\xi$ on $\GL_l(F_v)$ such that 
\begin{equation}
|W_v(g)|\le \xi(g), \quad g\in \GL_l(F_v).
\label{eq-convergence-gauge-arch}
\end{equation}
Combining \eqref{eq-convergence-gauge-arch} with Lemma~\ref{lemma-archimedean-gauge-estimate}, we conclude that the archimedean part of the product in \eqref{eq-convergence-eq2} is convergent. 

Now we are ready to prove \eqref{eq-convergence-eq2}. By Lemma~\ref{lemma-convergence-support-Ujmn} and due to the convergence of the local integral at each local place, it suffices to prove 
\begin{equation}
\prod_{v\not\in S} \int_{T_{m+n}(F_v)}\|t\|^{N_0}_v  |\det(a)|_v^{\Re(s_1)+c_1}|\det(b)|_v^{-\Re(s_2)+c_2}
|W_v^0\left(\bt_{m,n}(a,b)  \right)| dt<\infty.
\label{eq-convergence-eq5}
\end{equation}
Write $a, b$ as in \eqref{eq-convergence-abcoordinates}. Similar to \cite{BAS: Uodd}*{Lemma 4.3}, in order for the function
\begin{equation*}
   \bt_{m,n}(a,b) \mapsto  W_v^0\left(\bt_{m,n}(a,b)\right)
\end{equation*}
to be non-zero, we must have $a_1, \cdots, a_m, b_1, \cdots, b_n\in \mathcal{O}_{v}$. Thus, 
\begin{equation*}
    W_v^0\left(\bt_{m,n}(a,b)\right) = \Phi^0(a_1, \cdots, a_m, b_1, \cdots, b_n) W_v^0\left(\bt_{m,n}(a,b)\right),
\end{equation*}
where $\Phi^0$ is the characteristic function of $\mathcal{O}_{v}^{m+n}$. On the other hand, since $W_\phi^\psi=\prod_{v\in S} W_v \cdot \prod_{v\not\in S} W_v^0$ is the Whittaker function of a cusp form $\phi$, there are constants $B_0, C_0>0$, such that, for all $h\in \GL_l(\BA)$, 
\begin{equation*}
    | W_\phi^\psi(h)|\le C_0 \| h\|^{B_0}.
\end{equation*}
By \cite{MW}*{\S 1.2.2 (ii)}, there is a constant $C_1$ (which depends on $C_0$), such that,
for $t=\diag(a,b)\in \prod_{v\not\in S} T_{m+n}(\BA^S)$ with $a, b$ as in \eqref{eq-convergence-abcoordinates}, we have
\begin{equation*}
    \left| \prod_{v\not\in S} W_v^0 \left(\bt_{m,n}(a,b)\right) \right| \le C_1 \| a\|^{B_0} \| b\|^{B_0}.
\end{equation*}
Thus, the product of local integrals over $v\not\in S$ in \eqref{eq-convergence-eq5} is bounded by
\begin{equation*}
\begin{split}
C_1 & \prod_{v\not\in S} \int_{F_v^{m+n}} \Phi^0(a_1, \cdots, a_m, b_1, \cdots, b_n)    \|a\|^{N_0+B_0} \|b\|^{N_0+B_0}  \\
& \cdot |a_1a_2^2\cdots a_{m-1}^{m-1}a_m^m|_v^{\Re(s_1)+c_1} |b_1^nb_2^{n-1}\cdots b_{n-1}^2 b_n|_v^{\Re(s_2)-c_2}  \prod_{1\le i\le m} da_i \prod_{1\le j \le n} db_j.
\end{split}
\end{equation*}
This product converges for $\Re(s_1)\gg 0, \Re(s_2)\gg 0$ as desired (see, for example, \cite{BAS: Uodd}*{Lemma 4.5}). This finishes the proof of \eqref{eq-convergence-eq5}, and hence the convergence of the double integral in \eqref{eq-unfolding-eq1}, for $\Re(s_1)\gg 0, \Re(s_2)\gg 0$.

\section{The local integrals}\label{section: local theory}
In this section, let $F$ be a local field. Let $\psi$ be a nontrivial additive character of $F$. We still fix a positive integer $l$ and non-negative integers $m,n$ such that $m+n\le l-1$. For $0\le j\le l-m-n-1$, we set $k=l-m-n-1-j$. 
\subsection{Definition of the local zeta integrals}\label{subsection: defn of local zeta integral}
 Let $\pi$ be an irreducible generic representation of $\GL_l(F)$ and let $\CW(\pi,\psi)$ be its Whittaker model. Let $(\tau_1,V_{\tau_1})$ (resp. $(\tau_2,V_{\tau_2})$) be an irreducible generic representation of $\GL_m(F)$ (resp. $\GL_n(F)$). As in the last section, we write $\btau=(\tau_1,\tau_2)$ and $\wh \btau=(\tau_2,\tau_1)$. Let $\Bs=(s_1,s_2)\in \C^2$ be a pair of complex numbers. Then we can consider the induced representation 
$$\RI(\Bs,\btau):=\Ind_{P_{m,n}(F)}^{\GL_{m+n}(F)}(\tau_1||^{s_1-\frac{1}{2}}\boxtimes \tau_2||^{-s_2+1/2}).$$
We fix $\psi^{-1}$-Whittaker functionals $\lambda_{i}$ of $\tau_{i}$. Recall that a section $f_{\Bs}\in \RI(\Bs,\btau)$ is a function $f_\Bs: \GL_{m+n}(F)\to V_{\tau_1}\boxtimes V_{\tau_2}$ satisfying certain quasi-invariance properties. We consider the $\C$-valued function
$$\xi_{f_\Bs}:\GL_{m+n}(F)\times \GL_m(F)\times \GL_n(F)\to \C$$
defined by 
$$\xi_{f_\Bs}(h,a_1,a_2)=\lambda_{1}\boxtimes \lambda_{2}(\tau_{1}(a_1)\boxtimes \tau_{2}(a_2)(f_\Bs(h))).$$
Set $\CW(\Bs,\btau,\psi^{-1})=\wpair{\xi_{f_\Bs}: f_\Bs\in \RI(\Bs,\btau)}$. Note that an element $\xi_\Bs$ satisfies
\begin{align*}
\xi_{\Bs}\left(\bpm b_1 &\\ &b_2 \epm u h, a_1,a_2 \right)=|\det(b_1)|^{s_1+\frac{n-1}{2}} |\det(b_2)|^{-s_2-\frac{m-1}{2}}\xi_{\Bs}(h,a_1b_1,a_2b_2),
\end{align*}
for $a_1,b_1\in \GL_m(F), a_2,b_2\in \GL_n(F), u\in N_{m,n}(F), h\in \GL_{m+n}(F)$. In particular
\begin{align*}
\xi_\Bs\left(\bpm u_1&\\ &u_2 \epm uh, I_m,I_n \right)&=\psi^{-1}(u_1)\psi^{-1}(u_2)\xi_s(h,I_m,I_n),
\end{align*}
for $ u_1\in N_m(F), u_2\in N_n(F), u\in N_{m,n}(F), h\in \GL_{m+n}(F).$ We usually write $\xi_{\Bs}(h,I_m,I_n)$ as $\xi_{\Bs}(h)$ for simplicity.

 Similarly, we can consider the space $\CW(1-\wh\Bs, \wh \btau,\psi^{-1})=\wpair{\xi_{f_{1-\wh \Bs}}: f_{1-\wh \Bs}\in \RI(1-\wh \Bs,\wh\btau)}$. Let $M_{w_{m,n}}: \RI(\Bs,\btau)\to \RI(1-\wh \Bs,\wh\btau)$ be the standard intertwining operator defined by 
 $$M_{w_{m,n}}f_{\Bs}(h)=\int_{N_{n,m}(F)}f_{\Bs}(w_{m,n}uh)du.$$
This integral converges in a cone of $\Bs$ and can be extended meromorphically to $\C^2$, see \cite{MW}. Moreover, to obtain an element in $I(1-\wh \Bs, \wh \btau)$, we need to identify $\tau_1\boxtimes \tau_2$ with $\tau_2\boxtimes \tau_1$ in the definition of $M_{w_{m,n}}$. Equivalently, we need to compose the above integral with the map $\tau_1\boxtimes \tau_2\to \tau_2\boxtimes \tau_1$ defined by $v_1\boxtimes v_2\to v_2\boxtimes v_1$, for $v_1\in \tau_1, v_2\in \tau_2$. We adopt the convention that omits the above identification from the notation. The above intertwining operator gives a map 
$$M_{w_{m,n}}: \CW(\Bs,\btau,\psi^{-1})\to \CW(1-\wh \Bs,\wh\btau,\psi^{-1})$$
defined by 
$$M_{w_{m,n}}(\xi_\Bs)(h,a_1,a_2)=\int_{N_{n,m}(F)}\xi_\Bs (w_{m,n}uh, a_2,  a_1)du.$$

For $W\in \CW(\pi,\psi)$,  $\xi_\Bs\in \CW(\Bs,\btau,\psi^{-1})$, and for $j$ with $0\le j\le l-m-n-1$, we consider the local zeta integrals
\begin{equation}\label{eq: local zeta integral}
\Psi(W, \xi_\Bs; j):=\int_{N_{m+n}(F)\bs \GL_{m+n}(F)}\int_{\ov U^{j,m,n}(F)}W\left(\ov u \gamma_{m,n}\bpm h&\\ &I_{l-m-n} \epm \right)\xi_\Bs(h)d\ov udh,\end{equation}
where we recall that
\begin{align*}
\ov{U}^{j,m,n}&=\wpair{\ov u(x,y)=\bpm I_m&&&&\\ x&I_{j}&&&\\ &&1&&\\ &&&I_{k}&\\ &&&y& I_n  \epm: \substack{ x\in \Mat_{j\times m}\\ y\in \Mat_{n\times k}}},\\
\gamma_{m,n}&=\eta_{j,m,n}s_{j,m,n}^{-1}=\bpm I_m&&\\ &&I_{l-m-n}\\ &I_n& \epm.
\end{align*}
Here we remark that the natural numbers $m,n$ appearing in the local zeta integral \eqref{eq: local zeta integral} are determined by the section $\xi_{\Bs}$, which is an element of $\Ind_{P_{m,n}(F)}^{\GL_{m+n}(F)}(\tau_1||^{s_1-1/2}\otimes \tau_2||^{-s_2+1/2})$. In particular, if we take $\wt \xi_{1-\wh\Bs}\in \CW(1-\wh\Bs,\wh\btau,\psi^{-1})$, we should have
\begin{equation}
\label{eq: right side local zeta integral}
\Psi(W,\wt\xi_{1-\wh\Bs};j)=\int_{N_{m+n}(F)\bs \GL_{m+n}(F)}\int_{\ov U^{j,n,m}(F)}W\left(\ov u \gamma_{n,m}\bpm h&\\ &I_{l-m-n} \epm \right)\wt\xi_{1-\wh\Bs}(h)d\ov udh.
\end{equation}

\begin{remark} In this remark, we assume that $F$ is a global field. If $\phi=\otimes \phi_v$ is a cusp form on $\GL_l(\BA)$ and $f_{\Bs}=\otimes f_{\Bs,v}\in \RI(\Bs,\btau)$ is a pure tensor of a global section, then Theorem \ref{theorem: unfolding of the left side} and Theorem \ref{theorem: unfolding of the right side}  imply that
$$I_j(\phi,f_\Bs)=\prod_v \Psi(\rho(s_{j,m,n})W_v,\xi_{f_{\Bs,v}};j), \quad I_j(\phi,\wt f_\Bs)=\prod_v \Psi(\rho(s_{j,m,n})W_v,\xi_{\wt f_{\Bs,v}};j).$$
Here $\rho$ denotes the right translation.
\end{remark}

\begin{remark}\label{remark: comparison with JPSS integral1} 
In this remark, we consider the degenerate case when $m>0$ and $n=0$. In this case, $\btau=\tau_1$ is just a representation of $\GL_m(F)$, and $\Bs=s$ is a single complex number. Moreover, an element $\xi_\Bs$ has the form $\xi_\Bs(h)=W'(h)|\det(h)|^{s-1/2}$ and we have $M_{w_{m,0}}(\xi_\Bs)=\xi_{\Bs}.$ Thus 
\begin{align*}
\Psi(W,\xi_{\Bs};j)=\int_{N_m(F)\bs \GL_m(F)}\int_{ \Mat_{j\times m}(F)}&W\left(\bpm I_m &&\\ x&I_j&\\ &&I_{l-m-j} \epm \bpm h&\\ &I_{l-m} \epm \right)\\ & \cdot W'(h)|\det(h)|^{s-1/2}dxdh,
\end{align*}
and 
\begin{align*}
\Psi(W,M_{w_{m,0}}(\xi_\Bs);j)=\int_{N_m(F)\bs \GL_m(F)}\int_{ \Mat_{m\times k}(F)}&W\left(\bpm I_{j+1} &&\\ &I_k&\\ &y&I_{m} \epm \bpm &I_{l-m}\\ I_m& \epm \bpm h&\\ &I_{l-m} \epm \right)\\ & \cdot W'(h)|\det(h)|^{s-1/2}dydh.
\end{align*}
Here we notice that $\gamma_{m,0}=I_l$ while $\gamma_{0,m}=\bpm &I_{l-m} \\ I_m & \epm.$ A simple change of variable shows that 
\begin{align*}
\Psi(W,\xi_\Bs;j)=\int_{N_m(F)\bs \GL_m(F)}\int_{ \Mat_{j\times m}(F)}&W\left(\bpm h &&\\ x&I_j&\\ &&I_{l-m-j} \epm \right)  W'(h)|\det(h)|^{s-1/2-j}dxdh.
\end{align*}
One can compare the above integral with that defined by Jacquet--Piatetski-Shapiro--Shalika in \cite{JPSS} and observe that
\begin{equation}\label{eq: comparison with JPSS left side integral}
\Psi(W,\xi_{\Bs};j)=\Psi^{\mathrm{JPSS}}(s-j+\frac{l-m-1}{2}, W,W'; j),
\end{equation}
where $\Psi^{\mathrm{JPSS}}$ denotes the integral defined in \cite{JPSS}*{p.387}. On the other hand, for $W\in \CW(\pi,\psi)$, we denote $\wt W(g)=W(J_l {}^t\!g^{-1})$, which represents a Whittaker function of the contragredient representation $\wt \pi$ of $\pi$. It is easy to check that 
\begin{align*}
\Psi(W,M_{w_{m,0}}(\xi_\Bs);j)=\int_{N_m(F)\bs \GL_m(F)}\int_{ \Mat_{k\times m}(F)}&\wt W\left(\bpm h &&\\ y&I_k&\\ &&I_{j+1} \epm \bpm I_m &\\ & J_{l-m} \epm \right)\\
&\wt W'(h )|\det(h)|^{-s+1/2-k}dydh.
\end{align*}
Thus we get 
\begin{equation*}
\Psi(W, M_{w_{m,0}}(\xi_{\Bs}); j)= \Psi^{\mathrm{JPSS}}\left(1-(s-j+\frac{l-m-1}{2}) ,\rho\left( \bpm I_m &\\  & J_{l-m}\epm \right)\wt W, \wt W'; l-m-1 \right).
\end{equation*}
\end{remark}

\begin{remark}\label{remark: comparison with JPSS integral2} 
Similarly, in the degenerate case where $m=0$ and $n>0$, $\btau=\tau_2$ is just a representation of $\GL_n(F)$, $\Bs=s$ is a single complex number, and an element $\xi_{\Bs}\in \CW(\Bs,\btau,\psi^{-1})$ has the form $\xi_{\Bs}(h)=W^{\prime\prime}(h)|\det(h)|^{-s+1/2}$ where $W^{\prime\prime}\in \mathcal{W}(\tau_2,\psi^{-1})$.	 In this case, we have
\begin{equation*}
\Psi(W,\xi_\Bs;j) = \Psi^{\mathrm{JPSS}}(s+j-\frac{l-n-1}{2}, \rho\left(\bpm I_n&\\ &J_{l-n} \epm\right)  \wt W,\wt W^{\prime\prime}; l-n-j-1),
\end{equation*}	
and
\begin{equation*}
\Psi(W, M_{w_{0,n}}(\xi_\Bs);j) = \Psi^{\mathrm{JPSS}}(1-(s+j-\frac{l-n-1}{2}),    W,  W^{\prime\prime}; l-n-j-1).
\end{equation*}	
\end{remark}

\begin{remark}
If $l=2r+1$ and $m=n$ with $1\le m\le r$, then the integral $\Psi(W,\xi_{\Bs};r-m)$  is the local zeta integral of $\RU_{E/F}(2r+1)\times \Res_{E/F}(\GL_m)$ at split places as in \cite{BAS:Uodd}, where $E/F$ is a quadratic extension of global fields.
\end{remark}

\begin{proposition}\label{prop-convergence}
The local zeta integrals $\Psi(W,\xi_{\Bs}; j)$ are absolutely convergent for $\Re(s_i)\gg 0$ for $i=1,2$. Over nonarchimedean local fields, there exist $W$ and $\xi_\Bs$, such that the integral is absolutely convergent and equals 1, for all $\Bs$. Over archimedean fields, for any $\Bs$, there are choices of data $(W^i, \xi_{\Bs}^i)$ such that $\sum_i \Psi(W^i,\xi_{\Bs}^i; j)$ is holomorphic and nonzero in a neighborhood of $\Bs$.
\end{proposition}

\begin{proof}
For $n=0$, this was already proved in \cite{JPSS} over nonarchimedean local fields and in  \cite{JacquetShalika-archimedean} over archimedean fields. Very similar statements can be found in many other places in the literature, for example, \cite{Soudry1993}, \cite{Soudry1995}, \cite{GinzburgRallisSoudry1998}, \cite{BAS:Uodd},  and \cite{CFK2022}. We provide some details here for completeness.

First, we consider the case where $F$ is nonarchimedean. By the Iwasawa decomposition and the fact that smooth vectors are finite under the maximal compact subgroup, we get that $\Psi(W,\xi_{\Bs}; j)$ is a finite sum of integrals of the form
\begin{equation*}
\int_{T_{m+n}(F)}\int_{\ov U^{m,n}(F)}W^\prime(\bt_{m,n}(a,b)\ov u) d\ov u W_{\tau_1}(a)W_{\tau_2}(b)  |\det(a)|^{s_1+\frac{n-1}{2}-j}|\det(b)|^{-s_2-\frac{m-1}{2}+k} \delta_{B_{m+n}}(t)^{-1} dt
\end{equation*}
where $W^\prime\in \mathcal{W}(\pi,\psi)$, $W_{\tau_1}\in \mathcal{W}(\tau_1,\psi^{-1})$, $W_{\tau_2}\in \mathcal{W}(\tau_2, \psi^{-1})$, $t=\diag(a, b)$ with $a\in T_m(F), b\in T_n(F)$ and $\bt_{m,n}(a,b)=\diag(a,I_{l-m-n},b)$. Here the term $|\det(a)|^{-j} |\det(b)|^{k} $ comes from conjugating $\bt_{m,n}(a,b)$ to the left of $\ov u$ and making a change of variables on $\ov u$. By Lemma~\ref{lemma-convergence-support-Ujmn}, the last integral is a finite sum of integrals of the form
\begin{equation}
\int_{T_{m+n}(F)} W^\prime(\bt_{m,n}(a,b) )  W_{\tau_1}(a)W_{\tau_2}(b)  |\det(a)|^{s_1+\frac{n-1}{2}-j}|\det(b)|^{-s_2-\frac{m-1}{2}+k} \delta_{B_{m+n}}(t)^{-1} dt.
\label{eq-convergence-local-integral-finite-place1}
\end{equation}
Now we recall the asymptotic expansion of Whittaker functions \cite{JPSS}*{Section 2.5}. There is a finite set $X_l$ of functions on $T_l(F)$ such that for every $W\in \mathcal{W}(\pi,\psi)$ we have
\begin{equation*}
W(t)=\sum_{\chi\in X_l} \omega_{\pi}(t_l) \phi_\chi(t_1, t_2, \cdots, t_{l-1}) \chi(t)
\end{equation*}
where $t=\diag(t_1t_2\cdots t_l, t_2\cdots t_l, \cdots, t_{l-1}t_l, t_l)\in T_l(F)$ and $\phi_\chi\in \mathcal{S}(F^{l-1})$. Then for every $W\in \mathcal{W}(\pi,\psi)$, we have
\begin{equation}
|W(t)| \le \sum_{\eta\in Y_l}   \phi_\eta(t_1, t_2, \cdots, t_{l-1}) \eta(t)
\label{eq-convergence-asymptotic-expansion-Whittaker}
\end{equation}
where $\phi_\eta\in \mathcal{S}(F^{l-1})$ is non-negative and $\eta$ varies in another finite set $Y_l$ of finite functions on $T_l(F)$.
Applying the majorization \eqref{eq-convergence-asymptotic-expansion-Whittaker} to $W^\prime$ (and the analogous ones for $W_{\tau_1}$ and $W_{\tau_2}$), we obtain the absolute convergence of the integral \eqref{eq-convergence-local-integral-finite-place1}  for $\Re(s_i)\gg 0$ for $i=1,2$. Hence $\Psi(W,\xi_{\Bs}; j)$ is absolutely convergent for $\Re(s_i)\gg 0$ for $i=1,2$.

We continue to assume that $F$ is nonarchimedean. Since $N_{m+n}(F)T_{m+n}(F)\overline{N}_{m+n}(F)$ is an open dense subset of $\GL_{m+n}(F)$ whose complement has Haar measure zero, we may rewrite $\Psi(W,\xi_{\Bs}; j)$ as
\begin{equation} 
\label{eq-convergence-local-integral-Iwasawa}
\begin{split}
\int_{T_{m+n}(F)} \int_{\overline{N}_{m+n}(F)}  \int_{\ov U^{j,m,n}(F)} &W\left(\ov u \gamma_{m,n}\bpm t \overline{v} &\\ &I_{l-m-n} \epm \right) \xi_\Bs( \overline{v}, a, b)  \\
 &|\det(a)|^{s_1+\frac{n-1}{2}}  |\det(b)|^{-s_2-\frac{m-1}{2}} \delta_{B_{m+n}}(t)^{-1} d\ov u d\ov v dt,
\end{split}
\end{equation}
where $t=\diag(a, b)$ with $a\in T_m(F)$, $b\in T_n(F)$.
Similar to \cite{Soudry1993}*{Proposition 6.1}, we choose $\xi_\Bs$ to have support in $B_{m+n}(F)\cdot \mathcal{V}_1$, where $\mathcal{V}_1$ is a small open compact subgroup of $\GL_{m+n}(F)$, and such that $\xi_\Bs(u, b_1, b_2)=W_{\tau_1}(b_1)W_{\tau_2}(b_2)$ for $u\in \mathcal{V}_1$, $b_1\in T_{m}(F), b_2\in T_n(F)$. Here, $W_{\tau_i}\in \mathcal{W}(\tau_i, \psi^{-1})$ for $i=1, 2$. We choose $\mathcal{V}_1$ to be so small that $W$ is fixed by $\pi( \diag(\overline{v}, I_{l-m-n}))$ for $\overline{v}\in \mathcal{V}_1$.
Thus, $\Psi(W,\xi_{\Bs}; j)$ is equal to
\begin{equation*}
\begin{split}
\mathrm{vol}(\mathcal{V}_1\cap \overline{N}_{m+n}(F))  \cdot \int_{T_{m+n}(F)}   \int_{\ov U^{j,m,n}(F)}  &W\left(\ov u \gamma_{m,n}\bpm t  &\\ &I_{l-m-n} \epm \right)  W_{\tau_1}(a)W_{\tau_2}(b) \\
& |\det(a)|^{s_1+\frac{n-1}{2}}  |\det(b)|^{-s_2-\frac{m-1}{2}} \delta_{B_{m+n}}(t)^{-1} d\ov u   dt.
\end{split}
\end{equation*}
We conjugate $\diag(t, I_{l-m-n})$ to the left of $\ov u$ and make a change of variable in $\ov u$ to get
\begin{equation*}
\begin{split}
\mathrm{vol}(\mathcal{V}_1\cap \overline{N}_{m+n}(F))  \cdot \int_{T_{m+n}(F)}   \int_{\ov U^{j,m,n}(F)}  &\rho(\gamma_{m,n}) W\left(  \bpm a  & &\\ &I_{l-m-n}&\\ & & b \epm \ov u  \right)  W_{\tau_1}(a)W_{\tau_2}(b) \\
& |\det(a)|^{s_1+\frac{n-1}{2}-j}  |\det(b)|^{-s_2-\frac{m-1}{2}+k} \delta_{B_{m+n}}(t)^{-1} d\ov u   dt.
\end{split}
\end{equation*}
Now we choose $W$, $W_{\tau_1}$, and $W_{\tau_2}$ such that the function $$(a,b,\ov u)\mapsto \rho(\gamma_{m,n})W\left(  \bpm a  & &\\ &I_{l-m-n}&\\ & & b \epm \ov u  \right) W_{\tau_1}(a)W_{\tau_2}(b) $$ is the characteristic function of a small neighborhood of $(I_m, I_n, I_l)$. Thus, the integral can be made constant.

Now we assume $F$ is archimedean. Similar to \cite{Soudry1993}*{Lemma 5.2}, there is a  positive integer $A_0$, such that for any $\xi_{\Bs}$, there is a constant $c_{\Bs}>0$, such that 
\begin{equation*}
|\xi_{\Bs}(\diag(a, b)k) | \le c_\Bs |\det(a)|^{\Re(s_1)+\frac{n-1}{2}}  |\det(b)|^{-\Re(s_2)-\frac{m-1}{2}}  \|\diag(a, b)\|^{A_0},  
\end{equation*}
where $a\in T_m(F), b\in T_n(F)$, and $k$ is in the maximal compact subgroup $K_l$ of $\GL_l(F)$.  We then use the Iwasawa decomposition, \eqref{eq-convergence-gauge-arch}, and Lemma~\ref{lemma-archimedean-gauge-estimate} to conclude the absolute convergence of $\Psi(W,\xi_{\Bs}; j)$.

Now we prove the non-vanishing of the integrals when $F$ is archimedean. Write $\Psi(W,\xi_{\Bs}; j)$ in the form \eqref{eq-convergence-local-integral-Iwasawa}. Similar to \cite{GinzburgRallisSoudry1998}*{Proposition 6.7}, we can choose $\xi_{\Bs}$ to have support in $P_{m,n}(F)\cdot \ov N_{m+n}(F)$, and assume
\begin{align*}
\xi_{\Bs}\left(\bpm b_1 &\\ &b_2 \epm u \ov v, a_1,a_2 \right)=|\det(b_1)|^{s_1+\frac{n-1}{2}} |\det(b_2)|^{-s_2-\frac{m-1}{2}} \varphi_1(\ov v) W_{\tau_1}(a_1b_1)W_{\tau_2}(a_2b_2),
\end{align*}
for $a_1,b_1\in \GL_m(F), a_2,b_2\in \GL_n(F), u\in N_{m,n}(F), \ov v\in \ov N_{m+n}(F)$, $W_{\tau_i}\in \mathcal{W}(\tau_i,\psi^{-1})$ for $i=1, 2$, and $\varphi_1\in C_c^\infty(\ov N_{m+n}(F))$. With this choice, $\Psi(W,\xi_{\Bs}; j)$ is equal to an integral of the form
\begin{equation}
\label{eq-nonvanishing-local-integral-archimedean}
\begin{split}
\int_{T_{m+n}(F)} \int_{\overline{N}_{m+n}(F)}  \int_{\ov U^{j,m,n}(F)} &W\left(\ov u \gamma_{m,n}\bpm t \overline{v} &\\ &I_{l-m-n} \epm \right) \varphi_1(\overline{v})W_{\tau_1}(a) W_{\tau_2}(b) \\
 &|\det(a)|^{s_1+\frac{n-1}{2}}  |\det(b)|^{-s_2-\frac{m-1}{2}} \delta_{B_{m+n}}(t)^{-1} d\ov u d\ov v dt.
\end{split}
\end{equation}
We consider the $d\ov v$ integration first. This is a convolution of $\varphi_1$ against $W\in \mathcal{W}(\pi,\psi)$. By the Dixmier-Malliavin Theorem \cite{DixmierMalliavin1978}, a linear combination of the $d\ov v$ integrals represents a general element of $\mathcal{W}(\pi,\psi)$ (see \cite{GinzburgRallisSoudry1998}*{pp. 233} or \cite{Soudry1995}*{pp. 219}). Thus, a suitable linear combination of integrals of the form \eqref{eq-nonvanishing-local-integral-archimedean} gives an integral of the form
\begin{equation*}
\begin{split}
\int_{T_{m+n}(F)}   \int_{\ov U^{j,m,n}(F)} &W\left(\ov u \gamma_{m,n}\bpm t   &\\ &I_{l-m-n} \epm \right)  W_{\tau_1}(a) W_{\tau_2}(b) \\
 &|\det(a)|^{s_1+\frac{n-1}{2}}  |\det(b)|^{-s_2-\frac{m-1}{2}} \delta_{B_{m+n}}(t)^{-1} d\ov u   dt.
\end{split}
\end{equation*}
We conjugate $\diag(t, I_{l-m-n})$ to the left of $\ov u$ to get
\begin{equation*}
\begin{split}
\int_{T_{m+n}(F)}   \int_{\ov U^{j,m,n}(F)} &  \rho(\gamma_{m,n})W\left(  \bpm a  & &\\ &I_{l-m-n}&\\ & & b \epm \ov u  \right) W_{\tau_1}(a) W_{\tau_2}(b) \\
 &|\det(a)|^{s_1+\frac{n-1}{2}-j}  |\det(b)|^{-s_2-\frac{m-1}{2}+k} \delta_{B_{m+n}}(t)^{-1} d\ov u   dt.
\end{split}
\end{equation*}
Now we choose $W$ so that $\rho(\gamma_{m,n})W(t\ov u)=\rho(\gamma_{m,n})W(t)\varphi_2(\ov u)$ for $t\in B_l(F)$, $\ov u \in U^{j,m,n}(F)$ and $\varphi_2\in C_c^\infty(U^{j,m,n}(F))$. Then the above integral becomes
\begin{equation*}
\begin{split}
  \int_{\ov U^{j,m,n}(F)} \varphi_2(\ov u)d\ov u  \cdot   \int_{T_{m+n}(F)}  &  \rho(\gamma_{m,n})W\left(  \bpm a  & &\\ &I_{l-m-n}&\\ & & b \epm    \right) W_{\tau_1}(a) W_{\tau_2}(b) \\
 &|\det(a)|^{s_1+\frac{n-1}{2}-j}  |\det(b)|^{-s_2-\frac{m-1}{2}+k} \delta_{B_{m+n}}(t)^{-1}    dt.
\end{split}
\end{equation*}
The $d\ov u$ integral is a nonzero constant for appropriate $\varphi_2$. For appropriate $W, W_{\tau_1}, W_{\tau_2}$, the $dt$ integral is holomorphic and nonzero in a neighborhood of any given $\Bs$. This proves that there is a linear combination of the local integrals $\Psi(W,\xi_{\Bs}; j)$ which is holomorphic and nonzero in a neighborhood of any given $\Bs$.
\end{proof}

\subsection{Local functional equations}
\label{subsection-local-functional-equation}
\begin{proposition}\label{proposition: existence of gamma}
Let $F$ be a non-archimedean local field of characteristic different from 2. There exists a meromorphic function $\Gamma(\Bs,\pi\times\btau,\psi;j)$ such that
$$\Psi(W,M_{w_{m,n}}(\xi_{\Bs});j) =\Gamma(\Bs,\pi\times\btau,\psi;j)\Psi(W,\xi_{\Bs};j), $$
for any $W\in \CW(\pi,\psi)$ and $\xi_{\Bs}\in \CW(\Bs,\btau,\psi^{-1})$.
\end{proposition}

\begin{proof}
Recall that
$$Y_{j,m,n}=\wpair{\bpm u&*&*\\ &I_{m+n+1} &*\\ &&v \epm, u\in N_{j}, v\in N_{k}}$$
and we have defined a character $\psi_{j}$ on $Y_{j,m,n}(F)$, see \eqref{eq-psi}.
One can check that for any $y\in Y_{j,m,n}(F)$, we have
\begin{align}\label{eq: quasi-invariance1}
\begin{split}
\Psi(\rho(y)W,\xi_{\Bs};j)&=\psi_j(y)\Psi(W,\xi_{\Bs};j), 
\end{split}
\end{align}
and for any $h=\bpm a& b\\ c&d\epm \in \GL_{m+n}(F)$ with $a\in \Mat_{m\times m}(F)$,
\begin{align}\label{eq: quasi-invariance2}
\Psi\left(\rho\left(\bspm I_j& & & & \\ &a&&b& \\ &&1&& \\ &c&&d& \\ &&&&I_k \espm \right)W,\rho(h)\xi_{\Bs};j\right)=\Psi(W,\xi_{\Bs};j).
\end{align}
Let 
$$H_{j,m,n}=\wpair{\bspm u&*&*&*&*\\ &a&&b&*\\ &&1&&*\\ &c&&d&*\\ &&&&v \espm, u\in N_j, v\in N_k, \bpm a& b\\ c&d\epm\in \GL_{m+n} }=\GL_{m+n}\ltimes Y_{j,m,n}.$$
One can define a representation $\nu_{\Bs}$ of $H_{j,m,n}(F)$ by $\nu_{\Bs}|_{\GL_{m+n}(F)}=\RI(\Bs,\btau)$ and $\nu_{\Bs}|_{Y_{j,m,n}(F)}=\psi_j$. Then \eqref{eq: quasi-invariance1} and \eqref{eq: quasi-invariance2} imply that the bilinear form $(W,\xi_{\Bs})\mapsto \Psi(W,\xi_{\Bs};j)$ defines an element in 
$$\Hom_{H_{j,m,n}(F)}(\pi\otimes \nu_{\Bs},1).$$
Similarly, the bilinear form $(W,\xi_{\Bs})\mapsto \Psi(W, M(\Bs, \btau)\xi_{\Bs}; j)$ satisfies the same properties \eqref{eq: quasi-invariance1} and \eqref{eq: quasi-invariance2} and hence 
$(W,\xi_{\Bs})\mapsto \Psi(W, M(\Bs, \btau)\xi_{\Bs}; j)$
also defines an element in 
$$\Hom_{H_{j,m,n}(F)}(\pi\otimes \nu_{\Bs},1).$$
By the uniqueness of Bessel models (see \cite{GGP} and \cite{Chan2022}*{Corollary 5.11} when $F$ is of characteristic zero, and \cite{Mezer2023}*{Theorem 1.9} when $F$ is of positive characteristic different from 2), we have  
\begin{equation}\label{eq-uniqueness-of-Bessel}\dim_{\C}\Hom_{H_{j,m,n}(F)}(\pi\otimes \nu_{\Bs},1)\le 1\end{equation}
excluding a discrete set of $\Bs$. This proves the existence of the gamma factor. By Proposition \ref{prop-convergence}, there exist data $W, \xi_\Bs$ such that $\Psi(W, \xi_\Bs; j)$ is 1, which shows that $\Gamma(\Bs,\pi\times \btau,\psi; j)$ is meromorphic.
\end{proof}

\begin{remark}\label{remark: archimedean gamma}
If $F$ is archimedean, the local integrals still define elements in $$ \Hom_{H_{j,m,n}(F)}(\pi\otimes\nu_{\Bs},1).$$
It is known that this Hom space has dimension at most one when $j=0$ by \cite{Chen-Sun} and when $m=n$ by \cites{GGP, JiangSunZhu}. Thus, in these cases, we still have the local gamma factors $\Gamma(\Bs, \pi\times \btau,\psi; j)$. As pointed out in \cite{Chan2022} in the non-archimedean case, the general multiplicity one result (namely, when $m\ne n$) should also follow from the general framework of \cite{GGP}. But it seems that this is not recorded anywhere.
\end{remark}

\begin{remark}
\label{remark-gamma-factor-JPSS}
By Remark~\ref{remark: comparison with JPSS integral1} and Remark~\ref{remark: comparison with JPSS integral2}, we immediately obtain that
\begin{equation*}
\begin{split}
\Gamma((s_1,0), \pi\times(\tau_1, 0), \psi;j) =\omega_{\tau_1}(-1)^{l-1}\gamma( s_1-j+\frac{l-m-1}{2}, \pi\times\tau_1, \psi),
\end{split}
\end{equation*}	
and
\begin{equation*} 
\begin{split}
\Gamma((0, s_2), \pi\times(0, \tau_2), \psi;j) =\omega_{\tau_2}(-1)^{l-1}\gamma(s_2+j-\frac{l-n-1}{2},\wt \pi\times \wt \tau_2,\psi).
\end{split}
\end{equation*}		
\end{remark}

The gamma factor defined in Proposition~\ref{proposition: existence of gamma} is indeed just a product of the JPSS local gamma factors defined in \cite{JPSS}. More precisely, we have the following
\begin{proposition}
\label{proposition-gamma-comparison}
Let $F$ be a local non-archimedean field of characteristic different from 2. Then we have 
\[
\begin{split}	
\Gamma(\Bs, \pi\times(\tau_1, \tau_2), \psi;j) = \omega_{\tau_1}(-1)^{l-1}\omega_{\tau_2}(-1)^{l-1}  \cdot \frac{\gamma(s_1+\frac{k-j}{2},\pi\times\tau_1,\psi) \gamma(s_2+\frac{j-k}{2},\wt \pi\times \wt \tau_2,\psi)}{\gamma(s_1+s_2,\tau_1\times \wt \tau_2,\psi)}.
 \end{split}\]
\end{proposition}
If $l=2r+1$, $m=n$ and $j=r-m$, the gamma factor is just the local gamma factor for $\RU_{E/F}(2r+1)\times \Res_{E/F}(\GL_m)$ at split places, and the above relation with JPSS local gamma factors was proved in \cite{ChengWang}. To streamline the presentation and to avoid making the main body of the paper too long, we will defer the proof of Proposition~\ref{proposition-gamma-comparison} to Appendix~\ref{appendix}.

\subsection{Unramified calculation}

In this subsection, let $F$ be a non-archimedean local field with ring of integers $\CO$. Let $\varpi\in \CO$ be a fixed uniformizer,  and $q=|\CO/(\varpi)|$. Our goal in this subsection is to compute the local zeta integral \eqref{eq: local zeta integral} when everything is unramified. In particular, we assume that $\pi$ is unramified with Satake parameters $\mathbf{\alpha}=\diag(\alpha_1,\dots,\alpha_{l})\in \GL_l(\C)$ and $\tau_1$ (resp. $\wt \tau_2$) is unramified with Satake parameters $\mathbf{\beta}^1=\diag(\beta_1^1,\dots,\beta_m^1)\in \GL_m(\C)$ (resp. $\mathbf{\beta}^2=\diag(\beta_1^2,\dots,\beta_n^2)\in \GL_n(\C)$). Moreover, we assume that $W\in \CW(\pi,\psi)$ is the spherical Whittaker function normalized by $W(I_l)=1$, $\xi_\Bs$ is the spherical Whittaker function associated with the normalized spherical section $f_\Bs\in \RI(\Bs,\btau)$. Let $W_{\tau_i}\in \CW(\tau_i,\psi^{-1})$, $i=1,2$, be the spherical Whittaker function of $\tau_i$ normalized such that $W_{\tau_1}(I_m)=1$ and $W_{\tau_2}(I_n)=1$. From the definition of $\xi_{\Bs}$, we get $$\xi_{\Bs}(t)=|\det(a)|^{s_1+\frac{n-1}{2}}|\det(b)|^{-s_2-\frac{m-1}{2}}W_{\tau_1}(a)W_{\tau_2}(b)$$ 
for $t=\diag(a,b)\in T_{m+n}(F)$, with $a\in T_m(F), b\in T_n(F)$. 
By Iwasawa decomposition $\GL_{m+n}(F)=N_{m+n}(F)T_{m+n}(F)K_{m+n}$, where $K_{m+n}=\GL_{m+n}(\CO)$, we have
\begin{align*}\Psi(W,\xi_\Bs; j)&=\int_{T_{m+n}(F)}\int_{\ov U^{j,m,n}(F)}W(\ov u \gamma_{m,n}\diag(t,I_{l-m-n}))\xi_\Bs(t) \delta_{B_{m+n}}(t)^{-1}d\ov u dt\\
&=\int_{T_{m+n}(F)}\int_{\ov U^{m,n}(F)}W(\bt_{m,n}(a,b)\ov u)\xi_\Bs(t) |\det(a)|^{-j} |\det(b)|^{k}\delta_{B_{m+n}}(t)^{-1}d\ov u dt\\
&=\int_{T_{m+n}(F)}\int_{\ov U^{m,n}(F)}W(\bt_{m,n}(a,b)\ov u)W_{\tau_1}(a)W_{\tau_2}(b)\\
&\quad \cdot |\det(a)|^{s_1+\frac{n-1}{2}-j}|\det(b)|^{-s_2-\frac{m-1}{2}+k} \delta_{B_{m+n}}(t)^{-1}d\ov u dt
\end{align*}
where $t=\diag(a,b)$ with $a\in T_m(F), b\in T_n(F)$ and $\bt_{m,n}(a,b)=\diag(a,I_{l-m-n},b)$.  Here the term $|\det(a)|^{-j} |\det(b)|^{k} $ comes from a modulus character when we change variables on $\ov u$ and the term $\delta_{B_{m+n}}(t)^{-1}$ comes from the corresponding Haar measure when we use the Iwasawa decomposition. 
By Lemma \ref{lemma-convergence-support-Ujmn}, we have
\begin{align}\label{eq: unramified1}
\Psi(W,\xi_\Bs; j)=&\int_{T_{m+n}(F)}W(\bt_{m,n}(a,b))W_{\tau_1}(a)W_{\tau_{2}}(b)\\
& \quad \cdot |\det(a)|^{s_1+\frac{n-1}{2}-j}|\det(b)|^{-s_2-\frac{m-1}{2}+k}\delta_{B_{m+n}}^{-1}\left(\bpm a&\\ &b\epm \right)dadb \nonumber \\
=& \int_{T_{m+n}(F)}W(\bt_{m,n}(a,b^*))W_{\tau_1}(a)W_{\tau_{2}}(b^*) \nonumber \\
& \quad \cdot |\det(a)|^{s_1+\frac{n-1}{2}-j}|\det(b)|^{s_2+\frac{m-1}{2}-k}\delta_{B_{m+n}}^{-1}\left(\bpm a&\\ &b^*\epm \right)dadb  \nonumber
\end{align}
where $b^*=J_n ^t\!b^{-1}J_n^{-1}$, with $ J_n=\bpm &&1\\ &\iddots&\\ 1&& \epm$. Note that the function $b\mapsto W_{\tau_2}(b^*)$ is just the normalized Whittaker function of $\wt \tau_2$, namely, $W_{\tau_2}(b^*)=W_{\wt \tau_2}(b).$
We use the following notations following \cite{Jacquet-Shalika-EulerI}. For $m$-tuple $\mathbf{x}=(x_1,\dots,x_m)$, we write $\varpi^\bx=(\varpi^{x_1},\dots,\varpi^{x_m})$. 
Then \eqref{eq: unramified1} can be written as
\begin{align}\label{eq: unramified2}
\Psi(W,\xi_{\Bs};j)=&\sum_{\bx,\by}W(\varpi^{(\bx,0,\by^*)})W_{\tau_1}(\varpi^{\bx})W_{\wt\tau_2}(\varpi^\by)\delta_{B_{m+n}}^{-1}(\varpi^{(\bx,\by^*)})\\
&\cdot |\det(\varpi^\bx)|^{s_1+\frac{n-1}{2}-j} |\det(\varpi^{\by})|^{s_2+\frac{m-1}{2}-k} . \nonumber
\end{align}
Here $\bx$ (resp. $\by$) runs over all $m$-tuples (resp. $n$-tuples) of all integers, $\by^*=(-y_n,\dots,-y_1)$ for $\by=(y_1,\dots,y_n)$ and $(\bx,0,\by^*) $ denotes the $l$-tuple $(x_1,\dots,x_m,0,\dots,0,-y_n,\dots,-y_1)$ with $l-m-n$ zeros in the middle. Denote $T^+(m)$ the $m$-tuples of integers $\bx=(x_1,\dots,x_m)$ such that $x_1\ge x_2\ge\dots\ge x_m\ge 0$. Similarly, we define $T^+(n)$.
 By the Shintani-Casselman-Shalika formula \cites{Shintani, CasselmanShalika}, we have $W(\varpi^{(\bx,0,\by^*)})=0$ unless $\bx\in T^+(m)$ and $\by\in T^+(n)$. If $\bx\in T^+(m)$ and $\by\in T^+(n)$, we have $W(\varpi^{(\bx,0,\by^*)})=\delta_{B_{l}}^{1/2}(\varpi^{(\bx,0,\by^*)})S_{(\bx,0,\by^*)}(\alpha)$.
 Here $S_{(\bx,0,\by^*)}$ denotes the Schur polynomial associated with $(\bx,0,\by^*)$ (see \cite{Fulton-Harris}*{Appendix A}), or more explicitly, 
$$S_{(\bx,0,\by^*)}(\alpha)=\prod_{1\le i<j\le l}(\alpha_i-\alpha_j)^{-1} \cdot \det\bpm \alpha_1^{x_1+l-1} & \dots & \alpha_{l}^{x_1+l-1}\\ \vdots&&\vdots\\ \alpha_1^{-y_1}&\dots &\alpha_{l}^{-y_1} \epm.$$
Moreover, for $\bx\in T^+(m), \by\in T^+(n)$, we have
$$W_{\tau_1}(\varpi^{\bx})=\delta_{B_{m}}^{1/2}(\varpi^{\bx})S_{\bx}(\beta^1), \quad  W_{\wt \tau_2}(\varpi^{\by})=\delta_{B_n}^{1/2}(\varpi^{\by})S_{\by}(\beta^2).$$
We can check that
\begin{align*}
\delta_{B_m}(\varpi^{\bx})&=|\varpi^{x_1}|^{m-1}|\varpi^{x_2}|^{m-3}\dots |\varpi^{x_m}|^{-m+1},\\
\delta_{B_n}(\varpi^{\by})&=|\varpi^{y_1}|^{n-1}|\varpi^{y_2}|^{n-3}\dots |\varpi^{y_n}|^{-n+1},\\
\delta_{B_{m+n}}(\varpi^{(\bx,\by^*)})&=|\varpi^{x_1}|^{m-1}|\varpi^{x_2}|^{m-3}\dots|\varpi^{x_m}|^{-m+1}|\varpi^{y_n}|^{-n+1}\dots |\varpi^{y_1}|^{n-1}|\det(\varpi^{\bx})|^n|\det(\varpi^{\by})|^m,\\
\delta_{B_{l}}(\varpi^{(\bx,0,\by^*)})&=|\varpi^{x_1}|^{m-1}\dots |\varpi^{x_m}|^{1-m}|\varpi^{y_n}|^{1-n}\dots |\varpi^{y_1}|^{n-1}|\det(\varpi^{\bx})|^{l-m}|\det(\varpi^{\by})|^{l-n}.
\end{align*}
Combining the above formulas, \eqref{eq: unramified2} becomes
\begin{align} \label{eq: unramified3}
\Psi(W,\xi_{\Bs};j)&=\sum_{\substack{\bx\in T^+(m)\\ \by\in T^+(n)}}S_{(\bx,0,\by^*)}(\alpha)S_{\bx}(\beta^1)S_{\by}(\beta^2)|\det(\varpi^{\bx})|^{s_1+\frac{k-j}{2}}|\det(\varpi^{\by})|^{s_2+\frac{j-k}{2}}.
\end{align}

\begin{proposition}
\label{proposition-unramified-computation}
We keep the notations as above. Then
$$\Psi(W,\xi_\Bs; j)=\frac{L(s_1+\frac{k-j}{2},\pi\times \tau_1)L(s_2-\frac{k-j}{2},\wt \pi\times \wt\tau_2)}{L(s_1+s_2,\tau_1\times \wt \tau_2)}.$$
Recall that $k=l-m-n-1-j$ and thus $k-j=l-m-n-1-2j.$
\end{proposition}

If $n=0$, the above formula is the unramified calculation of the Jacquet--Piatetski-Shapiro--Shalika integral, see \cite{Jacquet-Shalika-EulerI}*{Proposition 2.4} and also \cites{Cogdell:Fields, Cogdell:IAS}. If $l=2r+1, m=n$ and $j=r-m=k$, the above unramified calculation is done  in \cite{GPS} (when $r=1$), in \cite{Tamir} (for general $r$ when $m=n=r$) with slightly different normalization, and in \cite{BAS:Uodd} (when $m=n<r$), where this was the unramified calculation of $L$-functions for $\RU_{2r+1,E/F}\times \Res_{E/F}(\GL_r)$ at split places for a quadratic extension $E/F$. 
\begin{proof}
Without loss of generality, we assume that $m\ge n$. Write $T_1=q^{-(s_1+\frac{k-j}{2})}, T_2=q^{-(s_2+\frac{j-k}{2})}$. For an $m$-tuple $\bx=(x_1,\dots,x_m)$, denote $|\bx|=\sum_{i=1}^m x_i$. An $m$-tuple $\bx\in T^+(m)$ can be identified with a partition of $|\bx|$ and can be represented by an Young diagram, see \cite{Fulton-Harris}*{\S 4} for example.  We can then write \eqref{eq: unramified3} as
\begin{equation} \label{eq: unramified4}
\Psi(W,\xi_{\Bs};j)=\sum_{\substack{\bx\in T^+(m)\\ \by\in T^+(n)}}S_{(\bx,0,\by^*)}(\alpha)S_{\bx}(\beta^1)S_{\by}(\beta^2)T_1^{|\bx|}T_2^{|\by|}.
\end{equation}
 On the other hand, we have
$$L(s_1+s_2,\tau_1\times \wt \tau_2)=\det(I-\beta^1\otimes \beta^2 T_1T_2)^{-1}=\sum_{e\ge 0}\Tr(\Sym^e(\beta^1\otimes \beta^2))(T_1T_2)^e.$$
Thus we get that 
\begin{align}\label{eq: unramified5}
\begin{split}
L(s_1+s_2,\tau_1\times \wt \tau_2)\Psi(W,\xi_{\Bs};j)=\sum_{\substack{\bx\in T^+(m), \by\in T^+(n), e\ge 0} } &S_{(\bx,0,\by^*)}(\alpha)S_{\bx}(\beta^1)S_{\by}(\beta^2)\\
&\cdot \Tr(\Sym^e(\beta^1\otimes \beta^2))T_1^{|\bx|+e}T_2^{|\by|+e}.
\end{split}
\end{align}
Since
$$L(s_1+\frac{k-j}{2},\pi\times \tau_1)=\sum_{c\ge 0}\Tr(\Sym^c(\alpha\otimes \beta^1))T_1^c,$$
and 
$$L(s_2+\frac{j-k}{2},\wt \pi\times \wt\tau_2)=\sum_{d\ge 0}\Tr(\Sym^d(\wt\alpha\otimes \beta^2))T_2^d,$$
where $\wt \alpha=\diag(a_1^{-1},\dots,a_l^{-1})$ is the Satake parameter for $\wt \pi,$ 
we get that
\begin{align}\label{eq: unramified6}
L(s_1+\frac{k-j}{2},\pi\times \tau_1)L(s_2+\frac{j-k}{2},\wt \pi\times \wt\tau_2)=\sum_{c\ge 0, d\ge 0}\Tr(\Sym^c(\alpha\otimes \beta^1))\Tr(\Sym^d(\wt\alpha\otimes \beta^2))T_1^c T_2^d.
\end{align}
Comparing \eqref{eq: unramified5} and \eqref{eq: unramified6}, in order to prove Proposition~\ref{proposition-unramified-computation}, it suffices to show 
\begin{align}\label{eq: unramified-reduction}
\begin{split}
\Tr(\Sym^c(\alpha\otimes \beta^1))\Tr(\Sym^d(\wt\alpha\otimes \beta^2))=\sum_{e\ge 0}\sum_{\substack{\bx\in T^+(m),\by\in T^+(n), e\ge 0\\ |\bx|=c-e, |\by|=d-e}}&S_{(\bx,0,\by^*)}(\alpha)S_{\bx}(\beta^1)S_{\by}(\beta^2)\\
&\cdot  \Tr(\Sym^e(\beta^1\otimes \beta^2)).
\end{split}
\end{align}
By \cite{Jacquet-Shalika-EulerI}*{Proposition 2.4},  we have 
$$\Tr(\Sym^e(\beta^1\otimes \beta^2))=\sum_{\bz\in T^+(n), |\bz|=e}S_{(\bz,0_{m-n})}(\beta^1)S_{\bz}(\beta^2).$$
Here $\bz=(z_1,\dots,z_n)$ can be identified with a partition of $e=|\bz|$ with at most $n$-parts (since $m\ge n$ by our assumption) and $S_{\bz}$ (resp. $S_{(\bz,0_{m-n})}$) is the Schur polynomial defined by $\bz$ with $n$ (resp. $m$) variables.
Similarly,
\begin{align*}\Tr(\Sym^c(\alpha\otimes \beta^1))&=\sum_{\mathbf{u}\in T^+(m), |\bu|=c}S_{(\bu,0_{l-m})}(\alpha)S_{\bu}(\beta^1),\\
\Tr(\Sym^d(\wt\alpha\otimes \beta^2))&=\sum_{\mathbf{v}\in T^+(n), |\bv|=d}S_{(\bv,0_{l-n})}(\wt \alpha)S_{\bu}(\beta^2).
\end{align*}
A simple matrix calculation shows that 
$$S_{(\bv,0_{l-n})}(\wt \alpha)=S_{(0_{l-n},\bv^*)}(\alpha).$$
See also \cite{Fulton-Harris}*{Exercise 15.50} for a representation theoretic explanation of this formula. Thus the left hand side of \eqref{eq: unramified-reduction} becomes
\begin{align*}
LHS=\sum_{\bu\in T^+(m), |\bu|=c}\sum_{\bv\in T^+(n), |\bv|=d}S_{(\bu,0_{l-m})}(\alpha) S_{(0_{l-n},\bv^*)}(\alpha)S_{\bu}(\beta^1)S_{\bv}(\beta^2),
\end{align*}
while the right side of \eqref{eq: unramified-reduction} becomes
\begin{align*}
RHS=\sum_{\bz\in T^+(n)} \sum_{\substack{\bx\in T^+(m),\by\in T^+(n), e\ge 0\\ |\bx|=c-|\bz|, |\by|=d-|\bz|}}&S_{(\bx,0,\by^*)}(\alpha)S_{\bx}(\beta^1)S_{(\bz,0)}(\beta^1)S_{\by}(\beta^2)S_{\bz}(\beta^2)
\end{align*}
By Littlewood-Richardson rule, see \cite{Fulton-Harris}*{(A.8)} or \cite{Mac}*{\S I.9}, we have 
\begin{align*}
S_{\bx}(\beta^1)S_{\bz}(\beta^1)=&\sum_{\bu\in T^+(m), |\bu|=c}c_{\bx, \bz}^\bu S_{\bu}(\beta^1),\\
S_{\by}(\beta^2)S_{\bz}(\beta^2)=&\sum_{\bv\in T^+(n), |\bv|=d} c_{\by, \bz}^\bv S_{\bv}(\beta^2),
\end{align*}
where in the first equation, $(\bz,0_{m-n})$ is identified with $|\bz|$ as a partition of $e=|\bz|$ with at most $n$ parts, and $c_{\bx, \bz}^\bu, c_{\by, \bz}^\bu$ are the Littlewood-Richardson coefficients as defined in \cite{Fulton-Harris}*{page 454} or \cite{Mac}*{\S I.9.2}. Thus 
\begin{align*}
RHS=\sum_{\bu\in T^+(m), |\bu|=c}\sum_{\bv\in T^+(n),|\bv|=d}\sum_{\substack{\bx\in T^+(m), \by,\bz\in T^+(n)\\ |\bx|+|\bz|=c, |\by|+|\bz|=d}}c_{\bx, \bz}^ \bu c_{\by, \bz}^ \bv S_{(\bx,0,\by^*)}(\alpha) S_{\bu}(\beta^1)S_{\bv}(\beta^2).
\end{align*}
Thus in order to prove \eqref{eq: unramified-reduction} and hence Proposition~\ref{proposition-unramified-computation}, it suffices to prove that for any $\bu\in T^+(m), \bv\in T^+(n)$, one has
\begin{align}\label{eq: unramified-reduction2}
S_{(\bu,0_{l-m})}(\alpha) S_{(0_{l-n},\bv^*)}(\alpha)=\sum_{\substack{\bx\in T^+(m), \by,\bz\in T^+(n)\\ |\bx|+|\bz|=c, |\by|+|\bz|=d}}c_{\bx, \bz}^ \bu c_{\by,\bz }^\bv S_{(\bx,0,\by^*)}(\alpha).
\end{align}
For $\bv=(v_1,\dots,v_n)\in T^+(n)$, we write $\wt \bv=(v_1,\dots,v_1,v_1-v_n,\dots,v_1-v_2,0)\in T^+(l)$. Then $S_{(0_{l-n},\bv^*)}(\alpha)=S_{\bv}(\alpha) D_{-v_1}(\alpha)$, where $D_{-v_1}(\alpha)=\det^{-v_1}(\alpha)$ following the notation of \cite{Fulton-Harris}*{\S 15.5}. Thus using Littlewood-Richardson rule again, we have
\begin{align*}
S_{(\bu,0_{l-m})}(\alpha) S_{(0_{l-n},\bv^*)}(\alpha)&=D_{-v_1}(\alpha)\sum_{\lambda\in T^+(l), |\lambda|=|\wt v|+|\bu|}c_{\wt \bv, \bu}^{ \lambda}S_\lambda(\alpha).
\end{align*}
Write $\lambda=(\lambda_1,\dots,\lambda_l)$. By the definition of Littlewood-Richardson coefficients, if $c_{\wt \bv, \bu}^\lambda\ne 0$, we must have $\lambda_{m+1}=\dots=\lambda_{l-n-1}=v_1$, which means that $S_{\lambda}\cdot D_{-v_1}=S_{(\lambda_1-v_1,\dots,\lambda_l-v_1)}$ must be of the form $S_{(\bx,0_{l-m-n},\by^*)}$ for $\bx\in T^+(m)$ and $\by\in T^+(n)$. Thus we get 
\begin{align*}
S_{(\bu,0_{l-m})}(\alpha) S_{(0_{l-n},\bv^*)}(\alpha)=\sum_{\bx\in T^+(m), \by\in T^+(n) } c_{\wt \bv, \bu}^{ \lambda}S_{(\bx,0,\by^*)},
\end{align*}
where $\lambda=(\lambda_1,\dots,\lambda_l)=(\bx,0,\by^*)+(v_1,\dots,v_1).$ Note that $|\bu|-|\bv|=|\bx|-|\by|$. Thus in order to prove \eqref{eq: unramified-reduction2}, it suffices to show that for any fixed $\bu,\bx \in T^+(m)$ and $\bv, \by\in T^+(n)$ with $|\bu|-|\bx|=|\bv|-|\by|$,
\begin{align}\label{eq: Tao's formula}
c_{\wt \bv,\bu}^{\lambda}=\sum_{\bz\in T^+(n)}c_{\bx, \bz}^{ \bu} c_{\by, \bz}^{ \bv},
\end{align}
where $ \lambda=(\lambda_1,\dots,\lambda_l)=(\bx,0,\by^*)+(v_1,\dots,v_1).$ The formula \eqref{eq: Tao's formula} was proved by Professor T. Tao in a MathOverflow answer \cite{Tao} using the hive model for Littlewood-Richardson coefficients introduced in \cite{Knutson-Tao}. A proof of \eqref{eq: Tao's formula} based on Tao's MathOverflow answer \cite{Tao} will be reproduced in \S\ref{subsection: Proof of Tao's formula} after we introduce some necessary notations and tools. 
 \end{proof}
\begin{remark}
Here we give an example of \eqref{eq: unramified-reduction2}. We take $l=4,m=2,n=1$ and $\bu=(2,1),\bv=(2)$. One can check that there are 3 choices of $\bz$, which are $\bz=(0),\bz=(1), \bz=(2)$, and correspondingly, there are 3 choices of $\by$ given by $\by=(2),\by=(1),\by=(0)$. When $\bz=(0)$, we must have $\bx=(2,1)$ and when $\bz=(2)$, we must have $\bx=(1)=(1,0)$. But when $\bz=(1)$, there are two choices of $\bx$, which are $\bx=(1,1)$ or $\bx=(2)=(2,0)$. One can check that in each case, $c_{\bx, \bz}^{ \bu}c_{\by, \bz}^{ \bv}=1$. Thus \eqref{eq: unramified-reduction2} becomes
$$S_{(2,1,0,0)}\cdot S_{(0,0,0,-2)}=S_{(2,1,0,-2)}+S_{(1,0,0,0)}+S_{(1,1,0,-1)}+S_{(2,0,0,-1)},$$
which could be checked directly using Littlewood-Richardson rule by noting that $S_{(0,0,0,-2)}=S_{(2,2,2,0)}\cdot D_{-2}$, where $D_{-2}=\det^{-2}$.
\end{remark}

\subsection{Proof of Tao's formula \eqref{eq: Tao's formula}}\label{subsection: Proof of Tao's formula} An integral \textit{n-hive} is an array of integers $a_{ij}$ for $0\le i, j, i+j\le n$ placed in the vertices of triangles of the following shape 
\begin{figure}[H]
 \includegraphics[width=0.35\textwidth ]{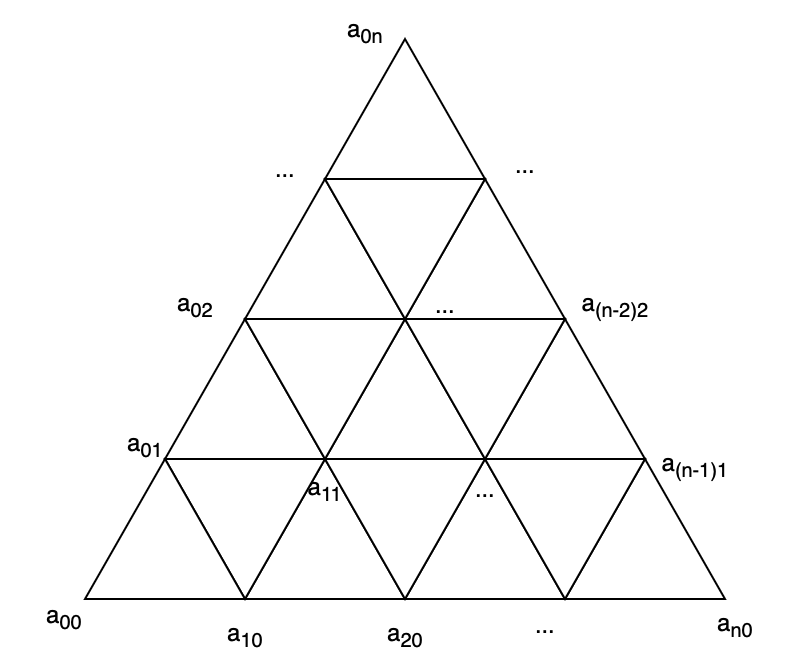} 
\caption{hive}
\label{fig:hive1}
\end{figure}
\noindent which satisfies all of the following rhombus inequalities: for each rhombus of the following types
$$ \includegraphics[width=0.35\textwidth ]{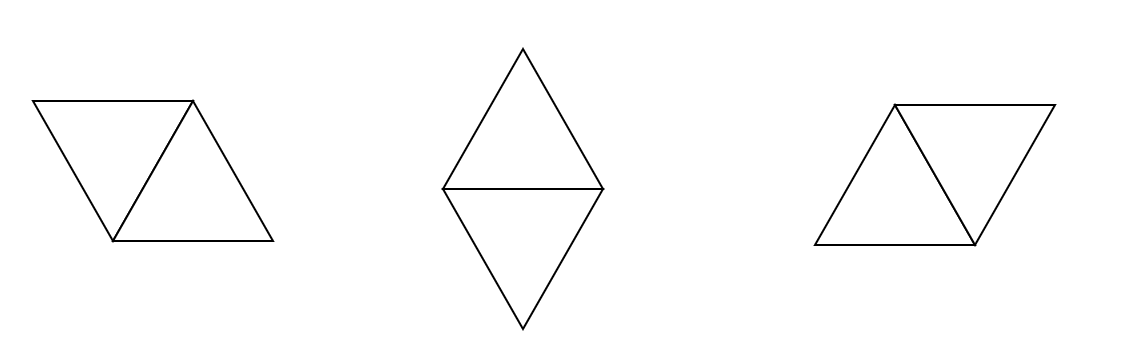}$$
the sum of the two integers at the obtuse vertices must be greater than or equal to the sum of the two integers at the acute vertices. 
\begin{theorem}[Knutson-Tao, \cite{Knutson-Tao}]\label{theorem: KN} Let $\bx=(x_1,\dots,x_n), \by=(y_1,\dots,y_n), \bz=(z_1,\dots,z_n)$ be partitions with $|\bz|=|\bx|+|\by|$, then $c_{\bx,\by}^{\bz}$ is the number of $n$-hives with boundary labels
$$ \includegraphics[width=0.35\textwidth ]{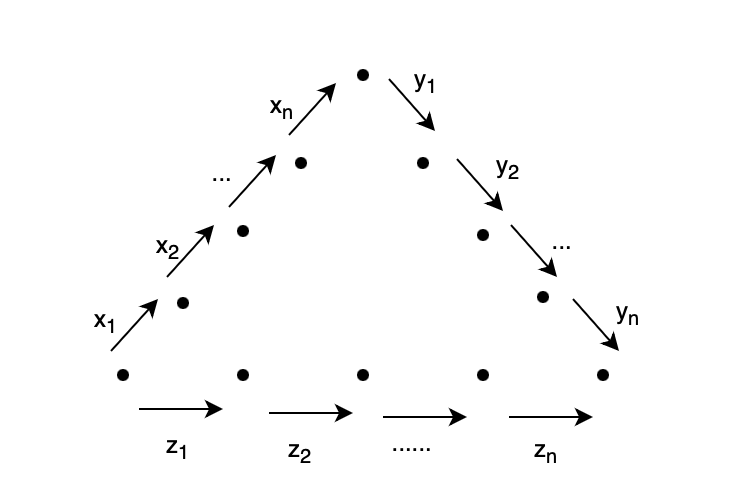}$$
Here the arrow and the number $x_i$ (resp. $y_j, z_k$) on the arrow indicates that the numbers increase by $x_i$ (resp. $y_j,z_k$) along the direction indicated by the arrow. One can normalize the above $n$-hive by assign any integer to any fixed vertex. 
\end{theorem}
We note that different normalization will give the same number of hives. The above theorem is proved in \cite{Knutson-Tao}. See also the appendix of \cite{Buch} for a different proof given by W. Fulton.
\begin{remark}
We give a simple example which also appeared in \cite{Buch}. We have $c_{(2,1),(2,1)}^{(3,2,1)}=2$, which can be computed in the following way. There are exactly two $3$-hives with boundary conditions given below,
$$\begin{matrix} &&&3&&&\\ &&&&& \\ &&3&&5&&\\ &&&&& \\ &2&&x&&6&\\ &&&&&& \\0&&3&&5&&6 \end{matrix}, $$
which are given by $x=4,5$.
\end{remark}

We temporarily call the following object an \textit{anti-n-hive}: an array of integers placed in the vertices of triangles of the shape as Figure \ref{fig:hive1} which satisfies the ``reverse" rhombus inequalities: for each rhombus below
$$ \includegraphics[width=0.35\textwidth ]{image/hive2.png}$$
the sum of the two integers at the obtuse vertices must be less than or equal to the sum of the two integers at the acute vertices.

For any $n$-hive, if we switch the sign of the number at each vertices, we will get an anti-$n$-hive. Note that, this process will change the boundary conditions, which gives us the following direct corollary.
\begin{corollary}\label{corollary: KN}
Let $\bx=(x_1,\dots,x_n), \by=(y_1,\dots,y_n), \bz=(z_1,\dots,z_n)$ be partitions with $|\bz|=|\bx|+|\by|$, then $c_{\bx,\by}^{\bz}$ is the number of anti-$n$-hives with boundary labels
$$ \includegraphics[width=0.35\textwidth ]{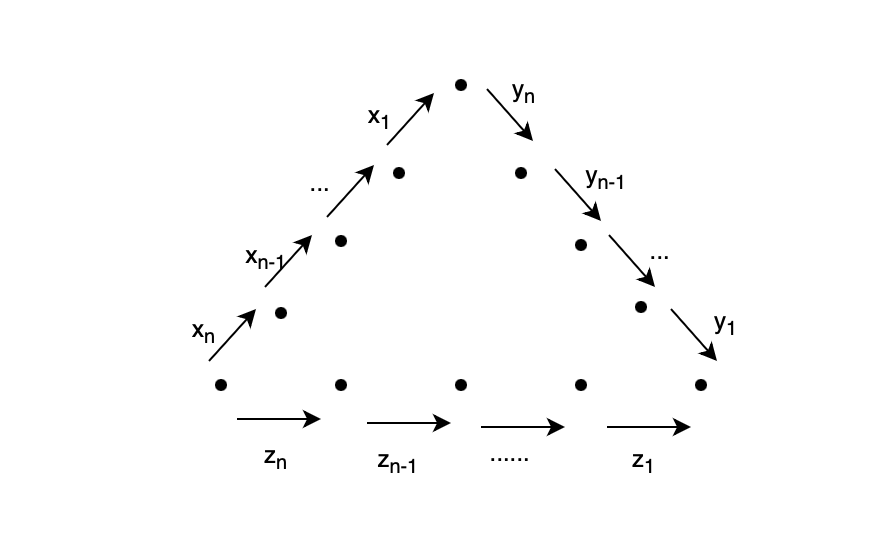}$$
Here the arrow and the number $x_i$ (resp. $y_j, z_k$) on the arrow indicates that the numbers increase by $x_i$ (resp. $y_j,z_k$) along the direction indicated by the arrow. One can normalize the above $n$-hive by assign any integer to any fixed vertex. 
\end{corollary}
Now we can prove Tao's formula \eqref{eq: Tao's formula}, which we restate it below.
\begin{proposition}
Let $l,m,n$ be non-negative integers with $l\ge m+n+1$ and $m\ge n$. Given $\bx, \bu\in T^+(m), \by, \bv\in T^+(n)$ with $|\bu|-|\bx|=|\bv|-|\by|\ge 0 $, then 
$$c_{\wt \bv, \bu}^{\lambda}=\sum_{\bz\in T^+(n)}c_{\bx,\bz}^{\bu} c_{\by,\bz}^{\bv}.$$
Here $\bu=(u_1,\dots,u_m), \bx=(x_1,\dots,x_m), \by=(y_1,\dots,y_n), \bv=(v_1,\dots,v_n)$, $\by^*=(-y_n,\dots,-y_2,-y_1)$, $\wt \bv=(0_{l-n},\bv^*)+(v_1,\dots,v_1)=(v_1,\dots,v_1,v_1-v_n,\dots,v_1-v_2,0)\in T^+(l)$, and $\lambda=(\bx,0_{l-m-n},\by^*)+(v_1,\dots,v_1)\in T^+(l)$. Moreover, $\bu$ in $c_{\wt \bv, \bu}^\lambda$ is viewed as an element in $T^+(l)$ in the obvious way, namely, $\bu=(\bu,0_{l-n})$.
\end{proposition}
\begin{proof}
By Theorem \ref{theorem: KN} and Corollary \ref{corollary: KN}, one can see that $c_{\wt \bv, \bu}^{\lambda}$ is the number of anti-$l$-hives with boundary conditions indicated below,
$$ \includegraphics[width=0.4\textwidth ]{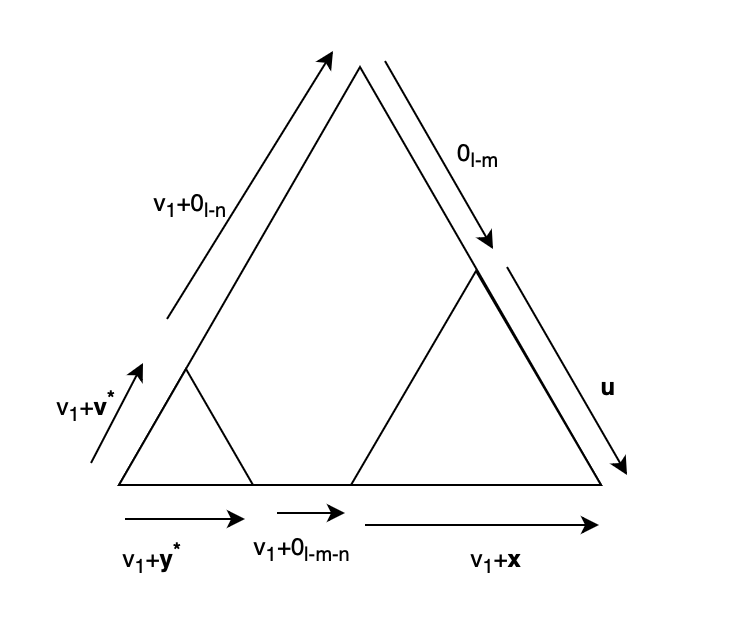}$$
where $v_1$ in the left side boundary and bottom boundary means $(v_1,\dots,v_1)\in T^+(l)$. Here the two interior line segments are not important here. For each hive above, we assume that its vertex integers are given by $(a_{ij})_{0\le i, j, i+j\le l}$ placed as in Figure \ref{fig:hive1}. Then $(a_{ij}-(i+j)v_1)_{0\le i, j, i+j\le l}$ is also an anti-$l$-hive which has the boundary conditions as indicated in the following Figure \ref{fig:hive6}. We also normalized the anti-$l$-hive so that the top vertex has value $0$.
\begin{figure}[H]
 \includegraphics[width=0.4\textwidth ]{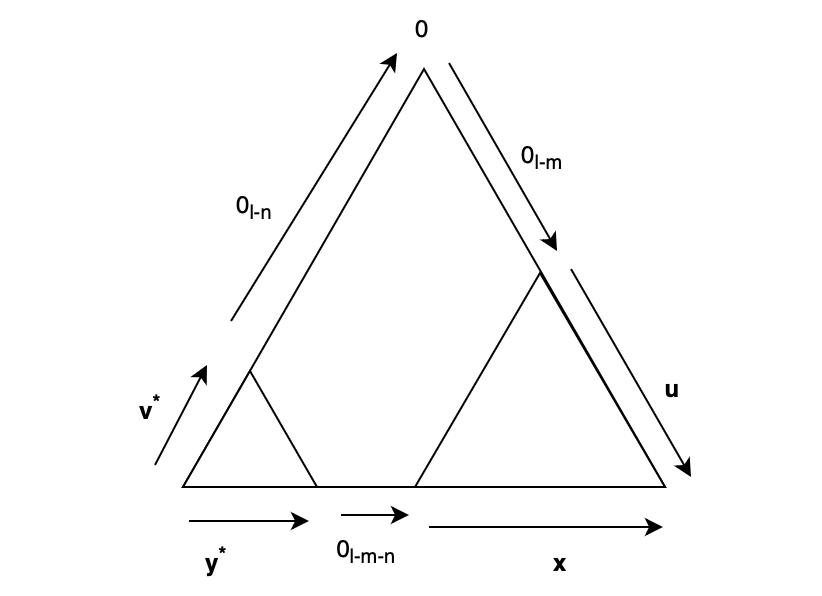} 
\caption{boundary condition for anti-hives which represents $c_{\wt \bv, \bu}^\lambda$}
\label{fig:hive6}
\end{figure}
\noindent Thus $ c_{\wt \bv, \bu}^\lambda$ is the number of anti-$l$-hives with boundary conditions as in Figure \ref{fig:hive6}. Using the reverse rhombus inequality, we can check that an anti-$l$-hive as above must vanish completely in the quadrilateral $ABEF$ (including each sides)  in Figure \ref{fig:hive7}. Moreover, inside the trapezoid $BCDE$, the values of the hive on each horizontal line are the same. In particular, this means that there exists a $\bz\in T^+(n)$ such that the boundary condition on $CB$ and $DE$ are both given by $\bz^*$.

\begin{figure}[H]
 \includegraphics[width=0.42\textwidth ]{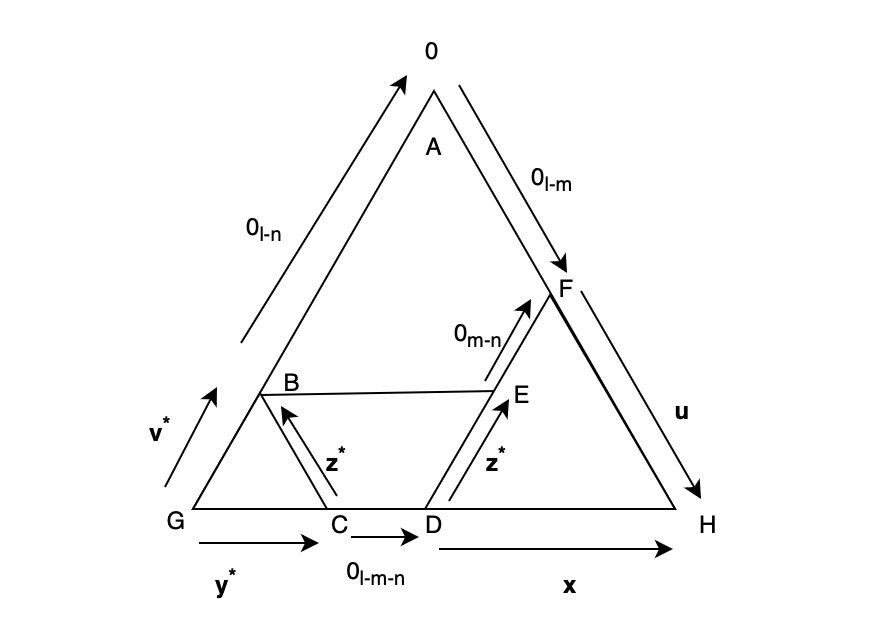} 
\caption{}
\label{fig:hive7}
\end{figure}

\noindent Thus such a hive is uniquely determined by its values in the anti-hives $BGC$ and $FDH$, with the indicated boundary conditions as in Figure \ref{fig:hive7}. Conversely, given anti-hives $BGC$ and $FDH$ with boundary conditions as in Figure \ref{fig:hive7}, we get an anti-hive with the boundary condition as in Figure \ref{fig:hive6} using a reverse process. Finally, note that the number of anti-hives $BGC$ is $c_{\by,\bz}^{\bv}$ and the number of anti-hives $FDH$ is $c_{\bx,\bz}^{\bu}$. Thus we get
$$c_{\wt \bv,\bu}^{\lambda}=\sum_{\bz\in T^+(n)} c_{\bx,\bz}^{\bu} c_{\by,\bz}^{\bv}.$$
This concludes the proof.
\end{proof}

\section{A local converse theorem}\label{section-preparation}
In the rest of this paper, we assume that $F$ is a non-archimedean local field of characteristic different from 2. Let $\CO$ be the ring of integers of $F$, $\fp$ be the maximal ideal of $\CO$ and let $\varpi\in \fp$ be a fixed uniformizer. The purpose of the rest of this paper is to prove the following
\begin{theorem}\label{theorem: main}
Let $l$ be a positive integer and let $\pi_1,\pi_2$ be two irreducible supercuspidal representations of $\GL_l(F)$ with the same central character. If $\Gamma(\Bs,\pi_1\times (\tau_1,\tau_2),\psi;0)=\Gamma(\Bs,\pi_2\times (\tau_1,\tau_2),\psi;0)$ for all irreducible generic representations $\tau_1$ (resp. $\tau_2$) of $\GL_m(F)$ (resp. $\GL_n(F)$) with $0\le n\le [l/2], 0\le m\le [l/2]$, then $\pi_1\cong \pi_2$.
\end{theorem}
\begin{remark}\label{remark: even case}
If $l=2r$ is even and $m=n=r$, we have not defined the gamma factor $\Gamma(\Bs,\pi\times (\tau_1,\tau_2),\psi ;0)$ yet, because our local zeta integral \eqref{eq: local zeta integral} and hence our local gamma factor defined from that in Proposition \ref{proposition: existence of gamma} require $m+n<l$. In the case if $l=2r, m=n=l$, the corresponding local gamma factor used in Theorem \ref{theorem: main} is the one defined from the local zeta integral of unitary group $\RU_{E/F}(2r)\times \Res_{E/F}(\GL_r)$ at a split place, see \cite{BAS: Uodd} and \cite{Morimoto}. Actually, the properties of this gamma factor is well studied. In particular, it has been shown that it is the product of Jacquet--Piatetski-Shapiro--Shalika local gamma factors after normalization, see \cite{Morimoto}. We will review its definition in \S\ref{subsection: review of gamma}.
\end{remark}

\begin{remark}\label{remark: JPSS gamma factor again}
Note that if $m=n=0$, then condition $\Gamma(\Bs,\pi_1\times (\tau_1,\tau_2),\psi;0)=\Gamma(\Bs,\pi_2\times (\tau_1,\tau_2),\psi;0) $  is empty. If $m>0$ and $n=0$, the corresponding gamma factor $\Gamma(\Bs,\pi_1\times (\tau_1,\tau_2),\psi;0) $ is exactly a Jacquet--Piatetski-Shapiro--Shalika local gamma factor up to a shift. 
\end{remark}

Here we recall Jacquet's local converse conjecture
\begin{conj}\label{conjecture: Jacquet}
Let $\pi_1,\pi_2$ be two irreducible generic representations of $\GL_l(F)$ with the same central character. If $\gamma(s,\pi_1\times \tau,\psi)=\gamma(s,\pi_2\times \tau,\psi)$ for any irreducible generic representation $\tau$ of $\GL_m(F)$ with $1\le m\le [l/2]$, then $\pi_1\cong \pi_2$.
\end{conj}
One can assume that $\pi_1,\pi_2$ are supercuspidal and remove the central character restriction after the work of \cite{JNS}. The above conjecture was proved in \cite{chai} and \cite{Jacquet-Liu} independently. In the next remark, we will explain that our Theorem \ref{theorem: main} indeed gives a new proof of Conjecture \ref{conjecture: Jacquet}.

\begin{remark}\label{remark: condition C(k)}
 We write $\CC(0)$ for the condition that the two representations $\pi_1,\pi_2$ have the same central character, which is always assumed. For $t\ge 1$, we denote $\CC(t):=\CC(t;\pi_1,\pi_2)$ the following condition for $\pi_1,\pi_2:$
\textit{
$$\Gamma(\Bs,\pi_1\times (\tau_1,\tau_2),\psi;0)=\Gamma(\Bs,\pi_2\times (\tau_1,\tau_2),\psi;0),$$
for any irreducible generic representation $\tau_1$ (resp. $\tau_2$) of $\GL_m(F)$ (resp. $\GL_n(F)$) with $0\le m,n\le t$.} To compare our result with Jacquet's local converse conjecture, we also denote by $\CC'(t):=\CC'(t; \pi_1,\pi_2)$ the condition: $\gamma(s,\pi_1\times\tau,\psi)=\gamma(s,\pi_2\times \tau,\psi)$ for any irreducible generic representation $\tau$ of $\GL_m(F)$ with $1\le m\le t$.    Note that $\gamma(1-s,\wt\pi\times \wt\tau,\psi)\gamma(s,\pi\times \tau,\psi^{-1})=1$, see \cite{JPSS}*{\S2.9 and \S3.1}. Thus, condition $\CC(t)$ is, in fact, equivalent to $\CC'(t)$ by Proposition \ref{proposition-gamma-comparison}. So our proof of Theorem \ref{theorem: main} gives a new proof of Jacquet's local converse conjecture.
\end{remark}

The proof of Theorem \ref{theorem: main} will be given in the next section. In the rest of this section, we introduce some necessary tools that will be used in the proof of Theorem \ref{theorem: main}.
\subsection{On the gamma factors for $\GL_{2r}\times (\GL_r,\GL_r)$}\label{subsection: review of gamma}
Recall that if $m+n\le l-1$, for generic representation $\pi$ of $\GL_l(F)$, $\tau_1$ (resp. $\tau_2$) of $\GL_m(F)$ (resp. $\GL_n(F)$), our local gamma factor $\Gamma(\Bs,\pi\times (\tau_1,\tau_2),\psi)$ is defined by the  local functional equation
$$\Psi(W,M_{w_{m,n}}(\xi_{\Bs});0)=\Gamma(\Bs,\pi\times (\tau_1,\tau_2),\psi;0)\Psi(W,\xi_{\Bs};0),$$
for all $W\in \CW(\pi,\psi)$ and $\xi_{\Bs}\in \CW(\Bs,(\tau_1,\tau_2),\psi^{-1}).$ See Proposition \ref{proposition: existence of gamma}. From now on, we will drop the index $j=0$ in the local gamma factor and simply write $\Gamma(\Bs,\pi\times (\tau_1,\tau_2),\psi):=\Gamma(\Bs,\pi\times (\tau_1,\tau_2),\psi;0)$. For $W\in \CW(\pi,\psi)$, $\rho(\gamma_{m,n}^{-1})W$ is also an element in $\CW(\pi,\psi)$. Thus, we have
\begin{align}\label{eq: modified LFE}
\Psi(\rho(\gamma_{m,n}^{-1})W,M_{w_{m,n}}(\xi_{\Bs});0)=\Gamma(\Bs,\pi\times (\tau_1,\tau_2),\psi)\Psi(\rho(\gamma_{m,n}^{-1})W,\xi_{\Bs};0),
\end{align}
for all $\xi_\Bs\in \CW(\Bs,(\tau_1,\tau_2),\psi^{-1})$. Here $\rho$ denotes the right translation, and $\gamma_{m,n}$ is the element in $\GL_l$ as defined after \eqref{eq: local zeta integral}.  The local functional equation \eqref{eq: modified LFE} is the one we will use to prove our local converse theorem.

As explained in Remark \ref{remark: even case}, we also need the local gamma factors for $\Gamma(\Bs, \pi\times (\tau_1,\tau_2),\psi)$ when $l=2r$ and $m=n=r$, which are not covered in our previous sections. This local gamma factor has been defined in \cite{BAS: Uodd} and studied in \cite{Morimoto}. We recall the definition now.

We first endowed $F^{2r}\oplus F^{2r}$ a symplectic structure $\pair{~,~}$ defined by 
$$\pair{(u_1,u_2),(v_1,v_2)}=2(u_1J_{2r}v_2^t-v_2J_{2r}u_2^t),$$
where $u_i, v_i\in F^{2r}$ are viewed as row vectors.
 For a nontrivial additive character $\psi$ of $F$ and for a character $\mu$ of $F^\times$, we can consider the Weil representation $\omega_{\psi^{-1},\mu, \mu^{-1}}$ of $\GL_{2r}(F)$, see \cite{Morimoto}*{\S 2.2}. Note that we used a slightly different normalization. The Weil representation $\omega_{\psi^{-1},\mu, \mu^{-1}}$ can be realized on the space $\CS(F^{r}\times F^{r})$, the Bruhat-Schwatz functions on $F^{2r}$. This is the Schr\"odinger model of the Weil representation. For example, we have the well-known formula
$$\left(\omega_{\psi^{-1}.\mu, \mu^{-1}}\left(\bpm I_r &X\\ &I_r \epm \right) \phi \right)(x,y)=\psi(x XJ_r y^t)\phi(x,y), X\in \Mat_{r\times r}(F).$$
In the following, we assume that $\mu$ is understood and omit it from the notation.

Now let $\pi$ be an irreducible generic representation of $\GL_{2r}(F)$, $(\tau_1,\tau_2)$ be a pair of irreducible generic representations of $\GL_r(F)$ and $\Bs=(s_1,s_2)$ be a pair of complex numbers. For $W\in \CW(\pi,\psi), \xi_\Bs\in \CW(\Bs,(\tau_1,\tau_2,),\psi^{-1})$, and $\phi\in \CS(F^{2r})$, we consider the local zeta integral
$$\Psi(W,\xi_{\Bs},\phi)=\int_{N_{2r}(F)\bs \GL_{2r}(F)}W(g)\xi_{\Bs}(g) (\omega_{\psi^{-1}}(g)\phi)(e_r,e_r)dg,$$
where $e_r\in F^r$ is the vector $(0,0,\dots,0,1)$. There exists a meromorphic function $\Gamma(\Bs, \pi\times(\tau_1,\tau_2),\mu,\psi)$ such that 
\begin{equation}\label{eq: LFE even case}\Psi(W, M_{w_{r,r}}\xi_{\Bs},\phi)=\Gamma(\Bs, \pi\times(\tau_1,\tau_2),\mu,\psi)\Psi(W,\xi_{\Bs},\phi)\end{equation}
for any $W\in \CW(\pi,\psi), \xi_{\Bs}\in \CW(\Bs, (\tau_1,\tau_2),\psi^{-1})$ and $\phi\in \CS(F^r\times F^r).$ Note that in \cite{BAS: Uodd} and \cite{Morimoto}, there is only a single complex variable involved in the local zeta integral and local gamma factor. Here we still use the two variable case.
\begin{remark}
In \cite{Morimoto}, the field $F$ is assumed to have characteristic zero. One can extend their results to odd characteristic local field $F$ without significant difficulty. See the work of Jo \cite{Jo} on similar situations. If we don't assume the existence of local gamma factors for a local field of positive characteristic, then we also have to require that our field $F$ be of characteristic zero in Theorem \ref{theorem: main} when $l$ is even.
\end{remark}

\subsection{Howe vectors}\label{subsection: Howe vector} Our strategy for the proof of Theorem \ref{theorem: main} is along the lines of that given in \cite{Sp(2r)} and \cite{U(2r+1)}. One basic tool for us is the partial Bessel functions associated with Howe vectors, as developed in \cite{Baruch}. Here we recall the basic construction. Let $\psi$ be a fixed unramified additive character of $F$ and we also view $\psi$ as a character of the maximal unipotent subgroup $N_l\subset \GL_l(F)$ in the usual way. For an integer $i\ge 0$, we consider the open compact subgroup $K_{\GL_l}^i:=I_l+\Mat_{l\times l}(\fp^i)$ of $\GL_l(F)$. Consider the character $\theta_i$ of $K_{\GL_l}^i$ defined by 
$$\theta_i(k)=\psi(\varpi^{-2i}(\sum_{s=1}^{l-1}k_{s,s+1})), \quad k=(k_{st})_{1\le s,t\le l}\in K_{\GL_l}^i.$$ 
One can check that $\theta_i$ is indeed a character of $K_{\GL_l}^i$. Consider the element 
$$d_i=\diag(\varpi^{-i(l-1)},\varpi^{-i(l-3)},\dots,\varpi^{i(l-3)},\varpi^{i(l-1)}),$$
and $H_l^i=d_iK^id_i^{-1}$, which is still an open compact subgroup of $\GL_l(F)$. One sees that $H_l^i$ has the form
$$H_l^i=\bpm 1+\fp^i& \fp^{-i} &\fp^{-3i} &\dots\\ \fp^{3i} &1+\fp^i &\fp^{-i} &\dots\\ \fp^{5i} &\fp^{3i} & 1+\fp^i&\dots\\ \dots&\dots&\dots&\dots \epm.$$
 We consider the character $\psi_i$ of $H_l^i$ defined by 
$$\psi_i(h):=\theta_i(d_i^{-1}hd_i), \quad h\in H_l^i.$$
For a subgroup $U\subset \GL_l(F)$, we denote $U^i:=U\cap H_l^i$. For example, $N_l^i$ denotes $N_l\cap H_l^i$. We also usually drop $l$ from the notation if $l$ is understood. It is easy to see that $\psi_i|_{N_l^i}=\psi|_{N_l^i}$.

 Let $(\pi,V)$ be an irreducible generic representation of $\GL_l(F)$ and for $v\in V$, we consider 
$$v_i=\frac{1}{\vol(N_l^i)}\int_{N_l^i}\psi_i^{-1}(u)\pi(u)vdu.$$
If $W\in \CW(\pi,\psi)$ is the Whittaker function associated with $v$, then we denote $W_i=W_{v_i}$. Note that 
$$W_i(u_1gu_2)=\psi(u_1)\psi_i(u_2)W_i(g), \quad \forall g\in \GL_l(F), u_1\in N_l, u_2\in N_l^i.$$
Actually, that exists a positive integer $C(v)>0$, such that  $W_i$ satisfies the additional quasi-invariance property \begin{equation}\label{eq: quasi-invariance of Howe vector} W_i(ugh)=\psi(u)\psi_i(h)W_i(g)\end{equation} for all $u\in N_l,g\in \GL_l(F),h\in H_l^i$ if $i>C(v)$, see \cite{Baruch}*{Lemma 3.2}. According to the proof \cite{Baruch}*{Lemma 3.2}, one can take $C(v)$ to be any integer such that $v$ is fixed by $\pi(K^{C(v)}_l)$.

Let $\omega$ be a character of $F^\times$ and we consider the space $C_c^\infty(\GL_l(F),\omega)$ consisting of smooth functions $f$ on $G$ such that $f$ is compactly supported modulo $Z_l$, the center of $\GL_l(F)$, and $f(zg)=\omega(z)f(g)$.  If $\pi$ is supercuspidal, let $\CM(\pi)$ be the space of matrix coefficients of $\pi$. Then $\CM(\pi)\subset C_c^\infty(\GL_l(F),\omega_\pi)$. For $f\in \CM(\pi)$, following \cite{CST}*{page 2089}, we consider the function
$$W^f(g)=\int_{N_l}\psi^{-1}(u)f(ug)du.$$
Note that the integral is convergent by assumption and defines an element in $\CW(\pi,\psi)$. Moreover, for an appropriate choice of $f$, we can assume that $W^f(I_l)=1$. See \cite{CST}*{page 2089-2090}. Thus, we can consider $W^f_i$. We also use the notation
$$\CB_i(g,f)=W^f_i(g), \quad g\in \GL_l(F).$$

\subsection{Weyl elements which support Bessel functions}\label{subsection: Weyl elements which support Bessel functions}
Let $\Delta=\Delta(\GL_l)$ be the set of simple roots of $\GL_l(F)$. Then $\Delta=\{\alpha_k: 1\le k\le l-1\}$, where 
$$\alpha_k(\diag(t_1,\dots,t_l))=t_{k}/t_{k+1},  \quad \diag(t_1,\dots,t_l)\in T_l(F).$$

 Let $\bW=\bW(\GL_l)$ be the Weyl group of $\GL_l(F)$. We sometimes identify $\bW$ with the permutation matrix in $\GL_l(F)$. Denote by $e$ the identiy element in $\bW$, which is represented by $I_l\in \GL_l(F)$. For $w\in \bW$, denote $C(w)=BwB$, where $B=B_l$ is the upper triangular subgroup of $\GL_{l}(F)$. There is a Bruhat order on $\bW$, which is recalled as follows. Given $w_1,w_2\in \bW$, then $w_1\le w_2$ (or $w_2\ge w_1$) if and only if $C(w_1)\subset \ov{C(w_2)}$. For $w\in \bW$, we denote $\Omega_w=\coprod_{w'\ge w}C(w')$. Then $C(w)$ is closed in $\Omega_w$ and $\Omega_w$ is open in $G$.

Let $\RB(\GL_l)=\wpair{w\in \bW(\GL_l): \alpha\in \Delta, w\alpha>0\implies w\alpha\in \Delta}$, which is the set of Weyl elements that can support partial Bessel functions.

 Let $w_0=J_{l}\in \GL_l(F)$, which represents the longest Weyl element of $\GL_l(F)$.  It is well-known that $w\in \RB(\GL_l)$ if and only if $w_0 w$ is the longest Weyl element of the Levi subgroup of a standard parabolic subgroup of $\GL_l(F)$. For $w\in \RB(\GL_l)$, let $P_w=M_wN_w$ be the corresponding parabolic subgroup such that $w_0 w=w_0^{M_w}$, where $M_w$ is the Levi subgroup of $P_w$ and $w_0^{M_w}$ is the longest Weyl element of $M_w$. Let $\theta_w$ be the subset of $\Delta$ which consists of all simple roots in $M_w$. Then we have the relation 
$$\theta_w=\wpair{\alpha\in \Delta| w\alpha>0}\subset \Delta.$$
The assignment $w\mapsto \theta_w$ is a bijection between $\RB(\GL_l)$ and subsets of $\Delta$. Moreover, it is known that the assignment $w\mapsto \theta_w$ is order-reversing, i.e., $w'\le w$ if and only if $\theta_w\subset \theta_{w'}$, see \cite{CPSS05}*{Proposition 2.11}. For example, we have $\theta_{w_0}=\emptyset$ and $\theta_{e}=\Delta$.

 Given a subset $\theta\subset \Delta$, we will write the corresponding Weyl element in $\RB(\GL_l)$ by $w_\theta$. For an integer $k$ with $1\le k\le l-1$, denote 
 $$\ov w_{k}=\bpm &I_{l-k}\\ I_k& \epm .$$
 
 \begin{lemma}\label{lem: bessel distance 1}
 For every $k$ with $1\le k\le l-1$, we have $\ov w_k=w_{\Delta-\wpair{\alpha_k}}.$
 \end{lemma}
 \begin{proof}
 We have $$w_{0}\ov w_k=\bpm J_k &\\ &J_{l-k} \epm,$$
 which is the longest Weyl element of the Levi subgroup 
 $$M_{\ov w_k}=\wpair{\bpm a&\\ &b \epm: a\in \GL_k(F), b\in \GL_{l-k}(F) }.$$
 The set of simple roots in $M_{\ov w_k}$ is $\Delta-\wpair{\alpha_k}$. Thus we have $\ov w_k\in \RB(\GL_l)$ and $\theta_{\ov w_k}=\Delta-\wpair{\alpha_k}$.
 \end{proof}

Denote $$\wt w_{n,m}=\bpm &&I_n\\ &I_{l-m-n}& \\ I_m&& \epm.$$
\begin{lemma}\label{lemma: bessel distance 2}
For positive integers $m,n$ with $1\le m+n\le l-1$, we have $\theta_{\wt w_{n,m}}=\Delta-\wpair{\alpha_m, \alpha_{l-n}}$.
\end{lemma}
\begin{proof}
We have $$w_0 \wt w_{n,m}=\begin{pmatrix}J_m &&\\ &J_{l-m-n}&\\&& J_n \end{pmatrix},$$
which is the longest Weyl element in the Levi subgroup 
$$M_{\wt w_{n,m}}=\begin{pmatrix}a&&\\ &b&\\ &&c \end{pmatrix},a\in \GL_n, b\in \GL_{l-m-n}, c\in \GL_m.$$
Thus $\theta_{\wt w_{n,m}}=\Delta-\wpair{\alpha_m, \alpha_{l-n}}.$
\end{proof}

Given $w,w'\in \RB(\GL_l)$ with $w>w'$, define (following Jacquet \cite{J})
$$d_B(w,w')=\max\wpair{m\mid \textrm{ there exist }w'_i\in \RB(\GL_l) \textrm{ with } w=w'_m>w'_{m-1}>\dots>w'_0=w'}.$$
The number $d_B(w,w')$ is called the Bessel distance of $w,w'$. By \cite{CPSS05}*{Proposition 2.1} and Lemma \ref{lem: bessel distance 1}, the set of elements in $\RB(\GL_l)$ which has Bessel distance 1 with the element $e\in \RB(\GL_l)$ are $\wpair{\ov w_k, 1\le k\le l-1}$, i.e.,
\begin{equation}\label{eq: Bessel distance 1}\wpair{w| d_B(w,e)=1}=\wpair{\ov w_k| 1\le k\le l-1}.\end{equation}
For $w,w'\in \bW$ with $w<w'$, we denote by $[w,w']$ the closed Bruhat interval $\wpair{w''\in \bW(\GL_l)| w\le w''\le w'}$.

\subsection{Cogdell-Shahidi-Tsai's theory on partial Bessel functions} In this subsection, we review certain basic properties of partial Bessel functions developed by Cogdell-Shahidi-Tsai recently in \cite{CST}.

For $w\in \RB(\GL_l)$, we denote 
\begin{equation}\label{eq: Aw}A_w=\wpair{a\in T_l(F)| \alpha(a)=1 \textrm{ for all }\alpha\in \theta_w}.\end{equation}
The set $A_w$ is in fact the center of $M_w$.
\begin{theorem}[Cogdell-Shahidi-Tsai] \label{theorem: CST} Let $\omega$ be a character of $F^\times$.
\begin{enumerate}
\item Let $w\in \bW$, $m>0$ and $f\in C_c^\infty(\Omega_w,\omega)$. Suppose $\CB_i(wa,f)=0$ for all $a\in A_w$. Then there exists $f_0\in C_c^\infty(\Omega_w-C(w), \omega)$, such that for sufficiently large $i$ depending only on $f$, we have $\CB_i(g,f)=\CB_i(g,f_0)$ for all $g\in \GL_l(F)$. 
\item Let $w\in \RB(\GL_l)$. Let $\Omega_{w,0}$ and $\Omega_{w,1}$ be $N_l\times N_l$ and $T_l$-invariant open sets of $\Omega_w$ such that $\Omega_{w,0}\subset \Omega_{w,1}$ and $\Omega_{w,1}-\Omega_{w,0}$ is a union of Bruhat cells $C(w')$ such that $w'$ does not support a Bessel function, i.e., $w'\notin \RB(\GL_l)$. Then for any $f_1\in C_c^\infty(\Omega_{w,1}, \omega)$ there exists $f_0\in C_c^\infty(\Omega_{w,0},\omega)$ such that for all sufficiently large $i$ depending only on $f_1$, we have 
$\CB_i(g,f_0)=\CB_i(g,f_1),$ for all $g\in \GL_l(F)$.
\end{enumerate}
\end{theorem}
\begin{proof}
Part (1) is \cite{CST}*{Lemma 5.13} and part (2) is \cite{CST}*{Lemma 5.14}. \end{proof}

\begin{corollary}\label{cor: germ expansion at level 0}
Let $f_1,f_2\in C_c^\infty(\GL_l(F),\omega)$ with $W^{f_1}(I_l)=W^{f_2}(I_l)=1$. Then there exist functions $f_{\ov w_k}\in C_c^\infty(\Omega_{\ov w_k},\omega)$ for all $k$ with $1\le k\le l-1$ such that for sufficiently large $i$ (depending only on $f_1,f_2$) we have 
$$\CB_i(g,f_1)-\CB_i(g,f_2)=\sum_{k=1}^{l-1} \CB_i(g,f_{\ov w_k}),\quad \forall g\in G.$$
\end{corollary}
This is essentially \cite{CST}*{Proposition 5.3}, see \cite{CST}*{page 2115} for a similar identity. Almost identical proofs in similar situations are given in \cite{Sp(2r)}*{Corollary 4.7} and \cite{U(2r+1)}*{Corollary 2.7}. We omit the proof here and just remark that each term in the expansion of the right side comes from the Weyl elements which has Bessel distance 1 from the trivial Weyl element $e\in \bW(\GL_l)$, namely the elements in the set \eqref{eq: Bessel distance 1}.

\subsection{Construction of certain sections of induced representations}\label{subsection: construction of section} Let $m,n$ be two positive integers and $\tau_1$ (resp. $\tau_2$) be an irreducible generic representation of $\GL_m(F)$ (resp. $\GL_n(F)$) and let $\Bs=(s_1,s_2)$. Consider 
$$N_{m,n}=\wpair{u_{m,n}(x)=\bpm I_m &x\\ &I_n \epm, x\in \Mat_{m\times n}}, \ov N_{m,n}=\wpair{\ov u_{m,n}(x):=\bpm I_m &\\ x&I_n \epm, x\in \Mat_{n\times m}},$$
and 
$$\ov N_{m,n}^k=\wpair{\ov u_{m,n}(x)\left| \bpm I_m &&\\ &I_{l-m-n} &\\ x&&I_n \epm\in H_l^k\right.}.$$
Here we identify $N_{m,n}$ etc. with its $F$-rational points and recall that $H_l^k$ is defined in Section \ref{subsection: Howe vector}.

Let $D$ be a compact open subset of $N_{m,n}$. For $x\in D$ and a positive integer $i$, we consider the set 
$$S(x,k)=\wpair{\ov y\in \ov N_{m,n}: \ov y x\in P_{m,n}\ov N_{m,n}^k}.$$

\begin{lemma}\label{lemma: preparation for construction of section}
\begin{enumerate}
\item For any positive integer $c$, there exists a positive integer $k_1=k_1(D,c)$ such that for all $k\ge k_1, x\in D, \ov y\in S(x,k)$, we can write 
$$\ov y x=u_{m,n}(x_1)\diag(a,b) \ov u_{m,n}(y_1),$$
with $a\in K^c_{\GL_m}, b\in K_{\GL_n}^c$. Here $u_{m,n}(x_1)\in N_{m,n}, \ov u_{m,n}(y_1)\in \ov N_{m,n}^k$. We recall that $K^c_{\GL_m}=I_m+\Mat_{m\times m}(\fp^c)$.
\item There exists an integer $k_2=k_2(D)$ such that $S(x,k)=\ov N_{m,n}^k$ for all $x\in D$ and $k\ge k_2$.
\end{enumerate}
\end{lemma}
\begin{proof}
This is an analogue of  \cite{Baruch}*{Lemma 4.1}, \cite{Sp(2r)}*{Lemma 5.1} and the proof is also similar. We provide a sketch below. For $x\in D$ and $\ov y\in S(x,k)$, we assume that $\ov y x= u_{m,n}(x_1)\diag(a,b) \ov u_{m,n}(y_1)$ for some $a\in \GL_m(F), b\in \GL_n(F), x_1\in \Mat_{m\times n}, y_1\in \Mat_{n\times m}$ with $\ov u_{m,n}(y_1)\in \ov N_{m,n}^k$. By abuse of notation, we also write $\ov y=\ov u_{m,n}(y), x=u_{m,n}(x)$. Then from the equation
$$\ov y^{-1} u_{m,n}(x_1)\diag(a,b)=x\ov u_{m,n}(-y_1), $$
we get
\begin{equation}\label{eq: a matrix equation in the construction of section}\bpm a& x_1b\\ -ya& (I_n-yx_1)b \epm=\bpm I_m-xy_1 &x\\ -y_1&I_n \epm.\end{equation}
We can solve that $a=I_m-xy_1$ and $b=I_n+y_1a^{-1}x$. Since when $x\in D$, the entries of $x$ are bounded, and the entries of $y_1$ go to zero as $k\to \infty$, we can take $k$ large enough such that $a=I_m-xy_1\in K_{\GL_m}^c$ and $ b=I_n+y_1a^{-1}x\in K_{\GL_n}^c.$ This proves (1).

By \eqref{eq: a matrix equation in the construction of section}, we have $y=y_1a^{-1} =y_1(I_m-xy_1)^{-1}=y_1(I_m+xy_1+(xy_1)^2+\dots)$. Again, since each entry of $x$ is bounded, we may take $k$ large such that the entries of $y_1(xy_1)^t$ are so large so that $\ov u_{m,n}(y_1(xy_1)^t)\in \ov N_{m,n}^k$  for $t\ge 0$. This shows that for $k$ large, we have $\ov u_{m,n}(y)\in \ov N_{m,n}^k$ and thus $S(x,k)\subset \ov N_{m,n}^k$ since $\ov y=\ov u_{m,n}(y)$ is arbitrarily chosen. See \cite{Sp(2r)}*{Lemma 5.1} for a similar and more detailed argument.

Take $x\in D$, we need to show $ \ov N_{m,n}^k \subset S(x,k)$ for $k$ large. As above, we write $x=u_{m,n}(x)$ by abuse of notation. We first assume that $k$ is so large such that if $\ov u_{m,n}(y)\in \ov N_{m,n}^k$, then $I_n+yx$ is invertible and $I_n-x(I_n+yx)^{-1}y$ is also invertible. This can be done because $x$ has bounded entries and $y$ has small entries if $\ov u_{m,n}(y)\in \ov N_{m,n}^k$ when $k$ large. Then we have $$\ov u_{m,n}(y)u_{m,n}(x)=u_{m,n}(x_1)\diag(a,b)\ov u_{m,n}(y_1),$$
with $ b=I_n+yx, a=I_n-b^{-1}y, x_1=xb^{-1}$ and $y_1=(I_n+yx)^{-1}y$. In particular, $ \ov u_{m,n}(y)u_{m,n}(x)\in P_{m,n}\ov N_{m,n}$. To show $ \ov u_{m,n}(y)\in S(x,k)$ for $k$ large, it suffices to show that one can choose $k$ large so that the above $ \ov u_{m,n}(y_1)\in \ov N_{m,n}^k$. Notice that $y_1=(I_n+yx)^{-1}y $ with bounded entries in $x$ and small entries in $y$, the argument is the same the above step. We are done.
\end{proof}

 Given $v_j\in V_{\tau_j}$, the space of $\tau_j$, for $j=1,2$, we consider the following $\tau_1\boxtimes \tau_2$-valued function on $\GL_{m+n}(F)$.
$$f_{\Bs}^{k,v_1,v_2}(g)=\left\{\begin{array}{ll}|\det(a)|^{s_1+\frac{n-1}{2}}|\det(b)|^{-s_2-\frac{m-1}{2}}\tau_1(a)v_1\boxtimes \tau_2(b)v_2, & \textrm{~if~} g=u_{m,n}(x)\diag(a,b)\ov u_{m,n}(y)\\
& \textrm{~with~}\ov u_{m,n}(y)\in \ov N_{m,n}^k,\\
0, & \textrm{ otherwise}.  \end{array}\right.$$
\begin{proposition}\label{proposition: definition of section}
For any $v_1,v_2$, there exists an integer $k_3(v_1,v_2)$ such that $f_\Bs^{k,v_1,v_2}$ defines a section in $\RI(\Bs,(\tau_1,\tau_2))$ for any $k\ge k_3(v_1,v_2)$.
\end{proposition}
\begin{proof}
This is an analogue of \cite{Sp(2r)}*{Lemma 5.2} and we only give a sketch of the proof. We first take a positive integer $c=c(v_1,v_2)$ such that $v_1$ is fixed by $K_{\GL_m}^c$ under the action of $\tau_1$ and $v_2$ is fixed by $K_{\GL_n}^c$ under the action of $\tau_2$. Now take $$k_3(v_1,v_2)=\max\wpair{c, k_1(K_{\GL_{m+n}}^c\cap N_{m,n},c), k_2(K_{\GL_{m+n}}^c\cap N_{m,n})}.$$
 For $k\ge k_3(v_1,v_2),$ we need to check 
\begin{equation}\label{eq: quasi-invariance of section}f_\Bs^{k,v_1,v_2}(u_{m,n}(x)\diag(a,b)g)=|\det(a)|^{s_1+\frac{n-1}{2}}|\det(b)|^{-s_2-\frac{m-1}{2}}\tau_1(a)\boxtimes \tau_2(b)f_\Bs^{k,v_1,v_2}(g),\end{equation}
for all $x\in \Mat_{m\times n}(F)$, $a\in \GL_m(F), b\in \GL_n(F), g\in \GL_{m+n}(F)$, and there exists an open compact subgroup $K'\subset \GL_{m+n}(F)$ such that \begin{equation}\label{eq: smoothness of section}f_\Bs^{k,v_1,v_2}(gh)=f_\Bs^{k,v_1,v_2}(g), \forall g\in \GL_{m+n}(F), h\in K'. \end{equation}
The first property \eqref{eq: quasi-invariance of section} is from the definition and we only address the second one \eqref{eq: smoothness of section}. 

Take a positive integer $t\ge k$ such that $\ov N_{m,n}\cap K_{\GL_{m+n}}^t\subset \ov N_{m,n}^k$. We take $K'=K_{\GL_{m+n}}^t$ in \eqref{eq: smoothness of section}. We have the decomposition
$$K_{\GL_{m+n}}^t=(K_{\GL_{m+n}}^t\cap N_{m,n})(K_{\GL_{m+n}}^t\cap M_{m,n})(K_{\GL_{m+n}}^t\cap \ov N_{m,n}).$$
For $h\in (K_{\GL_{m+n}}^t\cap \ov N_{m,n})$, we have $f_\Bs^{k,v_1,v_2}(gh)=f_\Bs^{k,v_1,v_2}(g) $ since $ h\in \ov N_{m,n}^k$ by assumption on $t$. For $h\in (K_{\GL_{m+n}}^t\cap  M_{m,n})$, we write $h=\diag(a_0,b_0)$. We first notice that $h^{-1}\ov N_{m,n}^k h\subset \ov N_{m,n}^k$, and thus $f_{\Bs}^{k,v_1,v_2}(g)=0$ if and only if $f_{\Bs}^{k,v_1,v_2}(gh)=0$. Next, we assume that $g=u_{m,n}(x)\diag(a,b)\ov u_{m,n}(y)$ with $ \ov u_{m,n}(y)\in \ov N_{m,n}^k$. Then $gh=u_{m,n}(x)\diag(aa_0,bb_0)\ov u_{m,n}(b_0^{-1}ya_0)$. Thus
\begin{align*}f_{\Bs}^{k,v_1,v_2}(gh)&=|\det(aa_0)|^{s_1+\frac{n-1}{2}}|\det(bb_0)|^{-s_2-\frac{m-1}{2}}\tau_1(aa_0)v_1\boxtimes \tau_2(bb_0)v_2\\
&=f_{\Bs}^{k,v_1,v_2}(g),
\end{align*}
where in the last step we used $\det(a_0)=\det(b_0)=1$ and $\tau_1(a_0)v_1=v_1,\tau_2(b_0)v_2=v_2$ (because $a_0\in K_{\GL_m}^t\subset K_{\GL_m}^c$ by the assumption $t\ge k\ge c$).
Finally, we take $h\in (K_{\GL_{m+n}}^t\cap N_{m,n})\subset K_{\GL_{m+n}}^c\cap N_{m,n}.$ Thus by Lemma \ref{lemma: preparation for construction of section}, we have $S(h,k)=S(h^{-1},k)=\ov N_{m,n}^k.$ In particular, for $ \ov u_{m,n}(y)\in \ov N_{m,n}^k$, we have $ \ov u_{m,n}(y) h\in P_{m,n}\ov N_{m,n}^k$ and $ \ov u_{m,n}(y) h^{-1}\in P_{m,n}\ov N_{m,n}^k$. Thus $f_{\Bs}^{k,v_1,v_2}(g)=0$ if and only if $f_{\Bs}^{k,v_1,v_2}(gh)=0$. Moreover, by Lemma \ref{lemma: preparation for construction of section} (1), we can write $\ov u_{m,n}(y) h=u_{m,n}(x_1)\diag(a_1,b_1)\ov u_{m,n}(y_1) $ with $a_1\in K_{\GL_m}^c, b_1\in K_{\GL_n}^c$. Thus for $g=u_{m,n}(x)\diag(a,b)\ov u_{m,n}(y)$, we have
$$gh=u_{m,n}(x)\diag(a,b)\ov u_{m,n}(y)h=u_{m,n}(x+ax_1b^{-1})\diag(aa_1,bb_1)\ov u_{m,n}(y_1).$$
From the definition, we see that $ f_{\Bs}^{k,v_1,v_2}(gh)=f_{\Bs}^{k,v_1,v_2}(g)$ because $\det(a_1)=\det(b_1)=1$, $\tau_1(a_1)v_1=v_1,$ and $\tau_2(b_1)v_2=v_2.$ This concludes the proof.
\end{proof}

We also consider the action of the intertwining operator $M_{w_{m,n}}$ on $f_{\Bs}^{i,v_1,v_2}$:
$$\wt f_{1-\wh \Bs}^{k,v_1,v_2}(g):=M_{w_{m,n}}(f_\Bs^{k,v_1,v_2})(g)=\int_{N_{n,m}(F)}f_\Bs^{k,v_1,v_2} (w_{m,n}ug)du.$$
\begin{lemma}\label{lemma: intertwining operator on section}
Let $D$ be an open compact subset of $N_{m,n}$. Then there is an integer $k_0(D,v_1,v_2)\ge k_3(v_1,v_2)$ such that for any $k\ge k_0(D,v_1,v_2)$, we have
$$\wt f_{1-\wh \Bs}^{k,v_1,v_2}(w_{m,n}^{-1}x)=\vol(\ov N_{m,n}^k)v_2\boxtimes v_1.$$
\end{lemma}
\begin{proof}
We take $c$ to be a common conductor of $v_1$ and $v_2$ (namely, $v_1$ is fixed by $\tau_1(K_{\GL_m}^c)$ and $v_2$ is fixed by $\tau_2(K_{\GL_n}^c)$) and we take  $k_0(D,v_1,v_2)=\max\wpair{k_3(v_1,v_2), k_1(D,c), k_2(D)}$. Assume $k\ge k_0(D,v_1,v_2)$. Then we have $S(x,k)=\ov N_{m,n}^k$ by Lemma \ref{lemma: preparation for construction of section}.
By definition, $$\wt f_{1-\wh \Bs}^{k,v_1,v_2}(w_{m,n}^{-1}x)=M_{w_{m,n}}(f_\Bs^{k,v_1,v_2})(g)=\int_{N_{n,m}(F)}f_\Bs^{k,v_1,v_2} (w_{m,n}uw_{m,n}^{-1}x)du.$$
For $u\in N_{n,m}$, we have $ \ov u:=w_{m,n}uw_{m,n}^{-1}\in \ov N_{m,n}.$ By the definition of $f_\Bs^{k,v_1,v_2} $, we have $f_\Bs^{k,v_1,v_2} (\ov u x)\ne 0$ if and only if $\ov u x\in P_{m,n}\ov N_{m,n}^k$ if and only if $\ov u\in S(x,k)=\ov N_{m,n}^k.$  Moreover, by Lemma \ref{lemma: preparation for construction of section} (1), we have 
$$\ov u x=u_{m,n}(x_1)\diag(a_1,b_1)\ov u_{m,n}(y_1),$$
with $x_1\in \Mat_{m\times n}(F), \ov u_{m,n}(y_1)\in \ov N_{m,n}^k$, $a_1\in K_{\GL_m}^c, b_1\in K_{\GL_n}^c$. By definition, we have
$$\wt f_{1-\wh \Bs}^{k,v_1,v_2}(w_{m,n}^{-1}x)=\vol(N_{m,n}^k)v_2\boxtimes v_1.$$
Here we recall that there is an extra operator that switches $v_1\boxtimes v_2$ to $v_2\boxtimes v_1$ in the definition of intertwining operators. This finishes the proof.
\end{proof}

In the above lemma, notice that $w_{m,n}^{-1}=w_{n,m}$. As we did in Subsection \ref{subsection: defn of local zeta integral}, we can consider the corresponding $\C$-valued function: $\xi_{\Bs}^{k,v_1,v_2}=\xi_{f_{\Bs}^{k,v_1,v_2} }\in \CW(\Bs,(\tau_1,\tau_2),\psi^{-1})$ and $\wt \xi_{1-\wh \Bs}= \xi_{\wt f_{1-\wh \Bs}^{k,v_1,v_2}}\in \CW(1-\wh \Bs, (\tau_2,\tau_1),\psi^{-1})$. By Lemma \ref{lemma: intertwining operator on section}, for $x\in D$ and $k\ge k_0(D,v_1,v_2)$, we have 
\begin{align}\label{eq: right side of the section}
\wt \xi_{1-\wh \Bs}^{k,v_1,v_2}(u_{n,m}(x_1)\diag(b,a)w_{n,m} x)=&\vol(\ov N_{m,n}^k)|\det(b)|^{1-s_2+\frac{m-1}{2}}|\det(a)|^{-(1-s_1)-\frac{n-1}{2}}\\
& \qquad W_{v_1}(a)W_{v_2}(b),\nonumber
\end{align}
for $x_1\in \Mat_{n\times m}(F), a\in \GL_m(F), b\in \GL_n(F).$ Here $W_{v_1}(a)=\lambda_1(\tau_1(a)v_1)$ for a fixed $\lambda_1\in \Hom_{N_m}(\tau_1,\psi^{-1})$ as in Subsection \ref{subsection: defn of local zeta integral}, and $W_{v_2}$ is defined similarly. Notice that $W_{v_1}\in \CW(\tau_1,\psi^{-1})$ and $W_{v_2}\in \CW(\tau_2,\psi^{-1}).$

\subsection{A result of Jacquet-Shalika}
\begin{proposition}\label{proposition: Jacquet-Shalika}
Let $W'$ be a smooth function on $\GL_n(F)$ which satisfies $W'(ug)=\psi(u)W'(g)$ for all $u\in N_n$ and for each $m$, the set $\wpair{g\in \GL_n(F)| W'(g)\ne 0, |\det(g)|=q^m} $ is compact modulo $U_{\GL_n}$. Assume, for all irreducible generic representation $\tau$ of $\GL_n(F)$ and for all Whittaker functions $W\in \CW(\tau,\psi^{-1})$, the following integral 
$$\int_{U_{\GL_n}\setminus \GL_n}W'(g)W(g)|\det(g)|^{s-k}dg$$
vanishes, where $k$ is a fixed number, then $W'\equiv 0$.
\end{proposition}

This is a corollary of \cite{Jacquet-Shalika-generic}*{Lemma 3.2}. The version as in Proposition \ref{proposition: Jacquet-Shalika} appeared in \cite{Henniart}*{Lemme 2.1}. See also \cite{chen}*{Corollary 2.1} or \cite{Baruch}*{Lemma 5.2} for a proof of the current version.

\section{Proof of the local converse theorem}\label{section-proof}

In this section, we prove Theorem \ref{theorem: main}. Let $\pi_1, \pi_2$ be two irreducible supercuspidal representations of $\GL_l(F)$ sharing the same central character $\omega$. For $j=1,2$, fix $f_j\in \CM(\pi_j)$ such that $W^{f_j}(I_l)=1$.

\begin{theorem}\label{theorem: inductive}
Let $m$ be an integer with $0\le m\le [l/2]$. Condition $\CC(m)$ implies that for each $j$ with $m+1\le j\le l-1-m$, there exist functions $f_{\ov w_j}\in C_c^\infty(\Omega_{\ov w_j}, \omega)$ such that for all $g\in \GL_l(F)$ and for all sufficiently large $i$ (depending only on $f_1, f_2$), we have
$$\CB_i(g,f_1)-\CB_i(g,f_2)=\sum_{j=m+1}^{l-m-1}\CB_i(g,f_{\ov w_j}).$$
\end{theorem}
We first show that Theorem \ref{theorem: inductive} implies Theorem \ref{theorem: main}.
\begin{proof}[Theorem \ref{theorem: inductive} implies Theorem \ref{theorem: main}]
Assume Theorem \ref{theorem: inductive}. Then condition $\mathcal{C}([l/2])$ gives $\mathcal{B}_i(g, f_1) = \mathcal{B}_i(g, f_2)$ for all $g \in \operatorname{GL}_l(F)$ and all sufficiently large $i$. Hence $W_i^{f_1} = W_i^{f_2}$, so $\mathcal{W}(\pi_1, \psi) \cap \mathcal{W}(\pi_2, \psi) \neq \emptyset$. Uniqueness of the Whittaker model then forces $\pi_1 \cong \pi_2$.
\end{proof}

\begin{remark}
For classical group analogues of Theorem \ref{theorem: inductive}, see \cite{Sp(2r)}*{Proposition 6.1} and \cite{U(2r+1)}*{Theorem 4.1}. Theorem \ref{theorem: inductive} appears to be stronger than Theorem \ref{theorem: main}. We expect it to be useful for the following question: given an integer $t$ with $t \le [l/2]$, classify all irreducible supercuspidal representations $\pi$ of $\operatorname{GL}_l(F)$ that are determined by the collection of gamma factors $\gamma(s, \pi \times \tau, \psi)$ as $\tau$ ranges over all irreducible generic representations of $\operatorname{GL}_m(F)$ with $1 \le m \le t$. 
\end{remark}

 We prove Theorem \ref{theorem: inductive} by induction. The base case $m=0$ of Theorem \ref{theorem: inductive} is precisely Corollary \ref{cor: germ expansion at level 0}. For the inductive step, we assume the following

\begin{indhypo}\label{indhypo1} Let $m$ be a positive integer with $m\le [l/2]$. Condition $\CC(m-1)$ implies that
there exist functions $f_{\ov w_j}\in C_c^\infty(\Omega_{\ov w_j}, \omega)$ for each $j$ with $m\le j\le l-m$ such that, 
\begin{equation}\label{eq: germ expansion at level m-1}\CB_i(g,f_1)-\CB_i(g,f_2)=\sum_{j=m}^{l-m}\CB_i(g,f_{\ov w_j}),
\end{equation}
for all $g\in \GL_l(F)$ and all sufficiently large $i$ depending only on $f_1, f_2$. 
\end{indhypo}

 Assuming the inductive hypothesis above, we employ a second inductive argument to show that condition $\mathcal{C}(m)$ implies the existence of functions $f_{\overline{w}_j} \in C_c^\infty(\Omega_{\overline{w}_j}, \omega)$ for each $j$ with $m+1 \le j \le l-1-m$ such that
 \begin{equation}\label{eq: germ expansion at level m}
 \mathcal{B}_i(g, f_1) - \mathcal{B}_i(g, f_2) = \sum_{j=m+1}^{l-m-1} \mathcal{B}_i(g, f_{\overline{w}_j})
 \end{equation}
holds for all sufficiently large $i$ (depending only on $f_1, f_2$) and for all $g \in \operatorname{GL}_l(F)$. The functions $f_{\overline{w}_j}$ appearing here may differ from those obtained from the $(m-1)$-th step in \eqref{eq: germ expansion at level m-1}; we do not distinguish them notationally.

To proceed with another induction argument, we introduce a refined condition $\mathcal{C}(m, n)$ for integers $0 \le n \le m$. \textit{We say that $\pi_1, \pi_2$ satisfy $\mathcal{C}(m, n)$ if they satisfy $\mathcal{C}(m-1)$ and additionally
 \[
 \Gamma(\mathbf{s}, \pi_1 \times (\tau_1, \tau_2), \psi) = \Gamma(\mathbf{s}, \pi_2 \times (\tau_1, \tau_2), \psi)
\]
 holds for all irreducible generic representations $\tau_1$ of $\operatorname{GL}_m(F)$ and $\tau_2$ of $\operatorname{GL}_k(F)$ with $0 \le k \le n$, as well as for all irreducible generic representations $\tau_2$ of $\operatorname{GL}_m(F)$ and $\tau_1$ of $\operatorname{GL}_k(F)$ with $0 \le k \le n$.} Note that condition $\CC(m,0)$ is strictly stronger than $\CC(m-1)$, while condition $\CC(m,m)$ coincides with $\CC(m)$.
 
 Recall that for any positive integers $j, m$ with $m+j<l$, we have defined an element 
$$\wt w_{j,m}=\bpm &&I_j\\ &I_{l-m-j}& \\ I_m && \epm$$
in \S\ref{subsection: Weyl elements which support Bessel functions}. Moreover, we know that $\wt w_{j,m}\in \RB(\GL_l)$ and $\theta_{\wt w_{j,m}}=\Delta-\wpair{\alpha_m,\alpha_{l-j}}$ by Lemma \ref{lemma: bessel distance 2}.

\begin{theorem}\label{theorem: second induction}
Let $m$ be a positive integer with $m\le [l/2]$, and let $n$ be an integer with $0\le n\le m$. Then condition $\CC(m,n)$ implies that there exist functions 
\begin{itemize}
\item $f_{\ov w_j}\in C_c^\infty(\Omega_{\ov w_j},\omega)$ for each $j$ with $m+1\le j\le l-m-1$;
\item $f_{j,m}'\in C_c^\infty(\Omega_{\wt w_{j,m}},\omega) $ for each $j$ with $n+1\le j\le m$; and
\item $f_{m,j}''\in C_c^\infty(\Omega_{\wt w_{m,j}},\omega)$ for each $j$ with $n+1\le j\le m$,
\end{itemize}
such that
\begin{equation}\label{eq: germ expansion at level mn+1}\CB_i(g,f_1)-\CB_i(g,f_2)=\sum_{j=m+1}^{l-m-1}\CB_i(g,f_{\ov w_j})+\sum_{j=n+1}^{m}\CB_i(g,f_{j,m}')+\sum_{j=n+1}^m \CB_i(g,f_{m,j}''),\end{equation}
for all $g\in \GL_{l}(F)$ and for all sufficiently large $i$ (depending only on $f_1,f_2$).
\end{theorem}

\begin{remark}\label{remark: the case when m=n}
If $n=m-1$, both $ f_{m,m}'$ and $f_{m,m}''$ are in $C_c^\infty(\Omega_{\wt w_{m,m}},\omega)$. By incorporating $f_{m,m}''$ into $f_{m,m}'$, the statement of Theorem \ref{theorem: second induction} reduces to the following: condition $\CC(m,m-1)$ implies the expansion 
$$ \CB_i(g,f_1)-\CB_i(g,f_2)=\sum_{j=m+1}^{l-m-1}\CB_i(g,f_{\ov w_j})+\CB_i(g,f'_{m,m}),$$
with certain $f_{\ov w_j}\in C_c^\infty(\Omega_{\ov w_j},\omega)$ and $f'_{m,m}\in C_c^\infty(\Omega_{\wt w_{m,m}}, \omega).$
\end{remark}

Note that by Theorem \ref{theorem: second induction}, condition $\CC(m,m)=\CC(m)$ implies that 
$$\CB_i(g,f_1)-\CB_i(g,f_2)=\sum_{j=m+1}^{l-m-1}\CB_i(g,f_{\ov w_j}), $$
which is exactly the desired expansion. Thus, Theorem \ref{theorem: second induction} implies Theorem \ref{theorem: inductive}, and hence Theorem \ref{theorem: main}. We prove Theorem \ref{theorem: second induction} in the remaining part of this section.

\subsection{Proof of the base case of Theorem \ref{theorem: second induction}}\label{subsection: base case}
In this subsection, we prove the base case of Theorem \ref{theorem: second induction}, namely, the case when $n=0$.

Let $k$ be a positive integer with $k<l$. Consider the parabolic subgroup $P_{k,l-k}$ of $\GL_l(F)$. Its Levi subgroup $M_{k,l-k}$ consists of block diagonal matrices, and a typical element is denoted by $$\bt_{k}(a,b):=\bpm a&\\ &b \epm, a\in \GL_k(F), b\in \GL_{l-k}(F).$$
For $y\in \Mat_{m\times (l-m-1)}(F),$
set $$u_1(y)=\bpm I_m &&y\\ &1&\\ &&I_{l-m-1} \epm.$$
\begin{lemma}\label{lemma: support} We fix the notations as in Inductive Hypothesis \ref{indhypo1}.
\begin{enumerate}
\item We have $\CB_i(h,f_{\ov w_j})=0, \forall h\in P_{k,l-k}$. In particular, the inductive hypothesis \eqref{eq: germ expansion at level m-1} implies that 
$$\CB_i(h,f_1)=\CB_i(h,f_2),$$
for all $h\in P_{k,l-k}$ and $i$ large.
\item For positive integer $j$ with $m+1\le j\le l-m$, we have
$$\CB_i(\ov w_m \bt_{m}(a,I_{l-m})u_1(y), f_{\ov w_j})=0, \forall a\in \GL_m(F), \forall y\in \Mat_{m\times (l-m-1)}(F).$$
In particular, the inductive hypothesis \eqref{eq: germ expansion at level m-1} implies that
\begin{align*}&\CB_i(\ov w_m \bt_{m}(a, I_{l-m})u_1(y),f_1)-\CB_i(\ov w_m \bt_{m}(a, I_{l-m})u_1(y),f_2)\\
&\qquad =\CB_i(\ov w_m \bt_{m}(a, I_{l-m})u_1(y),f_{\ov w_m}),\end{align*}
for all $a\in \GL_m(F), y\in \Mat_{m\times (l-m-1)}(F)$.
\item For any $a\in \GL_m(F)$, we can take $i$ large enough (which only depends on $f_{\ov w_m}$, and hence only on $f_1,f_2$), such that 
$$\CB_i(\ov w_m \bt_m(a,I_{l-m})u_1(y),f_{\ov w_m})=\left\{\begin{array}{ll} \CB_i(\ov w_m \bt_m(a,I_{l-m}),f_{\ov w_m}), & \mathrm{~if~} u_1(y)\in H_l^i,\\ 0, & \mathrm{otherwise}.  \end{array} \right. $$
\item For a fixed integer $k$ and $i$, the set $\wpair{a\in N_m(F)\bs \GL_m(F): \CB_i(\ov w_m \bt_m(a,I_{l-m}))\ne 0, |a|=q^k}$ is compact.
\end{enumerate}
\end{lemma}

\begin{proof}
(1) Recall that $$\CB_i(g,f_{\ov w_j})=\frac{1}{\vol(N_l^i)}\int_{N_l^i}\int_{N_l}f_{\ov w_j}(u_1gu_2)\psi^{-1}du_2du_1.$$
Since $\Supp(f_{\ov w_j})\subset \Omega_{\ov w_j}$, it suffices to show that $P_{k,l-k}\cap \Omega_{\ov w_j}=\emptyset$. Suppose that $P_{k,l-k}\cap \Omega_{\ov w_j}$ is not empty, then their intersection must contain a Bruhat cell, namely, there exists a $w\in \bW$ such that $w\ge \ov w_j$ and $C(w)\subset P_{k,l-k}$. Since $P_{k,l-k}$ is closed in $\GL_{l}$, we get $\ov {C(w)}\subset P_{k,l-k}$. The condition $w\ge \ov w_j$ implies that $C(\ov w_j)\subset \ov{C(w)}\subset P_{k,l-k}$. In particular, we have $\ov w_j\in P_{k,l-k}$. This is a contradiction.

(2) Consider the set 
\begin{align*}S&=\wpair{w\in \bW: w=\ov w_m \bt_{m}(a,I_{l-m}), \mathrm{ for ~ some~}a\in \GL_m}\\
&=\wpair{\ov w_m \bt_{m}(w',I_{l-m}): w'\in \bW(\GL_m)}.\end{align*}
Here we don't distinguish a Weyl element its rerepsentative.
 Denote $w_{\max}^m=\ov w_m \diag(J_m, I_{l-m})=\bpm &I_{l-m}\\ J_m& \epm$. Since the Weyl element in $\GL_m$ forms a Bruhat interval $[1,J_m]$, the set $S$ is in fact the Bruhat interval $[\ov w_m, w_{\max}^m]$. Since $$\wpair{\ov w_m \bt_{m}(a,I_{l-m})u_1(y), a\in \GL_m(F), y\in \Mat_{m\times (l-m-1)}(F)}\subset \cup_{w\in S} C(w),$$ it suffices to show that for any $w\in S$,  $C(w)\cap \Omega_{\ov w_j}=\emptyset$ if $m+1\le j\le l-m$. Suppose that $C(w)\cap \Omega_{\ov w_j}$ is non-empty, then $w\ge \ov w_i$. In particular, $ w_{\max}^m\ge \ov w_j$. Note that 
$$w_{0}w_{\max}=\bpm I_m&\\&J_{l-m} \epm,$$
which is the longest Weyl element of the Levi subgroup $$M_{w_{\max}^m}=\wpair{\diag(a_1,\dots,a_m, a): a_i\in \GL_1, a\in \GL_{l-m}}.$$ Note that the set $\theta_{w_{\max}^m}$ is the set of all Weyl elements in $M_{w_{\max}^m}$, which is $\Delta-\wpair{\alpha_1,\dots,\alpha_m}$. From $w_{\max}^m\ge \ov w_j$, we obtain $\theta_{w_{\max}^m}\subset \theta_{\ov w_j}$, namely, $\Delta-\wpair{\alpha_1,\dots,\alpha_m}\subset \Delta-\wpair{\alpha_j}$. This is impossible because $j>m$.

(3) This can be done using a root killing argument as in Lemma \ref{lemma-convergence-support-Ujmn}, or using a support argument as in \cite{Sp(2r)}*{Lemma 6.3 (3)}. Since the proof is similar to that of \cite{Sp(2r)}*{Lemma 6.3 (3)},  we omit the details.

(4) This is an analogue of \cite{Sp(2r)}*{Lemma 6.3 (4)} and the proof is similar. We omit the details.
\end{proof}

Notice that if $m>0, n=0$, we have defined a gamma factor $\Gamma(\Bs,\pi\times (\tau_1,0),\psi)$ for an irreducible generic representation $\tau_1$ of $\GL_m(F)$, which is just a shift of Jacquet--Piatetski-Shapiro--Shalika's local gamma factor. The pair notation $(\tau_1, 0)$ is used to maintain the form of a representation of $\GL_m(F)\times \GL_n(F)$ even when $n=0$ and hence the group $\GL_n(F)$ is trivial. See Remark \ref{remark: comparison with JPSS integral1}.

\begin{proposition}\label{prop: first step}
Condition $\CC(m,0)$ gives
\begin{equation}\label{eq: equality on bar wm}\CB_i(\ov w_m\bt_m(a,I_{l-m}),f_1)=\CB_i(\ov w_m \bt_m(a,I_{l-m}),f_2),\end{equation}
 and 
\begin{equation}\label{eq: equality on bar wl-m}\CB_i(\ov w_{l-m}\bt_{l-m}(I_{l-m},a),f_1)=\CB_i(\ov w_{l-m} \bt_m(I_{l-m},a),f_2) \end{equation}
for all $a\in \GL_m(F).$
\end{proposition}
This is roughly \cite{chen}*{Proposition 3.1}. Since our proof uses partial Bessel functions rather than the Kirillov model (as in \cite{chen}), we provide a detailed exposition below.
\begin{proof}
 For any irreducible generic representation $\tau_1$ of $\GL_m(F)$ and any $\xi_{\Bs}=W'|~|^{s-1/2}$ with $W'\in \CW(\tau_1,\ov \psi)$, we can consider the integral $\Psi(\rho(\gamma_{m,0}^{-1})\CB_i^f, \xi_{\Bs}; 0)$ for $f=f_1,f_2$, which is 
$$\Psi(\CB_i^f,\xi_\Bs;0)=\int_{N_m(F)\bs \GL_m(F)}\CB_i^f\left(\bt_m(a,I_{l-m})\right)  W'(a)|\det(a)|^{s-1/2}dh.$$
Here we notice that $\gamma_{m,0}=I_l$. See also Remark \ref{remark: comparison with JPSS integral1}. By inductive hypothesis \ref{indhypo1} and Lemma \ref{lemma: support} (1), we have
$$\CB_i^{f_1}\left(\bt_m(a,I_{l-m})\right)=\CB_i^{f_2}\left(\bt_m(a,I_{l-m}) \right). $$
Thus 
$$\Psi(\CB_i^{f_1},  \xi_{\Bs}; 0)=\Psi(\CB_i^{f_2},  \xi_{\Bs}; 0).$$
By the assumption on local gamma factors and the local functional equation \eqref{eq: modified LFE}, we have
$$\Psi( \CB_i^{f_1}-\CB_i^{f_2}, M_{w_{m,n}}(\xi_{\Bs});0)=0.$$
Plugin the definitions, see \eqref{eq: right side local zeta integral} or Remark \ref{remark: comparison with JPSS integral1}, we have
\begin{align*}
0=&\int_{[\GL_m]}\int_{ \Mat_{m\times (l-m-1)}}\left(\CB_i^{f_1}-\CB_i^{f_2} \right)\left(\bpm 1 &&\\ &I_{l-m-1}&\\ &y&I_{m} \epm \bpm &I_{l-m}\\ I_m& \epm \bpm a&\\ &I_{l-m} \epm \right)\\ & \qquad \cdot W'(a)|\det(a)|^{s-1/2}dydh\\
=&\int_{[\GL_m]}\int_{ \Mat_{m\times (l-m-1)}}\left(\CB_i^{f_1}-\CB_i^{f_2} \right)\left(\ov w_m \bt_m(a, I_{l-m})u_1(y)\right) W'(a)|\det(a)|^{s^*}dydh,
\end{align*}
where we identify an algebraic group over $F$ with its $F$-rational points, $[\GL_m] $ is the abbreviation of $N_m(F)\bs \GL_m(F)$ and $s^*=s-\frac{1}{2}+l-m-1$.  
By Lemma \ref{lemma: support} (2) and (3), we get 
\begin{align*}
\int_{N_m(F)\bs \GL_m(F)}\left( \CB_i^{f_1}-\CB_i^{f_2} \right)(\ov w_m \bt_m(a,I_{l-m}))W'(a)|\det(a)|^{s^*}dh=0.
\end{align*}
Note that this is true for all irreducible representation $\tau_1$ of $\GL_m(F)$ and for all $W'\in \CW(\tau_1,\psi^{-1})$. Thus by Proposition \ref{proposition: Jacquet-Shalika} and Lemma \ref{lemma: support} (4), we get that 
$$ \CB_i(\ov w_m \bt_m(a, I_{l-m}) ,f_1)=\CB_i(\ov w_m \bt_m(a, I_{l-m}) ,f_2).$$
To get the second assertion, we need to use the local gamma factor $\Gamma(\Bs,\pi\times (0,\tau_2),\psi)$ for a generic representation $\tau_2$ of $\GL_m(F)$. Here $\Bs=s$ is a complex number used to do twist on $\tau_2$. The calcluation is almost identical to the above. In fact, if we take $\xi_{\Bs}=W'|~|^{s-1/2}$ with $W'\in \CW(\tau_2,\psi^{-1})$, we can check that 
$$\Psi(\rho(\gamma_{0,m}^{-1})\CB_i^f, \xi_{\Bs};0)=\int_{[\GL_m]}\int\CB_i^f\left( \bpm 1 &&\\ &I_{l-m-1}&\\ &y&I_m  \epm\bt_{l-m}(I_{l-m},a)\right)W'(a)|\det(a)|^{s-1/2}dyda.$$
By Lemma \ref{lemma: support} (1), we have $\Psi(\rho(\gamma_{0,m}^{-1})\CB_i^{f_1}, \xi_{\Bs};0)=\Psi(\rho(\gamma_{0,m}^{-1})\CB_i^{f_2}, \xi_{\Bs};0) $. By the local functional equation \eqref{eq: modified LFE}, we get that
$$\Psi( \CB_i^{f_1}-\CB_i^{f_2}, M_{w_{m,n}}(\xi_{\Bs});0)=0.$$
By \eqref{eq: right side local zeta integral}, the above equation becomes
$$\int_{[\GL_m]}(\CB_i^{f_1}-\CB_i^{f_2})(\ov w_{l-m}\bt_{l-m}(I_{l-m},a))|a|^{s^*}da=0,$$
where $s^*$ is a translation of $s$ and its precise form is not important here. Then using Proposition \ref{proposition: Jacquet-Shalika} again, we get that 
$$(\CB_i^{f_1}-\CB_i^{f_2})(\ov w_{l-m}\bt_{l-m}(I_{l-m},a))=0, \forall a\in \GL_m(F). $$
This finishes the proof.
\end{proof} 

\begin{corollary}\label{cor: base case}
Assume condition $\CC(m,0)$. Then there exists 
\begin{enumerate}
\item[$\bullet$] $f_{\ov w_j}\in C_c^\infty(\Omega_{\ov w_j},\omega)$ for each $j$ with $m+1\le j\le l-m-1$;
\item[$\bullet$] $f_{j,m}'\in C_c^\infty(\Omega_{\wt w_{j,m}},\omega) $, for each $j$ with $1\le j\le m$; and
\item[$\bullet$] $f_{m,j}''\in C_c^\infty(\Omega_{\wt w_{m,j}},\omega)$, for each $j$ with $1\le j\le m$,
\end{enumerate}
such that
$$\CB_i(g,f_1)-\CB_i(g,f_2)=\sum_{j=m+1}^{l-m-1}\CB_i(g,f_{\ov w_j})+\sum_{j=1}^{m}\CB_i(g,f_{j,m}')+\sum_{j=1}^m \CB_i(g,f_{m,j}''),$$
for all $g\in \GL_{l}(F)$ and for all $i$ large enough depending only on $f_1,f_2$.
\end{corollary}
\begin{proof}
 By Lemma \ref{lemma: support} (2), inductive hypothesis \eqref{eq: germ expansion at level m-1} and \eqref{eq: equality on bar wm}, we get 
\begin{equation} \label{eq: vanishing of fm} \CB_i(\ov w_m \bt_m(a,I_{l-m}), f_{\ov w_m})=0.\end{equation}
As in the proof of Lemma \ref{lemma: support} (2), we consider $w_{\max}^m=\bpm &I_{l-m}\\ J_m&\epm$. Then for $w\in [\ov w_m, w_{\max}^m]$, we consider the set $A_w$ as defined in \eqref{eq: Aw}. From $w\le w_{\max}^m$, we know that $A_w\subset A_{w_{\max}^m}$ which is of the form $\diag(a_1,\dots, a_m, a I_{l-m})$, for $a_j, a\in \GL_1(F)$.  Moreover, we know that $w=\ov w_m \bt_m(w',I_{l-m})$. Thus, for any $a\in A_w$, we know that there exists an element $z=zI_{l}$ in the center of $\GL_l(F)$ and an element $g\in \GL_m(F)$ such that $wa=z\ov w_m \bt_m(b, I_{l-m})$. Thus from \eqref{eq: vanishing of fm}, we get that
\begin{equation}\label{eq: vanishing on bar wm}\CB_i(wa, f_{\ov w_m})=0,\end{equation}
for all $w\in [\ov w_m, w_{\max}^m]$ and all $a\in A_w$. Similarly, if we temporarily denote $w'_{\max}=\ov w_{l-m}\diag(I_{l-m}, J_m)$, then from \eqref{eq: equality on bar wl-m} we have
 \begin{equation}\label{eq: vanishing on bar wl-m}\CB_i(wa,f_{\ov w_{l-m}})=0,\end{equation} for all $w\in \RB(\GL_l)$ with $\ov w_{l-m}\le w\le w'_{\max},$ and all $a\in A_w$. The result in fact follows from \eqref{eq: vanishing on bar wm}, \eqref{eq: vanishing on bar wl-m} and Theorem \ref{theorem: CST} directly. We give some details about this implication below.

By the proof of Lemma \ref{lemma: support} and a simple calculation, we get that 
$$\begin{array}{lllll}
&\theta_{\ov w_m}&=\Delta-\wpair{\alpha_m},  &\theta_{w_{\max}^m}&=\Delta-\wpair{\alpha_1,\dots,\alpha_m},\\
&\theta_{\ov w_{l-m}}&=\Delta-\wpair{\alpha_{l-m}}, & \theta_{w_{\max}'}&=\Delta-\wpair{\alpha_{l-m},\dots,\alpha_{l-1}}.
\end{array}$$
Denote $$\Omega_{\ov w_m}^\circ=\bigcup_{\substack{w\in \RB(\GL_l), w>\ov w_m\\ d(w,\ov w_m)=1}}\Omega_w.$$
By applying Theorem \ref{theorem: CST} and \eqref{eq: vanishing on bar wm} to $\ov w_m$, we get a function $\ov f_m\in C_c^\infty(\Omega_{\ov w_m}^\circ,\omega)$ such that, after increasing $i$ if necessary, we have
$$\CB_i(g,f_{\ov w_m})=\CB_i(g,\ov f_m).$$
Note that the set $\wpair{w\in \RB(\GL_l): w>\ov w_m, d(w,\ov w_m)=1}=\wpair{w_{\Delta-\wpair{\alpha_m, \alpha_j}}, 1\le j\le l-1, j\ne m}.$ By a partition of unity argument on $f_m$, there exists a function $f_{\Delta-\wpair{\alpha_j,\alpha_m}}\in C_c^\infty(\Omega_{w_{\Delta-\wpair{\alpha_m,\alpha_j}}},\omega)$ such that
\begin{equation}\label{eq: temp expansion}\CB_i(g,f_{\ov w_m})=\CB_i(g, \ov  f_m)=\sum_{j\ne m}\CB_i(g,f_{\Delta-\wpair{\alpha_j, \alpha_m}}).   \end{equation}
We consider $j$ in 3 separate ranges. If $m+1\le j\le l-m-1$,  since $ w_{\Delta-\wpair{\alpha_m,\alpha_j}}\ge \ov w_j$, $f_{ {\Delta-\wpair{\alpha_j, \alpha_m}}}$ can be viewed as an element of $C_c^\infty(\Omega_{\ov w_j},\omega)$ and thus can be absorbed into $f_{\ov w_j}$ in \eqref{eq: germ expansion at level m-1}. In other words, we can assume that $f_{ {\Delta-\wpair{\alpha_j, \alpha_m}}}=0 $ after replacing $f_{\ov w_j}$ by $f_{\ov w_j}+f_{ {\Delta-\wpair{\alpha_j, \alpha_m}}}$ in \eqref{eq: germ expansion at level m-1}.  If $l-1\ge j\ge l-m$, we have $f_{ {\Delta-\wpair{\alpha_j, \alpha_m}}}\in C_c^\infty (\Omega_{\wt w_{l-j,m}}, \omega)$. We write $f_{ {\Delta-\wpair{\alpha_j, \alpha_m}}}$ as $f_{\wt w_{l-j},m}'$.  Thus \eqref{eq: temp expansion} becomes
\begin{equation}\label{eq: temp expansion 2}\CB_i(g,f_{\ov w_m})=\CB_i(g, \ov f_m)=\sum_{j=1}^{m-1}\CB_i(g,f_{\Delta-\wpair{\alpha_j, \alpha_m}})+\sum_{j=1}^{m}\CB_i(g,f_{\wt w_{j,m}}').   \end{equation}
If $j<m$, then $\ov w_m\le w_{\Delta-\wpair{\alpha_m,\alpha_j}}\le w_{\max} $, the formula \eqref{eq: vanishing on bar wm} and the above decomposition of $f_{\ov w_m}$ \eqref{eq: temp expansion} imply that 
$$\CB(w a, f_{\Delta-\wpair{\alpha_j, \alpha_m}})=0, w=w_{\Delta-\wpair{\alpha_m,\alpha_j}},  a\in A_{w}.$$
We then apply Theorem \ref{theorem: CST} to $w=w_{\Delta-\wpair{\alpha_m,\alpha_j}}$ and repeat the above process. We can get that for each $k$ with $k\ne j, m$, there exists a function $f_{\Delta-\wpair{\alpha_j,\alpha_k,\alpha_m}}\in C_c^\infty(\Omega_{w_{\Delta-\wpair{\alpha_j,\alpha_k,\alpha_m}}},\omega)$ such that
$$\CB(g, f_{\Delta-\wpair{\alpha_j, \alpha_m}})=\sum_{k\ne j, m}\CB(g, f_{\Delta-\wpair{\alpha_j, \alpha_k,\alpha_m}}). $$
Similarly as above, if $m+1\le k\le l-m-1$, we can assume that $  f_{\Delta-\wpair{\alpha_j, \alpha_k,\alpha_m}}=0$ after replacing $f_{\ov w_k}$ in \eqref{eq: germ expansion at level m-1} by $f_{\ov w_k}+ f_{\Delta-\wpair{\alpha_j, \alpha_k,\alpha_m}}$. If $l-1\ge k\ge l-m$, we have $  f_{\Delta-\wpair{\alpha_j, \alpha_k,\alpha_m}}\in C_c^\infty(\Omega_{\wt w_{l-k,m}},\omega)$. We can thus absorb $  f_{\Delta-\wpair{\alpha_j, \alpha_k,\alpha_m}}$ to $f_{\wt w_{l-k,m}}'$ in \eqref{eq: temp expansion 2} and assume that $f_{\Delta-\wpair{\alpha_j, \alpha_k,\alpha_m}}=0 $. Then \eqref{eq: temp expansion 2} becomes
\begin{equation}\label{eq: temp expansion 3}
\CB_i(g,f_{\ov w_m})=\CB_i(g, \ov f_m)=\sum_{1\le j<k\le m-1}\CB_i(g,f_{\Delta-\wpair{\alpha_j, \alpha_k \alpha_m}})+\sum_{j=1}^{m}\CB_i(g,f_{\wt w_{j,m}}')
\end{equation}
We continue to repeat the above process. In each time, we increase $i$ if necessary, and replacing $f_{\ov w_j}$ for $m+1\le j\le l-m-1$ in \eqref{eq: germ expansion at level m-1} and $f_{\wt w_{j,m}}'$  in \eqref{eq: temp expansion 2} by a new function in the same corresponding space if necessary. After repeating the above process at most $m$-times, we can get 
\begin{equation}\label{eq: temp expansion for bar wm} \CB_i(g,f_{\ov w_m})=\CB_i(g, \ov f_m)=\sum_{j=1}^{m}\CB_i(g,f_{\wt w_{j,m}}'), f_{\wt w_{j,m}}\in C_c^\infty(\Omega_{\wt w_{j,m}},\omega).\end{equation}

Similarly, using \eqref{eq: vanishing on bar wl-m} and Theorem \ref{theorem: CST}, there exists functions $f_{\wt w_{m,j}''}\in C_c^\infty(\Omega_{\wt w_{m,j}},\omega)$ such that 
\begin{equation}\label{eq: temp expansion for bar wl-m}\CB_i(g,f_{\ov w_{l-m}})=\sum_{j=1}^m\CB_i(g,f_{\wt w_{m,j}}'').\end{equation}
Now the result follows from the inductive hypothesis \eqref{eq: germ expansion at level m-1}, equations \eqref{eq: temp expansion for bar wm} and \eqref{eq: temp expansion for bar wl-m}.
\end{proof}

\subsection{Proof of Theorem \ref{theorem: second induction}}\label{subsection: proof in the odd case}
Note that Corollary \ref{cor: base case} gives the base case of Theorem \ref{theorem: second induction}. Given a positive integer $n$ with $1\le n\le m$, we assume that we have proved Theorem \ref{theorem: second induction} for $n-1$, namely, we assume the following

\begin{indhypo}\label{indhypo2} Condition $\CC(m,n-1)$ implies that there exist functions
\begin{enumerate}
\item[$\bullet$] $f_{\ov w_j}\in C_c^\infty(\Omega_{\ov w_j},\omega)$ for each $j$ with $m+1\le j\le l-m-1$;
\item[$\bullet$] $f_{j,m}'\in C_c^\infty(\Omega_{\wt w_{j,m}},\omega) $, for each $j$ with $n\le j\le m$; and
\item[$\bullet$] $f_{m,j}''\in C_c^\infty(\Omega_{\wt w_{m,j}},\omega)$, for each $j$ with $n\le j\le m$,
\end{enumerate}
such that
\begin{equation}\label{eq: germ expansion at level mn}\CB_i(g,f_1)-\CB_i(g,f_2)=\sum_{j=m+1}^{l-m-1}\CB_i(g,f_{\ov w_j})+\sum_{j=n}^{m}\CB_i(g,f_{j,m}')+\sum_{j=n}^m \CB_i(g,f_{m,j}''),\end{equation}
for all $g\in \GL_{l}(F)$ and for all $i$ large enough depending only on $f_1,f_2$. If $n=m$,  \eqref{eq: germ expansion at level mn} is reduced to
\begin{equation}\label{eq: germ expansion at level mm}\CB_i(g,f_1)-\CB_i(g,f_2)=\sum_{j=m+1}^{l-m-1}\CB_i(g,f_{\ov w_j})+\CB_i(g,f_{m,m}').\end{equation}
See Remark \ref{remark: the case when m=n}.
 \end{indhypo}

We first prepare a lemma. For $a\in \GL_m(F),b\in \GL_n(F)$, we denote 
$$\bt_{m,n}(a,b)=\diag(a, I_{l-m-n},b)$$
as before.

\begin{lemma}\label{lemma: support 2} We fix the notations as in the Inductive Hypothesis \ref{indhypo2}.
\begin{enumerate}
\item For each $k$ with $1\le k\le l-1$, then for $i$ large enough which only depends on $f_1,f_2$, and for any $h\in P_{k,l-k}$, we have
$$\CB_i(h, f_{j,m}')=0,  \CB_i(h,f_{m,j}'')=0, \forall j, n\le j\le m.$$
\item For any $a\in \GL_m(F), b\in \GL_n(F), y\in \Mat_{m\times (l-m-1)}(F)$, we have 
\begin{align*}
\begin{array}{lll}
\CB_i(\wt w_{n,m}\bt_m(a,b)u_1(y),f_{\ov w_j})&=0, & m+1\le j\le l-m-1,\\
\CB_i(\wt w_{n,m}\bt_m(a,b)u_1(y),f_{j,m}')&=0, & n< j\le m,\\
\CB_i(\wt w_{n,m}\bt_m(a,b)u_1(y),f_{m,j}'')&=0, &n\le j\le m, \mathrm{~if~} n<m.
\end{array}
\end{align*}
In particular, by \eqref{eq: germ expansion at level mn}, we have $$\CB_i(\wt w_{n,m}\bt_m(a,b)u_1(y),f_1)-\CB_i(\wt w_{n,m}\bt_m(a,b)u_1(y),f_2)=\CB_i(\wt w_{n,m}\bt_m(a,b)u_1(y),f_{n,m}'). $$

\item If $u_1(y)\notin H_l^i$, we have 
$$\CB_i(\wt w_{n,m}\bt_m(a,b)u_1(y),f_{n,m}')=0 $$
for $i$ large enough depending only on $ f_1,f_2$.
\item For $k_1, k_2\in \Z$, the set $$\wpair{(a,b)\in [\GL_m]\times [\GL_n]| \CB_i(\tilde w_{n,m}\bt_{m,n}(a,b), f_{n,m}')\ne 0, |\det(a)|=q^{k_1}, |\det(b)|=q^{k_2}}$$ is compact. Here $[\GL_m]$ stands for $N_m(F)\bs \GL_m(F).$
\end{enumerate}
\end{lemma}
This is an analogue of \cite{Sp(2r)}*{Lemma 6.3}.

\begin{proof} (1) The proof is the same as the proof of Lemma \ref{lemma: support} (1) by noticing that $\wt w_{m,j}\notin P_{k,l-k}$ and $\wt w_{j,m}\notin P_{k,l-k}.$

(2) The proof is also similar to the proof of Lemma \ref{lemma: support} (2) and we give some details here. We consider the set 
\begin{align*}
S_{m,n}&=\wpair{w\in \bW(\GL_l): w=\wt w_{n,m}\bt_{m,n}(a,b), \mathrm{~for~some~}a\in \GL_m, b\in \GL_n}\\
&=\wpair{\wt w_{n,m}\bt_{m,n}(w,w'), \mathrm{~for~some~} w\in \bW(\GL_m), w'\in \bW(\GL_n)}.
\end{align*}
Note that the Weyl elements in $\GL_m$ (resp. $\GL_n$) form a Bruhat interval $[1,J_m]$ (resp. $[1,J_n]$). Thus for any $w\in S_{m,n}$ we have $\wt w_{n,m}\le w\le \wt w_{\max}$, where $$\wt w_{\max}=\wt w_{n,m}\bt_{m,n}(J_m, J_n)=\bpm && J_n\\ &I_{l-n-m} &\\ J_m &&\epm. $$
Notice that 
$$\wpair{\wt w_{n,m}\bt_{m,n}(a,b) u_1(y): a\in \GL_m(F), b\in \GL_n(F), y\in \Mat_{m\times (l-m-1)}}\subset \cup_{w\in S_{m,n}}C(w).$$
We have 
\begin{align*}
\theta_{\wt w_{\max}}&=\Delta-\wpair{\alpha_1,\dots,\alpha_m, \alpha_{l-n},\dots,\alpha_{l-1}},\\
\theta_{\ov w_j}&=\Delta-\wpair{\alpha_j}, \\
\theta_{\wt w_{j,m}}&=\Delta-\wpair{\alpha_m, \alpha_{l-j}}\\
\theta_{\wt w_{m,j}}&=\Delta-\wpair{\alpha_j, \alpha_{l-m}}.
\end{align*}
From these relations, we can see that $C(\wt w_{\max})\cap \Omega_{\ov w_j}=\emptyset$, for all $j$ with $m+1\le j\le l-m-1$; $C(\wt w_{\max})\cap \Omega_{ \wt w_{j,m}}=\emptyset$, for all $j$ with $ n<j\le m$; and $C(\wt w_{\max})\cap \Omega_{w_j''}=\emptyset,$ for all $j$ with $n\le j\le m$ except the case $n=j=m$. As in the proof of Lemma \ref{lemma: support} (2), this gives the conclusion. The ``in particular" part follows from the expansion \eqref{eq: germ expansion at level mn} and \eqref{eq: germ expansion at level mm} in the inductive hypothesis \ref{indhypo2}.

(3) This is an analogue of \cite{Sp(2r)}*{Lemma 6.3 (3)} and the proof is similar. We omit the details.

(4) This is an analogue of \cite{Sp(2r)}*{Lemma 6.3 (4)}. We also omit the details here.
\end{proof}
\begin{proposition}\label{proposition: main inductive step}
Assume that $1\le n\le m\le [l/2]$ and $m+n\le l-1$. Condition $\CC(m,n)$ implies that 
\begin{equation}\label{eq: equality on mn}
\CB_i(\wt w_{n,m}\bt_{m,n}(a,b),f_1)=\CB_i(\wt w_{n,m}\bt_{m,n}(a,b),f_2),\end{equation}
and
 \begin{equation}\label{eq: equality on nm}\CB_i(\wt w_{m,n}\bt_{n,m}(b,a),f_1)=\CB_i(\wt w_{m,n}\bt_{n,m}(b,a),f_2),\end{equation}
for all $a\in \GL_m(F), b\in \GL_n(F).$
\end{proposition}
\begin{proof}
Given any irreducible generic representation $\tau_1$ of $\GL_m(F)$ and $\tau_2$ of $\GL_n(F)$, the assumption says that 
$$\Gamma(\Bs,\pi_1\times (\tau_1,\tau_2),\psi)=\Gamma(\Bs,\pi_2\times (\tau_1,\tau_2),\psi).$$
We use the local functional equation of the form in \eqref{eq: modified LFE}. We first compute 
$$\Psi(\rho(\gamma_{m,n}^{-1})\CB_i^f,\xi_{\Bs}^{k,v_1,v_2};0)$$
for the section $\xi_{\Bs}^{k,v_1,v_2}$ as defined in Subsection \ref{subsection: construction of section} and $f=f_1,f_2$. Here $v_j\in \tau_j$ are arbitrary vectors and we take $k\ge i$ large enough. 
We have 
\begin{align*}
\Psi(\rho(\gamma_{m,n}^{-1})\CB_i^f,\xi_{\Bs}^{k,v_1,v_2};0)&=\int_{[\GL_{m+n}]}\int_{\ov U^{0,m,n}}\CB_i^f\left(\ov u \gamma_{m,n}\bpm h&\\ &I_{l-m-n} \epm \gamma_{m,n}^{-1} \right)\xi_\Bs^{k,v_1,v_2}(h)d\ov u d h.
\end{align*}
Here $[\GL_{m+n}]$ stands for $N_{m+n}(F)\bs \GL_{m+n}(F)$ and we will use similar notation below. Since $N_{m,n}M_{m,n}\ov N_{m,n}$ is dense in $\GL_{m+n}(F)$, the above integral over $N_{m+n}(F)\bs \GL_{m+n}(F)$ can be replaced by $N_{m+n}\bs N_{m,n}M_{m,n}\ov N_{m,n}=(N_m\bs \GL_m\times N_n\bs \GL_n)\ov N_{m,n}$, where an algebraic group is identified with its $F$-rational points. For $h=\diag(a,b)\ov u_{m,n}(y_2)\in (N_m\bs \GL_m\times N_n\bs \GL_n)\ov N_{m,n} $ with $y_2\in \Mat_{n\times m}$, we can take the Haar measure $dh=|\det(a)|^{-n}|\det(b)|^{m}d\ov v da db.$ A simple calculation on the conjugation by $\gamma_{m,n}$ shows that
\begin{align*}
\Psi(\rho(\gamma_{m,n}^{-1})\CB_i^f,\xi_{\Bs}^{k,v_1,v_2};0)=&\int_{[\GL_m]\times [\GL_n]}\int_{\ov N_{m,n}}\int_{\ov U^{0,m,n}}\CB_i^f\left(\bt_{m,n}(a,b) \bpm I_m &&&\\ &1&&\\ &&I_{l-m-n-1} &\\ y_2 &&y_1 &I_n \epm \right)\\
&\quad \xi_\Bs^{k,v_1,v_2}\left(\diag(a,b)\ov u_{m,n}(y_2)\right) |\det(a)|^{-n}|\det(b)|^{l-n-1}dy_2 dy_1 d a db.
\end{align*}
If $\ov u_{m,n}(y_2)\notin \ov N_{m,n}^k$, then $\xi_\Bs\left(\diag(a,b)\ov u_{m,n}(y_2)\right)=0 $ by the definition of $\xi_{\Bs}^{k,v_1,v_2}$, see \S\ref{subsection: construction of section}.  If $\ov u_{m,n}(y_2)\in \ov N_{m,n}^k$, then $\bpm I_m &&\\ &I_{l-m-n}&\\ y_2 &&I_n \epm\in \ov N_l\cap H_l^i$ because $k\ge i$. See the definition of $\ov N_{m,n}^k$ in \S \ref{subsection: construction of section}. By \eqref{eq: quasi-invariance of Howe vector}, we have 
\begin{align*}&\CB_i^f\left(\bt_{m,n}(a,b) \bpm I_m &&&\\ &1&&\\ &&I_{l-m-n-1} &\\ y_2 &&y_1 &I_n \epm \right)\\
&\qquad =\CB_i^f\left(\bt_{m,n}(a,b) \bpm I_m &&&\\ &1&&\\ &&I_{l-m-n-1} &\\ &&y_1 &I_n \epm \right) .
\end{align*}
Note that by the expansion \eqref{eq: germ expansion at level mn}, Lemma \ref{lemma: support} (1) and Lemma \ref{lemma: support 2} (1), we have
\begin{align*}
&\CB_i^{f_1}\left(\bpm a&&\\ &I_{l-m-n} &\\ &&b \epm \bpm I_m &&&\\ &1&&\\ &&I_{l-m-n-1} &\\ &&y_1 &I_n \epm \right)\\
&\qquad =\CB_i^{f_2}\left(\bpm a&&\\ &I_{l-m-n} &\\ &&b \epm \bpm I_m &&&\\ &1&&\\ &&I_{l-m-n-1} &\\ &&y_1 &I_n \epm \right).
\end{align*}
Thus we get
\begin{align*}
\Psi(\rho(\gamma_{m,n}^{-1})\CB_i^{f_1},\xi_{\Bs}^{k,v_1,v_2};0)=\Psi(\rho(\gamma_{m,n}^{-1})\CB_i^{f_2},\xi_{\Bs}^{k,v_1,v_2};0).
\end{align*}
Then by the local functional equation \eqref{eq: modified LFE} and the assumption on the local gamma factors, we have
$$\Psi(\rho(\gamma_{m,n}^{-1})\CB_i^{f_1},\wt\xi_{1-\wh \Bs}^{k,v_1,v_2};0)=\Psi(\rho(\gamma_{m,n}^{-1})\CB_i^{f_2},\wt\xi_{1-\wh\Bs}^{k,v_1,v_2};0),$$
or 
\begin{equation}\label{eq: right side local zeta integral is zero}
\Psi(\rho(\gamma_{m,n}^{-1})(\CB_i^{f_1}-\CB_i^{f_2}),\wt\xi_{1-\wh\Bs}^{k,v_1,v_2};0)=0.
\end{equation}
Here $\wt\xi_{1-\wh\Bs}^{k,v_1,v_2} $ denotes $M_{w_{m,n}}(\xi_{1-\wh\Bs}^{k,v_1,v_2})$ as usual. In the following, we write $\CB_i^{f_1}-\CB_i^{f_2} $ as $\CB_i$ for simplicity. We have 
\begin{align*}
\Psi(\rho(\gamma_{m,n}^{-1})\CB_i,\wt \xi_{1-\wh\Bs}^v;0)
&=\int_{[\GL_{m+n}]}\int_{\ov U^{0,n,m}}\CB_i\left(\ov u \gamma_{n,m}\bpm h &\\ & I_{l-m-n} \epm \gamma_{m,n}^{-1}\right)\wt \xi_{1-\wh \Bs}^{k,v_1,v_2}(h)d\ov u dh.
\end{align*}
Since $N_{n+m}\bs P_{n,m}w_{n,m}N_{m,n}\subset N_{n+m}\bs \GL_{n+m}$ is open and dense, we can replace the integral above over $N_{n+m}\bs \GL_{n+m} $ by $ N_{n+m}\bs P_{n,m}w_{n,m}N_{m,n}$. If $h=\diag(b,a)w_{n,m}u_{m,n}(x)\in N_{n+m}\bs P_{n,m}w_{n,m}N_{m,n}$ with $a\in \GL_m, b\in \GL_n, x\in \Mat_{m\times n}$, we can take the quotient measure to be $$dh=|\det(b)|^{-m}|\det(a)|^n dx dadb.$$ Thus we have 
\begin{align}
\begin{split}\label{eq: computation of the right side}
\Psi(\rho(\gamma_{m,n}^{-1})\CB_i,\wt \xi_{1-\wh\Bs}^v;0)=&\int_{[\GL_n]\times [\GL_m]}\int_{U^{0,n,m}}\\
 &~  \CB_i\left( \ov u \gamma_{n,m} \bpm &b &\\ a &&\\ &&I_{l-m-n} \epm \bpm I_m &x& \\ &I_n &\\ &&I_{l-m-n} \epm \gamma_{m,n}^{-1}\right)\\
 &~\wt \xi_{1-\wh \Bs}^{k,v_1,v_2}( \diag(b,a)w_{n,m}u_{m,n}(x))|\det(b)|^{-m}|\det(a)|^n d\ov u dx dadb.
 \end{split}
\end{align}
A matrix calculation shows that 
\begin{align*}
&\gamma_{n,m} \bpm &b&\\ a&&\\ &&I_{l-m-n} \epm \bpm I_m &x& \\ &I_n &\\ &&I_{l-m-n} \epm \gamma_{m,n}^{-1}\\
&\quad =\bpm &&b\\ &I_{l-m-n} &\\ a&&ax \epm\\
&\quad =\wt w_{n,m} \bt_{m,n}(a,b)u_1'(ax),
\end{align*}
where 
\begin{align*}
u_2'(ax):=\bpm I_m &&ax\\ &I_{l-m-n} & \\ &&I_n \epm.
\end{align*}
 On the other hand, for $\ov u\in \ov U^{0,n,m}$, we can write 
$$\ov u=\bpm I_{n+1} &&\\ &I_{l-m-n-1} &\\ &y&I_m \epm, \mathrm{~for~} y\in \Mat_{m\times (l-m-n-1)}.$$
We have
\begin{align*}
&\ov u \gamma_{n,m} \bpm &b &\\ a &&\\ &&I_{l-m-n} \epm \bpm I_m &x& \\ &I_n &\\ &&I_{l-m-n} \epm \gamma_{m,n}^{-1}\\
&=\wt w_{n,m} \bt_{m,n}(a,b)u_1((a^{-1}y,ax)),
\end{align*}
where recall that
$$u_1((a^{-1}y,ax))=\bpm I_m &&a^{-1}y &ax\\ &1&&\\ &&I_{l-m-n-1} &\\ &&&I_n \epm. $$
After changing variables on $x$ and $y$, \eqref{eq: computation of the right side} becomes
\begin{align*}
\Psi(\rho(\gamma_{m,n}^{-1})\CB_i,\wt \xi^{k,v_1,v_2}_{1-\wh\Bs};0)=&\int_{[\GL_m]\times [\GL_n]}\int_{y\in \Mat_{m\times (l-m-n-1)}}\int_{x\in \Mat_{m\times n}}\CB_i(\wt w_{n,m}\bt_{m,n}(a,b)u_1((y,x)))\\
&\quad \wt \xi_{1-\wh \Bs}^{k,v_1,v_2}( \diag(b,a)w_{n,m}u_{m,n}(x))|\det(b)|^{-m}|\det(a)|^{l-m-n-1} dy dx db da.
\end{align*}
Set 
$$D_i=\wpair{(y,x)\in \Mat_{m\times (l-m-n-1)}\times \Mat_{m\times n}: u_1((y,x))\in H_l^i\cap N_l},$$
as in Lemma \ref{lemma: support 2} (3). By Lemma \ref{lemma: support 2} (2) and (3), we have 
$$\CB_i(\wt w_{n,m}\bt_{m,n}(a,b)u_1((y,x)))=0, \mathrm{~if~} ((y,x))\notin D_i.$$
If $(y,x)\in D_i$, by \eqref{eq: quasi-invariance of Howe vector}, we have 
$$\CB_i(\wt w_{n,m}\bt_{m,n}(a,b)u_1((y,x)))= \CB_i(\wt w_{n,m}\bt_{m,n}(a,b)).$$
Moreover, by Subsection \ref{subsection: construction of section}, in particular, \eqref{eq: right side of the section}, for $k\ge k_0(D,v_1,v_2)$, we have 
$$\wt \xi_{1-\wh \Bs}^{k,v_1,v_2}(\diag(b,a)w_{n,m} u_{m,n}(x))=\vol(\ov N_{m,n}^k)|\det(b)|^{1-s_2+\frac{m-1}{2}}|\det(a)|^{-(1-s_1)-\frac{n-1}{2}} W_{v_1}(a)W_{v_2}(b).$$
Thus we get
\begin{align*}
\Psi(\rho(\gamma_{m,n}^{-1})\CB_i,\wt \xi^{k,v_1,v_2}_{1-\wh\Bs};0)=&\vol(D_i)\vol(\ov N_{m,n}^k)\int_{[\GL_m]\times [\GL_n]}\CB_i(\wt w_{n,m}\bt_{m,n}(a,b))\\
& W_{v_1}(a) W_{v_2}(b)|\det(b)|^{s_2^*}|\det(a)|^{s_1^*}  db da,
\end{align*}
where $s_2^*=1-s_2-\frac{m+1}{2}, s_1^*=-(1-s_1)-\frac{n-1}{2}+l-m-n-1$. The explicit form of $s_1^*, s_2^*$ is not important here.
By \eqref{eq: right side local zeta integral is zero}, we get 
\begin{align*}
\int_{[\GL_m]\times [\GL_n]}\CB_i(\wt w_{n,m}\bt_{m,n}(a,b))W_{v_1}(a) W_{v_2}(b)|\det(b)|^{s_2^*}|\det(a)|^{s_1^*}  db da=0,
\end{align*}
Note that the above formula holds for every $v_1\in \tau_1,v_2\in \tau_2$. Thus by Proposition \ref{proposition: Jacquet-Shalika} and Lemma \ref{lemma: support 2} (4), we get that
\begin{equation*}
\CB_i(\wt w_{n,m}\bt_{m,n}(a,b))=0, \forall a\in \GL_m(F), b\in \GL_n(F).
\end{equation*}
This proves the first equation \eqref{eq: equality on mn}. The second equation \eqref{eq: equality on nm} follows from the same proof by switching $m$ and $n$ and using the local gamma factor $\Gamma(\Bs,\pi\times (\tau_2,\tau_1),\psi)$ for an irreducible generic representation $\tau_1$ of $\GL_m(F)$ and $\tau_2$ of $\GL_n(F)$. This finishes the proof. 
\end{proof}
\begin{remark}
Assuming $\pi$ is unitarizable\footnote{This entails no loss of generality for proving Jacquet's local converse conjecture, see \cite{JNS}.}, by \cite{chai}*{Proposition 3.3}, we have 
\begin{equation}\label{eq: chai's lemma} \ov{\CB_i(g,f)}=\CB_i(g^*,f) , \end{equation}
for $f=f_1,f_2$. Here $g^*=J_l{}^t\! g^{-1}J_l$. The equation \eqref{eq: equality on nm} can be deduced from \eqref{eq: equality on mn} using \eqref{eq: chai's lemma} because $(\wt w_{n,m}\bt_{m,n}(a,b))^*=\wt w_{m,n}\bt_{n,m}(b^*,a^*) $. The identity \eqref{eq: chai's lemma} reveals a symmetry between $\CB_i(\wt w_{n,m}\bt_{m,n}(a,b),f)$ and $\CB_i(\wt w_{m,n}\bt_{n,m}(b,a),f)$. In our approach, this symmetry is mirrored in the corresponding definition of local gamma factors: $\CB_i(\wt w_{n,m}\bt_{m,n}(a,b),f) $ arises naturally in $\Gamma(\Bs,\pi\times (\tau_1,\tau_2),\psi)$, while $\CB_i(\wt w_{m,n}\bt_{n,m}(b,a),f)$ arises naturally in $\Gamma(\Bs,\pi\times (\tau_2,\tau_1),\psi)$, where $\tau_1$ (resp. $\tau_2$) is an irreducible generic representation of $\GL_m(F)$ (resp. $\GL_n(F)$).
\end{remark}

\begin{corollary}\label{cor: induction}
Suppose that $1\le n\le m\le [l/2]$ and $n+m\le l-1$. Then condition $\CC(m,n)$ implies that there exist functions 
\begin{enumerate}
\item[$\bullet$] $f_{\ov w_j}\in C_c^\infty(\Omega_{\ov w_j},\omega)$ for each $j$ with $m+1\le j\le l-m-1$;
\item[$\bullet$] $f_{j,m}'\in C_c^\infty(\Omega_{\wt w_{j,m}},\omega) $, for each $j$ with $n+1\le j\le m$; and
\item[$\bullet$] $f_{m,j}''\in C_c^\infty(\Omega_{\wt w_{m,j}},\omega)$, for each $j$ with $n+1\le j\le m$,
\end{enumerate}
such that
\begin{equation}\CB_i(g,f_1)-\CB_i(g,f_2)=\sum_{j=m+1}^{l-m-1}\CB_i(g,f_{\ov w_j})+\sum_{j=n+1}^{m}\CB_i(g,f_{j,m}')+\sum_{j=n+1}^m \CB_i(g,f_{m,j}''),\end{equation}
for all $g\in \GL_{l}(F)$ and for all $i$ large enough depending only on $f_1,f_2$.
\end{corollary}
\begin{proof}
The proof is similar to the proof of Corollary \ref{cor: base case} and is just a simple application of Theorem \ref{theorem: CST}. We give some details here. By Lemma \ref{lemma: support 2} (2) and Proposition \ref{proposition: main inductive step},  condition $\CC(m,n)$ gives 
\begin{equation} \label{eq: vanishing of fnm'}\CB_i(\wt w_{n,m} \bt_{m,n}(a,b), f_{n,m}' )=0,\end{equation}
for any $a\in \GL_m(F), b\in \GL_n(F)$. As in the proof of Lemma \ref{lemma: support 2} (2), we consider 
$$\wt w_{\max}=\wt w_{n,m}\bt_{m,n}(J_m, J_n)=\bpm &&J_m \\ &I_{l-m-n}&\\ J_n && \epm.$$
From the description of $\RB(\GL_l)$ in terms of subsets of $\Delta$, we can check that any $w\in \RB(\GL_l)$ with $\wt w_{n,m}\le w\le \wt w_{\max}$ has the form $\wt w_{n,m}\bt_{m,n}(w_1,w_2)$ for certain $w_1\in \bW(\GL_m), w_2\in \bW(\GL_n)$. Moreover, for any such $w$, we have $A_w\subset A_{\wt w_{\max}}$. From the definition \eqref{eq: Aw}, we see that any element $t\in A_w$ has the form 
$$z\bt_{m,n}(t_1,t_2),$$
with $z=zI_l$ in the center of $\GL_l(F)$, a diagonal element $t_1$ in $\GL_m$, and another diagonal element $t_2$ in $\GL_n$. Thus \eqref{eq: vanishing of fnm'} implies that 
\begin{equation}\label{eq: vanishing at wnm}\CB_i(w t, f_{n,m}')=0,\end{equation}
for all $w$ with $\wt w_{n,m}\le w\le w_{\max}$ and all $t\in A_w$. If we denote $w_{\max}'=\wt w_{m,n}\bt_{n,m}(J_n, J_m)$, then from \eqref{eq: equality on nm}, one can obtain that 
\begin{equation}\label{eq: vanishing at wmn}
\CB_i(w t, f_{n,m}')=0, \forall w\in [\wt w_{m,n}, w_{\max}'], t\in A_w.
\end{equation}
Similar as in the proof of Corollary \ref{cor: base case}, the result follows from Theorem \ref{theorem: CST}, \eqref{eq: vanishing at wnm} and \eqref{eq: vanishing at wmn}. Since this argument is almost identical to the proof of Corollary \ref{cor: base case}, we omit the details.
\end{proof}
For odd $l=2r+1$, the proof of Theorem \ref{theorem: second induction} (and hence the proofs of Theorems \ref{theorem: inductive} and \ref{theorem: main}) is complete. For even $l=2r$, by Corollary \ref{cor: induction}, condition $\CC(r,r-1)$ implies that
\begin{align}\label{eq: inductive condition when l=2r}
\CB_i(g,f_1)-\CB_i(g,f_2)=\CB_i(g,f_{r,r}'),
\end{align}
for some $f_{r,r}'\in C_c^\infty(\Omega_{\wt w_{r,r}},\omega)$. In \S\ref{subsection: even case}, we show that condition $\CC(r,r)$ forces $f'_{r,r}$ 	
 to vanish after suitably increasing $i$, thereby completing the proof of Theorems \ref{theorem: inductive} and \ref{theorem: main} in the even case.

\subsection{Conclude the proof when $l$ is even}\label{subsection: even case}
In this final subsection, we assume that $l=2r$ is even. Recall that for a character $\mu$ of $F^\times$, we have a Weil representation $\omega_{\psi^{-1},\mu,\mu^{-1}}$ of $\GL_{2r}(F)$, see \S\ref{subsection: review of gamma} or \cite{Morimoto}*{\S 2.2}.  For a positive integer $c$, we consider the function $\phi^c\in \CS(F^r\times F^r)$ defined by 
$$\phi^c(x,y)=\chi_{\fp^{(2r-1)c}}(x_1)\dots \chi_{\fp^{3c}}(x_{r-1})\chi_{1+\fp^c}(x_r)\chi_{\fp^{(2r-1)c}}(y_1)\dots \chi_{\fp^{3c}}(y_{r-1})\chi_{1+\fp^c}(y_r),$$
for $x=(x_1,x_2,\dots,x_r)\in F^r, y=(y_1,\dots,y_r)\in F^r$. Here, for a set $A\subset F$, $\chi_A$ denotes the characteristic function of $A$. 
\begin{proposition}\label{proposition: induction in the even case}
Condition $\CC(r,r)$ implies that 
$$\CB_i(w_{r,r}\bt_{r,r}(a,b),f_{r,r}')\omega_{\psi^{-1}}(w_{r,r})\phi^c(e_rb,e_ra^*)=0,$$
for any $a,b\in \GL_r(F)$, and for large $c>i$. Here $a^*=J_r{}^t\! a^{-1} J_r.$
\end{proposition}
\begin{proof}
The calculation below is similar to the case given in \cite{Sp(2r)}*{\S 7}. We give a sketch of the proof below. The corresponding local zeta integrals and local functional equations were recalled in  \S\ref{subsection: review of gamma}. Similarly as the calculation in Proposition \ref{proposition: main inductive step}, we have
\begin{align*}
\Psi(\CB_i^{f_1},\xi_{\Bs}^{k,v_1,v_2},\phi^c)=\Psi(\CB_i^{f_2},\xi_{\Bs}^{k,v_1,v_2},\phi^c).
\end{align*}
Thus, by the assumption on local gamma factors, we have
$$\Psi(\CB_i^{f_1},\wt\xi_{\Bs}^{k,v_1,v_2},\phi^i)=\Psi(\CB_i^{f_2},\wt\xi_{\Bs}^{k,v_1,v_2},\phi^i). $$
Again, we denote $\CB_i=\CB_i^{f_1}-\CB_i^{f_2}$ for simplicity, and we get $\Psi(\CB_i,\wt\xi_{1-\wh\Bs}^{k,v_1,v_2},\phi^c)=0 .$ On the other hand, by definition, we have 
\begin{align*}
\Psi(\CB_i,\wt\xi_{1-\wh\Bs}^{k,v_1,v_2},\phi^c)&=\int_{[\GL_{2r}]} \CB_i(g)\omega_{\psi^{-1}}(g)\phi^c(e_r,e_r)\wt \xi_{1-\Bs}^{k,v_1,v_2}(g)dg\\
&=\int_{[\GL_r]\times [\GL_r]}\int_{N_{r,r}}\CB_i(w_{r}\bt_{r}(a,b)u_{r}(x))\omega_{\psi^{-1}}(w_{r}\bt_{r}(a,b)u_{r}(x))\phi^c(e_r,e_r)\\
&\qquad \qquad \wt \xi_{1-\Bs}^{k,v_1,v_2}(w_r \bt_r(a,b) u_r(x))|\det(a)|^{r} |\det(b)|^{-r}dxdadb.
\end{align*}
Here, for simplicity, we write $\bt_{r,r}(a,b)=\diag(a,b)$ as $\bt_r (a,b)$, $w_{r,r}=\bpm &I_r\\ I_r & \epm$ as $w_r$ and $u_{r,r}(x)=\bpm I_r &x\\ &I_r \epm$ as $u_r(x)$. By Lemma \ref{lemma: support 2} (2) and (3), we have 
$$\CB_i(w_r \bt_r(a,b) u_r(x))=0, \mathrm{~if~} u_r(x)\notin N_{r,r}\cap H_{2r}^i.$$
If $u_r(x)\in N_{r,r}\cap H_{2r}^i$ and $k\gg 0$, by \eqref{eq: right side of the section}, we still have 
$$\wt \xi_{1-\wh\Bs}^{k,v_1,v_2}(w_r \bt_r(a,b) u_r(x))=\vol(\ov N_{r,r}^k)|\det(b)|^{1-s_2+\frac{r-1}{2}}|\det(a)|^{-(1-s_1)-\frac{r-1}{2}}W_{v_1}(a)W_{v_2}(b).$$
If $c>i$, from the Weil representation formula \cite{Morimoto}*{\S 2.2}, we can check that 
$$\omega_{\psi^{-1}}(u_r(x))\phi^c=\psi^{-1}(x)\phi^c, u_r(x)\in N_{r,r}\cap N_{2r}^i,$$
see \cite{Sp(2r)}*{Lemma 5.5} for a very similar calculation. Here $\psi$ is viewed as a character of the maximal unipotent subgroup $N_l$. Thus we get 
\begin{align*}
\omega_{\psi^{-1}}(w_{r}\bt_{r}(a,b)u_{r}(x))\phi^c(e_r,e_r)=\psi^{-1}(x)\mu(\det(ab))|\det(a)\det(b^{-1})|^{1/2}(\omega_{\psi^{-1}}(w_r)\phi^c)(e_rb,e_ra^*),
\end{align*}
 see \cite{Morimoto}*{\S 2.2} for the corresponding Weil representation formulas. On the other hand, for $u_r(x)\in N_{r,r}\cap H_{2r}^i$, by \eqref{eq: quasi-invariance of Howe vector}, we get that 
$$\CB_i(w_{r}\bt_{r}(a,b)u_{r}(x))=\psi(x)\CB_i(w_{r}\bt_{r}(a,b)). $$
Combining the above calculations, we get that 
$$\int_{[\GL_r]\times [\GL_r]}\CB_i(w_r\bt_r(a,b))\omega_{\psi^{-1}}(w_r)\phi^c(e_rb,e_ra^*)W_{v_1}(a)W_{v_2}(a)|\det(a)|^{s_1^*}|\det(b)|^{-s_2^*}dadb=0. $$
Here, $s_1^*$ and $-s_2^*$ are certain translations of $s_1,-s_2$ respectively, and $[\GL_r]$ stands for $N_r(F)\bs \GL_r(F)$. Now the result follows from Proposition \ref{proposition: Jacquet-Shalika}.
\end{proof}

 \begin{corollary}
Under condition $\CC(r,r)$, we have $\CB_i(g,f_1)=\CB_i(g,f_2)$ for all sufficiently large $i$ depending only on $f_1,f_2$.
\end{corollary}
\begin{proof}
The proof follows a similar argument as the proof of Corollary \ref{cor: induction}. Set 
$$w_{\max}=w_{r,r}\bt_{r,r}(J_r,J_r)=\bpm &J_r \\J_r &\epm,$$
which is indeed the longest Weyl element of $\GL_{2r}$. For an Weyl element $w\in \RB(\GL_{2r})$ such that $w_{r,r}\le w\le w_{\max}$, we can check that it has the form $w_{r,r}\bt_{r,r}(w_1,w_2)$ for some $w_1,w_2\in \bW(\GL_{r})$. We claim that $\CB_i(tw,f_{r,r}')=0$ for all $t\in T_{2l}(F)$ and all $w$ with $w_{r,r}\le w\le \wt w_{\max}$. We write  $t=\diag(a_1,\dots,a_{2r})\in T_{2r}(F)$. Since $\CB_i(~,f_{r,r}')$ has a central character, we can assume that $a_{r+1}=1$.

 From $w\ge w_{r,r}$, we have $\theta_w\subset \theta_{w_{r,r}}=\Delta-\wpair{\alpha_r}$. In particular, we have $\alpha_r\notin \theta_w$ and thus $\beta:=-w(\alpha_r)>0$. For a root $\gamma$, we fix an embedding $x_\gamma: F\to N_{2r}$ such that $\Im(x_\gamma)$ is the root space of $\beta$. Pick $y\in \fp^{(2\mathrm{ht}\beta+1)i}$, where $\mathrm{ht}(\beta)$ denotes the height of $\beta$. Then $x_{-\beta}(y)\in H_{2r}^i$, see \S\ref{subsection: Howe vector}. For, we have
$$twx_{-\beta}(y)=x_{\alpha_r}(\alpha(t)y)tw.$$
By \eqref{eq: quasi-invariance of Howe vector}, we get that $\CB_i(twx_{-\beta}(y),f_{r,r}')=\psi(\alpha_r(t)y) \CB_i(tw,f_{r,r}')$. Thus if $\CB_i(tw,f_{r,r}')\ne 0 $, we get that $\alpha_r(t)y\in \CO$ for any $y\in \fp^{(2\mathrm{ht}\beta+1)i}$, which implies that $a_r=\alpha_r(t)\in  \fp^{-(2\mathrm{ht}\beta+1)i}$.  If $\alpha_r(t)\in  \fp^{-(2\mathrm{ht}\beta+1)i}$, we write 
$$tw=tw_{r,r}\bt_{r,r}(w,w')=w_{r,r}\bt_{r,r}(t_1w,t_2w'),$$
for some $w,w'\in \bW(\GL_r).$ Here $t_2=\diag(a_1,\dots,a_r), t_1=\diag(a_{r+1},\dots,a_{2r})$. By Proposition \ref{proposition: induction in the even case}, we get that
\begin{equation}\label{eq: partial Bessel times Weil is zero}\CB_i(tw, f_{r,r}')\omega_{\psi^{-1}}(w_{r,r})\phi^c(e_r t_2 w', e_r t_1^*w^*)=0.\end{equation}
Write $v_1=e_r t_2 w'=[0,0,\dots,0,a_r]w'=[v_{11},\dots,v_{1r}]$, where only one $v_{1j}$ is nonzero, which is $a_r$. Moreover, we write $v_2=e_rt_1^*w^*=[0,\dots,0,1]w^*=[v_{21},\dots,v_{2r}]$, where only one entry $v_{2j}$ is nonzero, which is 1. From the Weil representation formula, we can take $c$ large enough such that $\omega_{\psi^{-1}}(w_{r,r})\phi^c(e_r t_2 w', e_r t_1^*w^*)\ne 0 $, see \cite{Sp(2r)}*{Lemma 5.5 (2)} for the detailed calculation in a similar situation. From \eqref{eq: partial Bessel times Weil is zero}, we get $\CB_i(tw,f_{r,r}')=0$ for any $t\in T_{2r}(F), w\in \RB(\GL_{2r})$ with $w_{r,r}\le w\le w_{\max}$. A direct application of Theorem \ref{theorem: CST} shows that $\CB_i(g,f_{r,r}')=0$ after increasing $i$ if necessary. This finishes the proof.
\end{proof}
This finishes the proof of Theorem \ref{theorem: second induction}, and thus Theorem \ref{theorem: inductive} and Theorem \ref{theorem: main}.

\begin{remark}
Suppose that $F$ is a finite field. Let $l,m,n$ be non-negative integers with $m+n<l$. Let $\pi$ be an irreducible supercuspidal representation of $\GL_l(F)$, $\tau_1,\tau_2$ be irreducible generic representations of $\GL_m(F)$ and $\GL_n(F)$ respectively. Then for $W\in \CW(\pi,\psi)$ and $f\in \Ind_{P_{m,n}(F)}^{\GL_{m+n}(F)}(\tau_1\boxtimes\tau_2)$, we can still define the local zeta integral $\Psi(W,f)$ and the local gamma factor $\Gamma(\pi\times (\tau_1,\tau_2),\psi)$ as in \S\ref{section: local theory}. As in the $p$-adic case, modulo a normalization factor, this gamma factor should be the product of the gamma factors $\gamma(\pi\times \tau_1,\psi)$ and $\gamma(\wt\pi\times \wt \tau_2,\psi)$, where these factors were developed in \cite{Ro} by imitating Jacquet--Piatetski-Shapiro--Shalika's theory \cite{JPSS}. The arguments in the preceding sections can be adapted to give a new proof of the finite field analogue of Jacquet’s local converse conjecture, originally established in \cite{Nien}. Analogous results for classical groups and the exceptional group $G_2$ over finite fields were proved in \cites{Liu-Zhang: classical, HL: finite, Liu-Zhang: G2}.
\end{remark}

\appendix
\section{Factorization of the gamma function}
\label{appendix}
In this appendix, we prove Proposition~\ref{proposition-gamma-comparison}. By Remark~\ref{remark-gamma-factor-JPSS}, it suffices to prove

\begin{proposition}
\label{prop-factorization}
Let $F$ be a non-archimedean local field of characteristic different from 2.
Let $\pi, \tau_1, \tau_2$ be irreducible generic representations of $\GL_{l}(F)$, $\GL_m(F)$, and $\GL_n(F)$ respectively, with $l>m+n$. Then we have
\begin{equation}
\label{eq-Thm3.1-main}
\Gamma(\Bs, \pi\times(\tau_1, \tau_2), \psi;j) = \frac{\Gamma((s_1+\frac{n}{2},0), \pi\times(\tau_1, 0), \psi;j)  \Gamma((0,s_2+\frac{m}{2}), \pi\times(0, \tau_2), \psi;j)}{\gamma(s_1+s_2,\tau_1\times \wt \tau_2,\psi)}.
\end{equation}
\end{proposition}

The rest of this appendix is devoted to proving Proposition~\ref{prop-factorization}. We start by discussing the multiplicativity of the normalized intertwining operator $M^*(\Bs, \btau)$. Let $\tau_{1,j}$ (resp. $
\tau_{2,j})$ be representations of $\GL_{m_j}$ (resp. $\GL_{n_j}$) for $j=1, 2$, with $m_1+m_2=m$ (resp. $n_1+n_2=n$), and we first assume that
\begin{equation*}
	\tau_1=\Ind_{P_{m_1, m_2}}^{\GL_m}(\tau_{1,1}\boxtimes \tau_{1,2}), \quad \tau_2=\Ind_{P_{n_1, n_2}}^{\GL_n}(\tau_{2,1}\boxtimes \tau_{2,2}).
\end{equation*}
For a partition $l=r_1+r_2+\cdots+r_k$, let $P_{r_1, r_2, \cdots, r_k}$ be the standard parabolic subgroup of $\GL_{l}$ whose Levi part $M_{P_{r_1, r_2, \cdots, r_k}}$ is isomorphic to $\GL_{r_1}\times\GL_{r_2}\times\cdots\times\GL_{r_k}$, and denote its unipotent radical by $N_{r_1, r_2, \cdots, r_k}$.  

We consider the induced representation 
\begin{equation*} \begin{split}\Ind_{P_{m_1, m_2+n_1, n_2}}^{\GL_{m+n}} \left( \mathcal{W}(\tau_{1,1},\psi^{-1})||^{s_1-\frac{1}{2}} \otimes \Ind_{P_{m_2, n_1}}^{\GL_{m_2+n_1}} ( \mathcal{W}(\tau_{1,2},\psi^{-1})||^{s_1-\frac{1}{2}}  \otimes \mathcal{W}(\tau_{2,1},\psi^{-1})||^{-s_2+\frac{1}{2}}    )    \right. \\\left. \otimes \mathcal{W}(\tau_{2,2},\psi^{-1})||^{-s_2+\frac{1}{2}}   \right)\end{split}\end{equation*}
and we denote its underlying space by $\V( (\tau_{1,1},\tau_{1,2}\otimes \tau_{2,1}, \tau_{2,2}), (s_1, s_1, s_2, s_2))$. The underlying space $\V( (\tau_{1,1},\tau_{1,2}\otimes \tau_{2,1}, \tau_{2,2}), (s_1, s_1, s_2, s_2))$ consists of smooth functions 
$$\varphi_\Bs: \GL_{m+n}\times \GL_{m_2+n_1}\times M_{P_{m_1, m_2, n_1, n_2}}\to \mathbb{C}$$
such that for $h\in \GL_{m+n}$, $h_1\in \GL_{m_2+n_1}, m\in M_{P_{m_1, m_2, n_1, n_2}}$, $h_0\in \GL_{m_2+n_1}$, $a_i\in \GL_{m_i}$, and $b_i\in \GL_{n_i}$, we have
\begin{itemize}
\item 
\begin{equation*}
\begin{split}
\varphi_\Bs &\left( \bpm a_1 &* &*\\ &h_0 & *\\ && b_2 \epm h, h_1, m \right) =  \delta_{P_{m_1, m_2+n_1, n_2}}^{\frac{1}{2}} \left( \bpm a_1 & &\\ &h_0 & \\ && a_2 \epm \right)  \\
    &\quad |\det(a_1)|^{s_1-\frac{1}{2}} |\det(b_2)|^{-s_2+\frac{1}{2}} \varphi_\Bs \left(   h, h_1 h_0, m \bpm a_1 & &\\ &I_{m_2+n_1} & \\ && a_2 \epm  \right), \end{split}\end{equation*}
\item 
\begin{equation*}
\begin{split}
	\varphi_\Bs \left(  h, \bpm a_2 &* \\ &b_1\epm h_1, m \right) =& \delta^{\frac{1}{2}}_{P_{m_2, n_1}}\left( \bpm a_2 & \\ &b_1\epm \right) |\det(a_2)|^{s_1-\frac{1}{2}} |\det(b_1)|^{-s_2+\frac{1}{2}} \\
	& \varphi_\Bs \left(  h, h_1, m \bpm I_{m_1} && &\\ &a_2 &&\\ &&b_1&\\ &&&I_{n_2} \epm \right),
\end{split}
\end{equation*}

\item for fixed $(h,h_1)\in \GL_{m+n}\times \GL_{m_2+n_1}$, the function $m\mapsto \varphi_\Bs(h, h_1, m)$ belongs to the space $\mathcal{W}(\tau_{1,1},\psi^{-1})\otimes  \mathcal{W}(\tau_{1,2},\psi^{-1})  \otimes \mathcal{W}(\tau_{2,1},\psi^{-1}) \otimes \mathcal{W}(\tau_{2,2},\psi^{-1}).$ 
\end{itemize}
For any $\varphi_\Bs\in \V( (\tau_{1,1},\tau_{1,2}\otimes \tau_{2,1}, \tau_{2,2}), (s_1, s_1, s_2, s_2))$, 
we define an element $$\xi_{\varphi_{\Bs}} \in \Ind_{P_{m,n}}^{\GL_{m+n}}(\mathcal{W}(\tau_1,\psi^{-1}) ||^{s_1-\frac{1}{2}}\otimes \mathcal{W}(\tau_2,\psi^{-1})||^{-s_2+\frac{1}{2}})$$ by a Jacquet integral
\begin{equation}
\label{eq-GLl-Jacquet-integral}
\begin{split}
\xi_{\varphi_{\Bs}}(g, a, b)= & |\det(a)|^{-s_1-\frac{n-1}{2}} |\det(b)|^{s_2+\frac{m-1}{2}}  \\&\cdot \int_{Z_{m_2,m_1,n_2,n_1} } \varphi_{\Bs}\left( \omega   z \diag(a, b)g, I_{m_2+n_1}, I_{M_{P_{m_1, m_2, n_1, n_2}}}  \right) \psi(z)dz
\end{split}
\end{equation}
where $g\in \GL_{m+n}$, $a\in \GL_m$, $b\in \GL_n$, and 
\begin{equation*}
Z_{m_2,m_1,n_2,n_1}=\left\{ \begin{pmatrix} I_{m_2}& z_1&&\\ &I_{m_1} &&\\ &&I_{n_2} & z_2 \\ &&&I_{n_1}\end{pmatrix} \right\}, \quad \omega=\begin{pmatrix} & I_{m_1} &&\\ I_{m_2} & &&\\ &&&I_{n_1} \\ &&I_{n_2} & \end{pmatrix}.
\end{equation*}
This integral may not converge absolutely. To rectify this, we may twist the representations $\tau_{i,j}$ by auxiliary complex parameters $\zeta_{i,j}$ for $i, j=1, 2$, and there is a cone where the integral \eqref{eq-GLl-Jacquet-integral} is absolutely convergent. In the following, we will skip the parameter $\zeta_{i,j}$.

For $\theta_{\Bs}=\theta_{s_1, s_2}\in   \Ind_{P_{m_2,n_1}}^{\GL_{m_2+n_1}}  (\CW(\tau_{1,2},\psi^{-1})||^{s_1-\frac{1}{2}} \otimes    \CW(\tau_{2,1},\psi^{-1})||^{-s_2+\frac{1}{2}})$, recall that we have the standard intertwining operator $M(\Bs, (\tau_{1,2},\tau_{2,1}))$ defined by 
\begin{equation*}
\begin{split}
\left( M(\Bs, (\tau_{1,2},\tau_{2,1}))\theta_{\Bs} \right) (g, a_2, a_1)  &= \int_{N_{n_1, m_2}}\theta_{\Bs}(w_{m_2, n_1}ug,a_1,a_2)du
\end{split}
\end{equation*}
where $w_{m_2,n_1}= \begin{pmatrix} &I_{m_2} \\ I_{n_1} & \end{pmatrix}$. Then $$M(\Bs, (\tau_{1,2},\tau_{2,1}))\theta_{\Bs}\in   \Ind_{P_{n_1,m_2}}^{\GL_{m_2+n_1}}   (\CW(\tau_{2,1},\psi^{-1}) ||^{-s_2+\frac{1}{2} } \otimes   \CW(\tau_{1,2},\psi^{-1}) ||^{ s_1-\frac{1}{2}}).$$
The normalized intertwining operator is $M^*(\Bs, (\tau_{1,2},\tau_{2}))=\gamma(s_1+s_2, \tau_{1,2}\times\wt \tau_2,\psi)M(\Bs, (\tau_{1,2},\tau_{2}))$, and it is defined to satisfy the functional equation
\begin{align}
\begin{split}
\int_{N_{n_1,m_2}}&\theta_{\Bs}(w_{m_2,n_1}ug,I_{m_2}, I_{n_1}) \psi^{-1}(u) du \\
&= \int_{N_{m_2,n_1}} (M^*(\Bs, (\tau_{1,2},\tau_{2,1}))\theta_{\Bs}) (w_{n_1,m_2}ug,I_{n_1}, I_{m_2}) \psi^{-1}(u) du.
\end{split}
\label{eq-normalized-intertwining-GL-eq1}
\end{align}
In the above integral, the measure is taken to be the product measure of self dual Haar measure with respect to $\psi$.

Similarly, we have intertwining operators $M^*(\Bs, (\tau_{1,1}, \tau_{2,1}))$, $M^*(\Bs, (\tau_{1,2}, \tau_{2,2}))$, and $M^*(\Bs, (\tau_{1,1}, \tau_{2,2}))$.
Given any $\varphi_{\Bs} \in \V_{P_{m_1, m_2+n_2, n_1}}^{\GL_{m+n}} ((\tau_{1,1}, \tau_{1,2}),(\tau_{2,1}, \tau_{2,2}), (s_1, s_1), (-s_2, -s_2))$, we can apply the intertwining operators repeatedly to get the following sections:
\begin{equation*}
\begin{split}
\varphi_{\Bs}^{\prime} := M^*(\Bs, (\tau_{1,2}, \tau_{2,1}) )\varphi_{\Bs}  \in\V_{P_{m_1, m_2+n_2, n_1}}^{\GL_{m+n}} ((\tau_{1,1}, \tau_{2,1}),(\tau_{1,2}, \tau_{2,2}), (s_1, 1-s_2), (1-s_1, s_2)),
\end{split}
\end{equation*}
\begin{equation*}
\begin{split}
\varphi_{\Bs}^{\prime\prime} := M^*(\Bs, (\tau_{1,2}, \tau_{2,2})) M^*(\Bs, (\tau_{1,1},\tau_{2,1}))  \varphi_{\Bs}^\prime  \in  \V_{P_{n_1,m_1+n_2,m_2}}^{\GL_{m+n}} ((\tau_{2,1}, \tau_{1,1}),(\tau_{2,2}, \tau_{1,2}), (1-s_2, s_1), (s_2, 1-s_1)),
\end{split}
\end{equation*}
\begin{equation*}
\begin{split}
\varphi_{\Bs}^{\prime\prime\prime} := M^*(\Bs,(\tau_{1,1},\tau_{2,2}) )\varphi_{\Bs}^{\prime\prime}  \in   \V_{P_{n_1, n_2+m_1, m_2}}^{\GL_{m+n}} ((\tau_{2,1}, \tau_{2,2}),(\tau_{1,1}, \tau_{1,2}), (1-s_2, 1-s_2), (1-s_1, 1-s_1)).
\end{split}
\end{equation*}
By the multiplicativity of intertwining operators \cite{Shahidi-on-certain-L-functions}, we have
\begin{equation}
M^*(\Bs, (\tau_1, \tau_2)) \xi_{\varphi_{\Bs}}=\xi_{\varphi_{\Bs}^{\prime\prime\prime} }.
\label{eq-GL-multiplicativity-intertwining-operators}
\end{equation}
Note that \eqref{eq-GL-multiplicativity-intertwining-operators} still makes sense when $m_1=n_1=0$ (so that $m_2=m, n_2=n$).

\begin{proof}[Proof of Proposition~\ref{prop-factorization}]
Recall that $\Psi(W, \xi_{\Bs};j)$ is given by 
\begin{equation*}
\Psi(W, \xi_{\Bs}; j) =\int_{N_{m+n}\bs \GL_{m+n}}\int_{\ov U^{j,m,n}}W\left(\ov u \gamma_{m,n}\bpm h&\\ &I_{l-m-n} \epm \right)\xi_{\Bs}(h)d\ov udh.
\end{equation*}
We factor the $dh$ integration over $\GL_m\times \GL_n=\{ \diag(h_1, h_2): h_1\in \GL_m, h_2\in \GL_n\}$, to obtain
\begin{equation*}
\begin{split}
\int_{(\GL_m\times \GL_n ) N_{m+n}\bs \GL_{m+n}} \int_{N_n\backslash \GL_n}  \int_{N_m\backslash \GL_m}  \int_{\ov U^{j,m,n}}W\left(\ov u \gamma_{m,n} \bpm h_1 &&\\&h_2&\\ & &I_{l-m-n} \epm \bpm h&\\ &I_{l-m-n} \epm \right)  \\
|\det(h_1)|^{(s_1+\frac{n}{2})-\frac{1}{2}} |\det(h_2)|^{-(s_2+\frac{m}{2})+\frac{1}{2}} \xi_{\Bs}(h, h_1, h_2)d\ov u dh_1 dh_2 dh.
\end{split}
\end{equation*}
Note that $\gamma_{m,0}=I_{l}$, $\gamma_{0,n}=\bpm &I_{l-n} \\ I_n &\epm$, and
\begin{equation*}
\gamma_{m,n} \bpm h_1 &&\\&h_2&\\ & &I_{l-m-n} \epm \gamma_{m,n}^{-1}= \bpm h_1 &&\\&I_{l-m-n} &\\ & &h_2 \epm.
\end{equation*}
We decompose $\ov U^{j,m,n}$ into $\ov U^{j,m,0} \ov U^{j,0,n}$. Then $\Psi(W, \xi_{\Bs};j)$ is equal to 
\begin{equation*}
\begin{split}
\int_{(\GL_m\times \GL_n ) N_{m+n}\bs \GL_{m+n}} & \int_{N_n\backslash \GL_n}     \int_{\ov U^{j,0,n}}\int_{N_m\backslash \GL_m}  \int_{\ov U^{j,m,0}} \\
&W\left(\ov u_1   \gamma_{m,0}  \bpm h_1 &\\&I_{l-m} \epm  \ov u_2 \gamma_{0,n}  
\bpm h_2&\\ &I_{l-n} \epm  \gamma_{0,n}^{-1} \gamma_{m,n}   \bpm h &\\&I_{l-m-n} \epm \right)  \\
& \ \ \ \ \ \ \ \ \ \ \ \ |\det(h_1)|^{(s_1+\frac{n}{2})-\frac{1}{2}} |\det(h_2)|^{-(s_2+\frac{m}{2})+\frac{1}{2}} \xi_{\Bs}(h, h_1, h_2)d\ov u_1 dh_1 d\ov u_2 dh_2  dh.
\end{split}
\end{equation*}
Now we apply the intertwining operator $M^*((s_1, 0), (\tau_1, 0))$ to the inner $d\ov u_1 dh_1$ integral to obtain
\begin{equation}
\Gamma((s_1+\frac{n}{2}, 0), \pi\times (\tau_1, 0), \psi;j) \Psi(W, \xi_{\Bs}; j)=\Psi(W, M^*((s_1, 0), (\tau_1, 0))\xi_{\Bs}; j).
\label{eq-GLl-eq1}
\end{equation}
By the same argument as in the proof of \cite{Kaplan2015}*{(7.1), Section 8} (see also \cite{Kaplan2013}*{Section 4} and \cite{Morimoto}*{Section 4}), one can conclude the proof of Proposition~\ref{prop-factorization} as follows. We denote the following sections:
 $\xi_{\Bs}^\prime=M^*((s_1, 0), (\tau_1, 0))\xi_{\Bs}$, $\xi_{\Bs}^{\prime\prime}=M^*((s_1, s_2), (\tau_1, \tau_2)) M^*((0,0),(0,0)) \xi_{\Bs}^\prime$, and $\xi_{\Bs}^{\prime\prime\prime}=M^*((0,s_2), (0,\tau_2)) \xi_{\Bs}^{\prime\prime}$. By the similar functional equation \eqref{eq-GL-multiplicativity-intertwining-operators} for $M^*((0,0),(0,0))$ and $M^*((s_1, s_2), (\tau_1, \tau_2))$, one has $\Psi(W,  \xi_{\Bs}^{\prime\prime}; j)=\Psi(W,  \xi_{\Bs}^{\prime}; j).$
The proof of \eqref{eq-GLl-eq1} immediately gives
\begin{equation}
\Gamma\left((0,s_2+\frac{m}{2}),\pi\times(0, \tau_{2}),\psi;j \right) \Psi(W, \xi_{\Bs}^{\prime\prime}; j) =  \Psi(W,  \xi_{\Bs}^{\prime\prime \prime}; j).
\label{eq-GL-multiplicativity-reduction-2}
\end{equation}
Since $M^*(\Bs, (\tau_1, \tau_2)) \xi_{\Bs}=\xi_{\Bs}^{\prime\prime\prime}$, we have
\begin{equation*}
\begin{split}
	 &\Gamma\left((0,s_2+\frac{m}{2}),\pi\times(0, \tau_{2}),\psi;j\right) \Gamma\left((s_1+\frac{n}{2}, 0),\pi\times(\tau_{1}, 0),\psi;j\right) \Psi(W, \xi_{\Bs}; j) \\
	  =& \Gamma\left((0,s_2+\frac{m}{2}),\pi\times(0, \tau_{2}),\psi;j\right) \Psi(W, \xi_{\Bs}^\prime; j) \\
	  =& \Gamma\left((0,s_2+\frac{m}{2}),\pi\times(0, \tau_{2}),\psi;j\right) \Psi(W, \xi_{\Bs}^{\prime\prime}; j) \\
	  =&  \Psi(W, \xi_{\Bs}^{\prime\prime\prime}; j)	 \\
	  =& \Psi(W, M^*(\Bs, (\tau_1, \tau_2)) \xi_{\Bs}; j)\\
	  =& \gamma(s_1+s_2, \tau_1\times\wt \tau_2,\psi) \Psi(W, M(\Bs, (\tau_1, \tau_2)) \xi_{\Bs}; j). 
\end{split}
\end{equation*}
Thus we obtain 
$$
\Gamma(\Bs, \pi\times(\tau_1, \tau_2), \psi;j) = \frac{\Gamma((s_1+\frac{n}{2},0), \pi\times(\tau_1, 0), \psi;j)  \Gamma((0,s_2+\frac{m}{2}), \pi\times(0, \tau_2), \psi;j)}{\gamma(s_1+s_2, \tau_1\times\wt \tau_2,\psi)}.
$$
This completes the proof of Proposition~\ref{prop-factorization}.
\end{proof}

\end{document}